\documentclass[a4paper, 11pt]{article}
\usepackage{header2}
\usepackage{permutation}

\usepackage[font={small}]{caption}
\begin{document}

\title
{\makebox[0pt][c]{A Combinatorial Approach to Rauzy-type Dynamics III:}\\
\makebox[0pt][c]{The Sliding Dynamics, Diameter and Algorithm}}

\author{Quentin De Mourgues\\
%$^{\dagger, \ddagger}$\\
  \multicolumn{1}{p{.7\textwidth}}{\centering\emph{\small 
\rule{0pt}{14pt}%
LIPN, Universit\'e Paris 13\\
99, av.~J.-B.~Cl\'ement, 93430 Villetaneuse, France\\
{\tt quentin.demourgues@lipn.fr}}}
}

\maketitle

% !!!! latex abstract is buggy and breaks a lot of formatting in
% following text !!!!
% http://latex-community.org/forum/viewtopic.php?t=3996

\begin{center}
\begin{minipage}{.9\textwidth}
\noindent
\small
{\bf Abstract.} Rauzy-type dynamics are group actions on a collection of combinatorial
objects. The first and best known example (the Rauzy dynamics) concerns an action on
permutations, associated to interval exchange transformations (IET)
for the Poincar\'e map on compact orientable translation surfaces. The
equivalence classes on the objects induced by the group action 
have been classified by Kontsevich and Zorich in \cite{KZ03} and correspond bijectively to the connected components of the strata of the moduli space of abelian differentials.

\qquad In a paper \cite{BoiCM} Boissy proposed a Rauzy-type dynamics that acts on a subset of the permutations (the standard permutations) and conjectured that the Rauzy classes of this dynamics are exactly the Rauzy classes of the Rauzy dynamics restricted to standard permutations. In this paper, we apply the labelling method  introduced in \cite{D18} to classify this dynamics thus proving Boissy's conjecture.

\qquad Finally, this paper conclude our serie of three papers on the study of the Rauzy dynamics by presenting two new results on the Rauzy classes: An quadratic algorithm for outputting a path between two connected permutations and a tight $\Theta(n)$ bound on the diameter of the Rauzy classes for the alternating distance.
\end{minipage}

\end{center}

\newpage

\tableofcontents

\newpage

%%%%%%%%%%%%%%%%%%%%%%%%%%%%%%%%%%%%%%%%%%%%%%%%%%%%%%%

%\input{pp02_permutdiags.tex}
\section*{First part}
\label{sec:intros}
\addcontentsline{toc}{section}{\nameref{sec:intros}}

% \newpage
% \newpage

%%%%%%%%%%%%%%%%%%%%%%%%%%%%%%%%%%%%%%%%%%%%%%%%%%%%%%%
% \part{Introduction}

%%%%%%%%%%%%%%%%%%%%%%%%%%%%%%%%%%%%%%%%%%%%%%%%%%%%%%%

\section{Definition of the extended Rauzy classes}
%We will now apply the labelling method defined in the previous section to give an original proof of the classification of the Rauzy dynamics. Other proofs have been achieved in \cite{Boi12} (Boissy uses the classification proof for the extended Rauzy dynamics appearing in \cite{KZ03}) using geometric methods and in \cite{Fic16} and \cite{DS17} using combinatorial methods.
%
%The extended Rauzy classes (classified in \cite{KZ03}) are of interest in the translation surface field as they are in one-to-one correspondance with the connected components of the strata of the moduli space of abelian differentials (see \cite{Vee82}). As for the non-extended Rauzy classes, they are in one-to-one correspondance with the connected components of the strata of the moduli space of abelian differentials with a marked zero as shown in \cite{Boi12}.

The extended Rauzy classes are the equivalence classes of the extended Rauzy dynamics on the set of irreducible permutations. Those classes can be characterised by three invariants: the cycle invariant, the sign invariant and the hyperelliptic class. The characterisation theorem was proven by Kontsevich and Zorich in \cite{KZ03} and a formulation of it is theorem \ref{thm.Main_theorem_2}. The study of the Rauzy dynamics was initiated by \cite{Rau79} and \cite{Vee82} in the context of interval exchange transformations and the moduli space of abelian differentials.

 In the next two subsections we define the extended Rauzy dynamics and the three invariants.

As outlined in the abstract, the first goal of this article is to classify a dynamics introduced by Boissy in \cite{BoiCM} that we will call the sliding dynamics. Boissy conjectured that the Rauzy classes of this dynamics are the restriction of the extended Rauzy classes (of the extended Rauzy dynamics) to the subset of standard permutations. We introduce the dynamics in subsection \ref{sec_sliding_dyn}

%-------------------------------------------------------
\subsection{The extended Rauzy dynamics}
% : dynamics and extended dynamics on
%   permutations, and dynamics on matchings}
\label{ssec.3families}
%-------------------------------------------------------

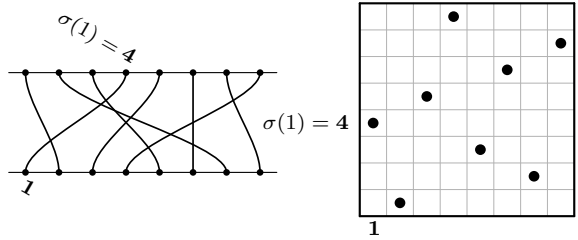
\begin{figure}[bt!]
\[
\begin{array}{cc}
\s=[41583627] \in \kS_8
&
\s=[41583627] \in \kS_8
\\
\textrm{diagram representation}
&
\textrm{matrix representation}
\\
\raisebox{12.5pt}{%
\setlength{\unitlength}{10pt}
\begin{picture}(11,3)
\put(0,0){\includegraphics[scale=2.5]{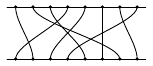}}
\put(.9,-.15){\rotatebox{-30}{\scriptsize{$\bm{1}$}}}
\put(3.6,5.4){\rotatebox{-30}{\makebox[0pt][c]{\scriptsize{$\s(1)=\bm{4}$}}}}
\end{picture}
}
&
\rule{0pt}{90pt}\raisebox{0pt}{%
\setlength{\unitlength}{10pt}
\begin{picture}(8,8)
\put(0,0){\includegraphics[scale=1]{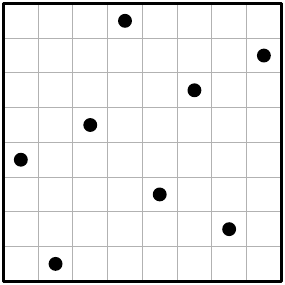}}
\put(0.4,-.6){\scriptsize{$\bm{1}$}}
\put(-.3,3.4){\makebox[0pt][r]{\scriptsize{$\s(1)=\bm{4}$}}}
\end{picture}
}
\end{array}
\]
\caption{\label{fig_match_representation} Diagram and matrix epresentations of permutations.}
\end{figure}

Let $\kS_n$ denote the set of permutations of size $n$.  Let us call
$\coxe$ the permutation $\coxe(i)=n+1-i$.

We say that $\s \in \kS_n$ is \emph{irreducible} if $\coxe \s$ doesn't
leave stable any interval $\{1,\ldots,k\}$, for $1 \leq k < n$,
i.e.\ if $\{\s(1), \ldots, \s(k)\} \neq \{n-k+1,\ldots,n\}$ for any
$k=1,\ldots,n-1$. Let us call $\kSirr_n$ 
the corresponding sets of irreducible configurations.

Finally a standard permutation $\s \in \kSirr_n$ is a permutation verifying $\s(1)=1$. We let $\kS^{st}_n$ be the set of standard permutations.

We represent permutation as
arcs in a horizontal strip, connecting $n$ points at the bottom
boundary to $n$ points on the top boundary (as in
Figure~\ref{fig_match_representation}, left).  Both sets of
points are indicised from left to right. We use the name of
\emph{diagram representation} for such representations.

We will also often represent permutations as grids filled with one
bullet per row and per column (and call this \emph{matrix
  representation} of a permutation). We choose here to conform to the
customary notation in the field of Permutation Patterns, by adopting
the algebraically weird notation, of putting a bullet at the
\emph{Cartesian} coordinate $(i,j)$ if $\s(i)=j$, so that the identity
is a grid filled with bullets on the \emph{anti-diagonal}, instead
that on the diagonal. An example is given in figure
\ref{fig_match_representation}, right.

Let us define a special set of permutations (in cycle notation)
\begin{subequations}
\label{eqs.opecycdef}
\begin{align}
\gamma_{L,n}(i)
&=
(i-1\;i-2\;\cdots\;1)(i)(i+1)\cdots(n)
\ef;
& 
\gamma_{L',n}(i)
&=\gamma_{L,n}(i)^{-1}
\ef;
\\
\gamma_{R,n}(i)
&=
(1)(2)\cdots(i)(i+1\;i+2\;\cdots\;n)
\ef;
&
\gamma_{R',n}(i)
&=\gamma_{R,n}(i)^{-1}
\ef;
\end{align}
\end{subequations}
i.e., in a picture
\begin{align*}
\gamma_{L,n}(i):
&
\quad
\setlength{\unitlength}{8.75pt}
\begin{picture}(8,2)(-.1,.8)
\put(-.3,0.3){\includegraphics[scale=1.75]{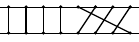}}
\put(0,0){$\scriptstyle{1}$}
\put(3,0){$\scriptstyle{i}$}
\put(7,0){$\scriptstyle{n}$}
\end{picture}
&
\gamma_{L',n}(i):
&
\quad
\setlength{\unitlength}{8.75pt}
\begin{picture}(8,2)(-.1,.8)
\put(-.3,0.3){\includegraphics[scale=1.75]{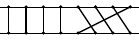}}
\put(0,0){$\scriptstyle{1}$}
\put(3,0){$\scriptstyle{i}$}
\put(7,0){$\scriptstyle{n}$}
\end{picture}\\
\gamma_{R,n}(i):
&
\quad
\setlength{\unitlength}{8.75pt}
\begin{picture}(8,2)(-.1,.8)
\put(-.3,0.3){\includegraphics[scale=1.75]{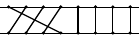}}
\put(0,0){$\scriptstyle{1}$}
\put(4,0){$\scriptstyle{i}$}
\put(7,0){$\scriptstyle{n}$}
\end{picture}
&
\gamma_{R',n}(i):
&
\quad
\setlength{\unitlength}{8.75pt}
\begin{picture}(8,2)(-.1,.8)
\put(-.3,0.3){\includegraphics[scale=1.75]{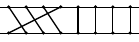}}
\put(0,0){$\scriptstyle{1}$}
\put(4,0){$\scriptstyle{i}$}
\put(7,0){$\scriptstyle{n}$}
\end{picture}
\end{align*}

\begin{figure}[tb!]
\begin{align*}
L\;\Big(
\;
\raisebox{14.2pt}{\includegraphics[scale=1.75, angle=180]{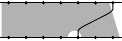}}
% \raisebox{-9pt}{\includegraphics[scale=1.75]{FigFol/Figure1_fig_PPlr3.pdf}}
\;
\Big)
&=
\raisebox{27.2pt}{\includegraphics[scale=1.75, angle=180]{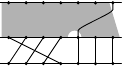}}
% \raisebox{-9pt}{\includegraphics[scale=1.75]{FigFol/Figure1_fig_PPlr4.pdf}}
\rule{0pt}{38pt}
&
L'\;\Big(
\;
\raisebox{-9pt}{\reflectbox{\includegraphics[scale=1.75]{FigFol/Figure1_fig_PPlr1_simp.pdf}}}
\;
\Big)
&=
\raisebox{-24pt}{\reflectbox{\includegraphics[scale=1.75]{FigFol/Figure1_fig_PPlr2_simp.pdf}}}
\\
R\;\Big(
\;
\raisebox{-9pt}{\includegraphics[scale=1.75]{FigFol/Figure1_fig_PPlr1_simp.pdf}}
\;
\Big)
&=
\raisebox{-24pt}{\includegraphics[scale=1.75]{FigFol/Figure1_fig_PPlr2_simp.pdf}}
&
R'\;\Big(
\;
\raisebox{14.2pt}{\reflectbox{\includegraphics[scale=1.75, angle=180]{FigFol/Figure1_fig_PPlr1_simp.pdf}}}
% \raisebox{-9pt}{\includegraphics[scale=1.75]{figure/fig_PPlr3.pdf}}
\;
\Big)
&=
\raisebox{29.2pt}{\reflectbox{\includegraphics[scale=1.75, angle=180]{FigFol/Figure1_fig_PPlr2_simp.pdf}}}
% \raisebox{-9pt}{\includegraphics[scale=1.75]{figure/fig_PPlr4.pdf}}
\rule{0pt}{38pt}
\end{align*}
\caption{\label{fig.defDyn3} The extended Rauzy dynamics.}
\end{figure}

% !!!! HAKING FOR THE SPECIFIC PAGING OF FIGURES !!!!
% \newpage

The extended Rauzy dynamics $\permsex_n$ is defined as follows
\begin{description}
\item[\phantom{$\matchs_n$}\gostrR{$\permsex_n$}\;:]\ The space of
  configuration is $\kSirr_n$, irreducible permutations of size
  $n$. There are four generators, $L$ $L'$, $R$ and $R'$. Let $\s \in \kSirr_n$  
\begin{align*}
L(\s)
&=\gamma_{L,n}(\s(1))\circ \s 
\ef;
& L'(\s)
&=\s\circ \gamma_{L',n}(\s^{-1}(1))
\ef;
\\
R(\s)
&=\s\circ \gamma_{R,n}(\s^{-1}(n))
\ef;
&
R'(\s)
&=\gamma_{R',n}(\s(n))\circ\s
\ef;
\end{align*}
%% \end{description}
%% ppppppppppppppppppp
%% \begin{description}
\end{description}

The extended Rauzy classes (i.e. the equivalence classes on $\kSirr_n$ induced by the $\permsex_n$ dynamics) were characterised by Kontsevich and Zorich in \cite{KZ03} (combinatorial proofs can also be found in \cite{Fic16} \cite{DS17} \cite{D18}) and the theorem can be formulated as follows:

\begin{table}[b!!]
\[
\begin{array}{r||c|cccccc}
n & \Id & \multicolumn{6}{c}{\textrm{non-exceptional classes}}
\\
\hline
4 & 3- & \\
5 & 22 \\
6 & 5+ & 5- \\
7 & 33+ & 24 & 33-
\end{array}
\]
\caption{\label{tab.smallsizeThm2}List of invariants $(\lam,s)$ for
  $n \leq 7$, for which the corresponding class in the $\permsex_n$
  dynamics exists. We shorten $s$ to $\{-,+\}$ if valued $\{-1,+1\}$,
  and omit it if valued~0.}
\end{table}

\begin{theorem}\label{thm.Main_theorem_2}

Besides an exceptional class $\Id$ which has invariants described in table \ref{table_invariant_id_tree}, the number of extended Rauzy classes with cycle invariant $\lam$ (no $\lam_i=1$) depends on the
number of even elements in the list $\{\lambda_i\}$, and is, for 
$n \geq 8$,
\begin{description}
% \item[zero,] if some $\lam_i = 1$, as we are considering primitive classes;
\item[zero,] if
% among $r$ and the $\lam_i$'s 
there is an odd number of even
elements;
\item[\phantom{zero,}\gostrR{one,}] if there is a positive even number of even elements. In which case the value of the sign invariant is 0.
\item[\phantom{zero,}\gostrR{two,}] if there are no even elements at all. The two classes have
non-zero opposite sign invariant.
\end{description}
%% if both $r$ and all $\lam_i$ are odd, and $\lam_i>1$ for
%% all $i$, and they have
%% opposite sign invariant, and there exists exactly one primitive class 
%% with cycle invariant
%% $(\lam,r)$, if both $r$ and all $\lam_i$ are odd, and $\lam_i>1$
For $n \leq 7$ the number of classes with given cycle
invariant may be smaller than the one given above, and the list in
Table \ref{tab.smallsizeThm2} gives a complete account.
\end{theorem}

where the invariants and the exceptional class are defined next section.

%-------------------------------------------------------
\subsection{Definition of the invariants}
\label{ssec.invardefs}
%-------------------------------------------------------

In this section we recall the definition of the
invariants, the proof of their invariance can be found in \cite{DS17}.

% -------------------------------------------------------
\subsubsection{Cycle invariant}
\label{sssec.cycinv}
%-------------------------------------------------------

\noindent
Let $\s$ be a permutation, identified with its diagram.  An edge of
$\s$ is a pair $(i^-,j^+)$, for $j=\s(i)$, where $-$ and $+$ denote
positioning at the bottom and top boundary of the diagram.  Perform
the following manipulations on the diagram: (1)~replace each edge with
a pair of crossing edges; more precisely, replace each edge endpoint,
say $i^-$, by a pair of endpoints, $i_\ell^-$ and $i_r^-$ ($i_\ell$ on the left), then introduce the edges $(i_\ell^-,j_r^+)$ and
$(i_r^-, j_\ell^+)$.  (2)~connect by an arc the points $i_r^\pm$ and
$(i+1)_\ell^\pm$, for $i=1,\ldots,n-1$, both on the bottom and the top of
the diagram; (3)~connect by an path the top-right and bottom-left (respectively the top-left and bottom-right enpoinds  $1_\ell^+$ and $n_r^-$.)
endpoints, $n_r^+$ and $1_\ell^-$ We call this path the bottom path (respectivelly the top path).
\begin{figure}[tb!]
\begin{center}
\begin{tabular}{cp{0cm}c}
\raisebox{12.75pt}{\includegraphics[scale=1]{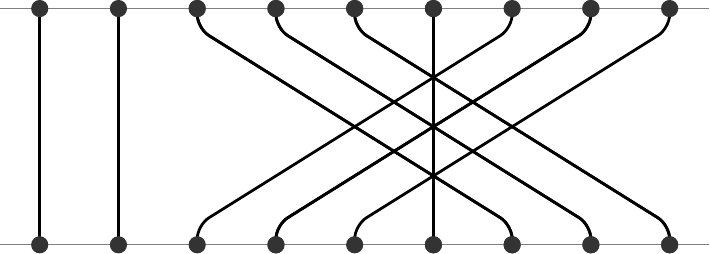}}
&&
\includegraphics[scale=1]{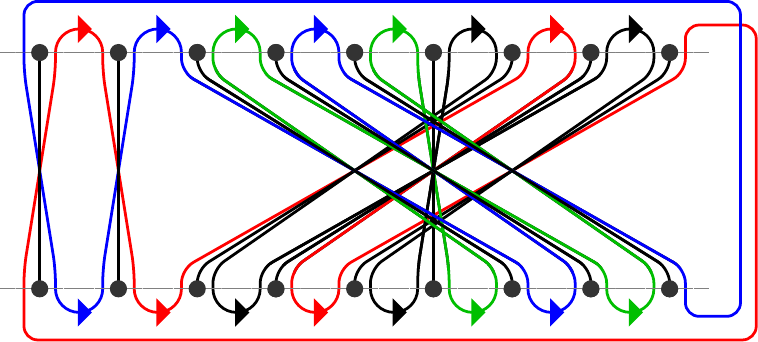}
\end{tabular}
\end{center}
\caption{\label{fig.exLamMatch}Left: a standard permutation,
  $\s=[\,127896345\,]$. Right: the construction of the cycle
  invariant. Different cycles are in different colour. The length of a
  cycle, defined as the number of top (or bottom) arcs, is
  thus 2 for red,blue,black and green. Thus $\lam=\{2,2,2,2\}$. Moreover, the blue cycle is the top principal cycle and the red cycle is the bottom principal cycle.}
\end{figure}

The resulting structure is composed of a number of closed cycles. 
We call the cycle that goes through the top path the top principal cycle and the one that goes through the bottom path the bottom principal cycle. Note that a cycle can be both the top and bottom principal cycle in which case we say that it is the principal cycle.

Define the length of an (open) path as the number
of top (or bottom) arcs (connecting a white endpoint to a black
endpoint) in the path. These numbers are always positive integers (for
$n>1$ and irreducible permutations) and $\lam = \{ \lam_i \}$, the
collection of lengths of the cycles, will be called the \emph{cycle
  invariant} of $\s$. Define $\ell(\s)$ as the number of cycles in
$\s$ minus one.  See
Figure~\ref{fig.exLamMatch}, for an example.

Note that this quantity does \emph{not} coincide with the ordinary
path-length of the corresponding paths. The path-length of a cycle of
length $k$ is $2k$, unless it goes through the top or bottom path, in which
case it is $2k+1$ (if it goes through one of the two) or $2k+2$ (if is goes through both). 

In the interpretation within the geometry of translation surfaces, the
cycle invariant is exactly the collection of conical singularities in
the surface (we have a singularity of $2k\pi$ on the surface, for
every cycle of length $\lam_i=k$ in the cycle invariant, see
article \cite{DS17}).

It is easily seen that
\be
\sum_i \lam_i = n-1
\ef,
\label{eq.size_inv_cycle}
\ee
this formula is called the 
\emph{dimension formula}. 
% \emph{Gauss-Bonnet formula}. 
Moreover, in
the list $\lam=\{\lam_1,\ldots,\lam_{\ell}\}$, there is an even number of
even entries. This is part of theorem \ref{thm_arf_value} stated next section.

We have
\begin{proposition}
\label{prop.cycinv1}
The quantity $(\lam)$ is invariant in the $\permsex$ dynamics.
\end{proposition}
\noindent
For a proof, see \cite{DS17} section 3.1.

We have also shown in \cite{DS17} appendix B that cycles of length 1 have an especially simple behaviour and can thus be omitted from the classification theorem. Thus \textbf{all the classes we consider in this article have a cycle invariant $\lam$ with no parts of length 1.} 

We recall that a cycle of length one correspond in a permutation $\s$ to two edges $(i,j)$ and $(i+1,j-1)$ in case it is not a principal cycle. If it is a principal cycle the notation are slightly more complicated: 

If the top principal cycle has length 1 then we have the three edges $(k,1),(k+1,j),(n,j+1)$.
If the bottom principal cycle has length 1 then we have the three edges $(1,k),(j,k+1),(j+1,n)$.
and if the principal cycle has length 1 then we have the four edges $(1,k),(n,k+1),(j,1),(j+1,n)$.

%-------------------------------------------------------
\subsubsection{Sign invariant}
\label{ssec.arf_inv_intro}
%-------------------------------------------------------

\begin{figure}[tb!]
\[
\begin{array}{cp{4mm}c}
\includegraphics[scale=3]{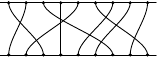}
&&
\includegraphics[scale=3]{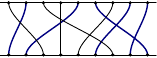}
%% \\
%% \includegraphics[scale=1.6]{FigFol/Figure1_fig_corde_ex1.pdf}
%% &&
%% \includegraphics[scale=1.6]{FigFol/Figure1_fig_corde_ex2.pdf}
\end{array}
\]
\caption{\label{fig.exarf}Left: an example of permutation,
  $\s=[\,251478396\,]$. Right: an example of subset $I=\{1,2,6,8,9\}$
  (labels are for the bottom endpoints, edges in $I$ are in
  blue). There are two crossings, out of the maximal number
  $\binom{|I|}{2}=10$, thus $\chi_I=8$ in this case, and this set
  contributes $(-1)^{|I|+\chi_I}=(-1)^{5+8}=-1$ to $A(\s)$.}
 %% On the
 %%  bottom part of the figure, the same configurations seen as chord
 %%  diagrams.}
\end{figure}

\noindent
For $\s$ a permutation, let $[n]$ be identified to the set of edges
(e.g., by labeling the edges w.r.t.\ the bottom endpoints, left to
right). For $I \subseteq [n]$ a set of edges, define $\chi(I)$ as the
number of pairs $\{i',i''\} \subseteq I$ of non-crossing edges.  Call
\be
\Abar(\s) := \sum_{I \subseteq [n]} (-1)^{|I| + \chi(I)}
\label{eq.arf1stDef}
\ee
the 
\emph{Arf invariant} of $\s$ (see Figure \ref{fig.exarf} for an
example). Call $s(\s)= \mathrm{Sign}(\Abar(\s)) \in \{-1,0,+1\}$ the
\emph{sign} of $\s$.

Both the quantity $\Abar(\cdot)$ and $s(\cdot)$ are invariant for the dynamics $\perms$. The proof can be found in \cite{DS17} section 4.2. 
 There exists an important relationship between the arf invariant and the cycle invariant that we describe below. The proof of this theorem was done in \cite{D18}, theorem 13.

\begin{theorem}\label{thm_arf_value}
Let $\s$ be a permutation with cycle invariant $\lambda$ and let $\ell$ be the number of cycles minus 1 of $\s$ i.e. $\ell=|\lam|-1$.
\begin{itemize}
\item The list $\lambda$ has an even number of even parts. 
\item $\Abar(\s)=\begin{cases}\pm2^{\frac{n+\ell}{2}}& 
% \text{ if the number of even parts in the list $\lambda\bigcup\{r\}$
%   is 0.}\\
\text{ if there are no even parts in the list $\lambda$}\\
0 &\text{ otherwise.} \end{cases}$
\end{itemize}
\end{theorem}

as a consequence we have:
\begin{proposition}\label{pro_sign_inv}
The sign of $\s$ can be written as
$s(\s)=2^{-\frac{n+\ell}{2}}\Abar(\s)$, where $\ell$ is the number of
cycles minus 1 of $\s$.
\end{proposition}

\noindent

% -------------------------------------------------------
\subsection{Exceptional classes}
\label{sssec.ididp}
% -------------------------------------------------------

The invariants described above allow to characterise all classes
for the dynamics on irreducible configurations, \emph{with one
  exception}. This exceptional class is called $\Id_n$ and is the class containing the identity permutation $id_n$.

$\Id_n$ is called the \emph{`hyperelliptic class'},
because the Riemann surface associated to $\Id_n$ is hyperelliptic. 

The class was studied in details in Appendix C.1 of \cite{DS17}. In this article we will only need to know the standard permutations of $\Id_n$. 
\begin{lemma}\label{lem_except_perms}
The standard permutations of $\Id_n$ are :
\begin{center}
\raisebox{15pt}{$id_i=$}\includegraphics[scale=.5]{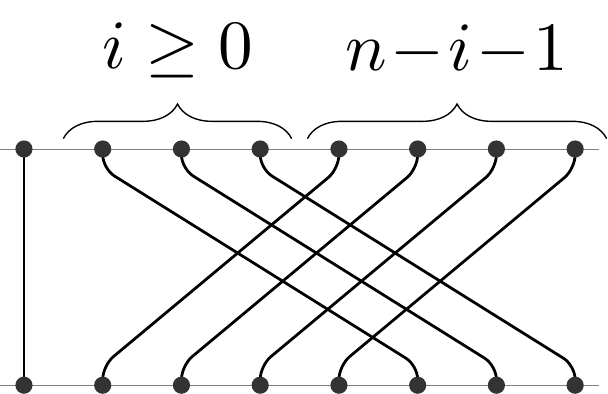}
\end{center}
for $0\leq i \leq n-1$. Clearly $id_0=id$ and $id_i=L^i(id)=L^{\prime\, n-i-1}(id).$
\end{lemma}
\begin{table}[t!!]
\[
\begin{array}{|c|c|c|}
\cline{2-3}
\multicolumn{1}{c|}{}
&
\textrm{$n$ even}
&
\textrm{$n$ odd}
\\
\hline
\lam \text{ of } \Id_n & 
\{n-1\}&
\raisebox{-5pt}{\rule{0pt}{15pt}}%
\{\frac{n-1}{2},\frac{n-1}{2}\}\\
\hline
\end{array}
\]
\[
\begin{array}{|c|cccccccc|}
\cline{2-9}
\multicolumn{1}{c|}{n \mod 8}
&
0&1&2&3&4&5&6&7
\\ \hline
s \text{ of } \Id_n & 
+&0&-&-&-&0&+&+
\\ 
\hline
\end{array}
\]
\caption{Cycle, rank and sign invariants of the exceptional class. The
  sign \hbox{$s \in \{-1,0,+1\}$} is shortened into $\{-,0,+\}$.
\label{table_invariant_id_tree}}
\end{table}

The cycle and sign invariants of the class depend from its size
mod~4, and are described in Table~\ref{table_invariant_id_tree}.  

\section{The Sliding dynamics\label{sec_sliding_dyn}}

It is known since Rauzy \cite{Rau79} that every extended Rauzy class contains standard permutations. The sliding dynamics introduced by Boissy in \cite{BoiCM} is defined on the set of standard permutation $\kS^{St}_n$ as follows:

\begin{figure}[b!]
\begin{align*}
S_e\;\Big(
\;
\raisebox{-11pt}{\includegraphics[scale=.4]{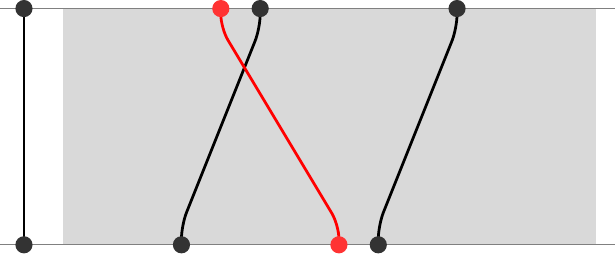}}
% \raisebox{-9pt}{\includegraphics[scale=1.75]{FigFol/Figure1_fig_PPlr3.pdf}}
\;
\Big)
&=
\raisebox{-11pt}{\includegraphics[scale=.4]{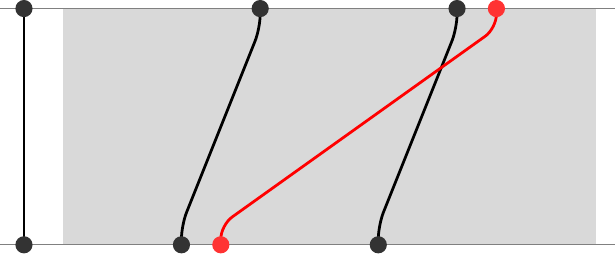}}
% \raisebox{-9pt}{\includegraphics[scale=1.75]{FigFol/Figure1_fig_PPlr4.pdf}}
\rule{0pt}{38pt}
&
S_e\;\Big(
\;
\raisebox{-11pt}{\includegraphics[scale=.4]{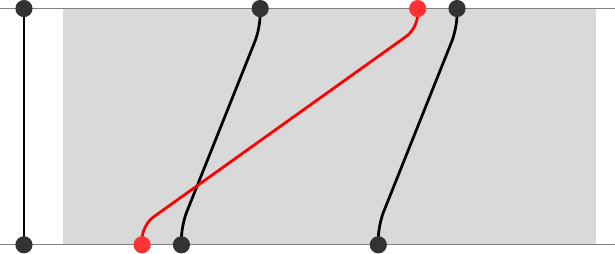}}
\;
\Big)
&=
\raisebox{-11pt}{\includegraphics[scale=.4]{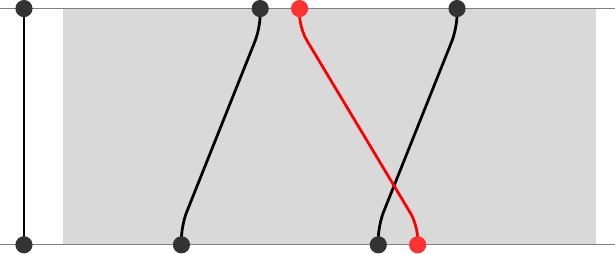}}\\
S_e\;\Big(
\;
\raisebox{-11pt}{\includegraphics[scale=.4]{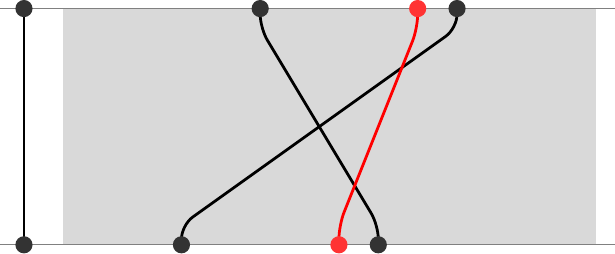}}
% \raisebox{-9pt}{\includegraphics[scale=1.75]{FigFol/Figure1_fig_PPlr3.pdf}}
\;
\Big)
&=
\raisebox{-11pt}{\includegraphics[scale=.4]{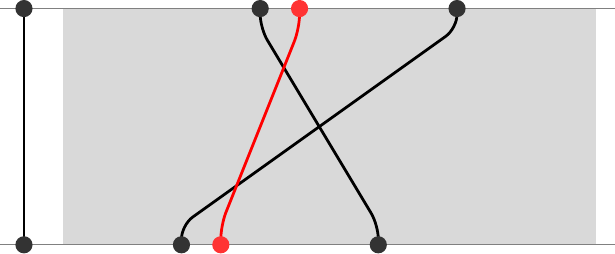}}
% \raisebox{-9pt}{\includegraphics[scale=1.75]{FigFol/Figure1_fig_PPlr4.pdf}}
\rule{0pt}{38pt}
&
S_e\;\Big(
\;
\raisebox{-11pt}{\includegraphics[scale=.4]{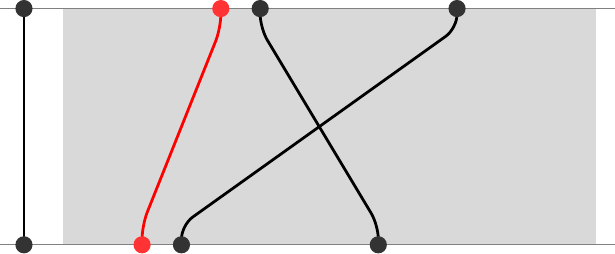}}
\;
\Big)
&=
\raisebox{-11pt}{\includegraphics[scale=.4]{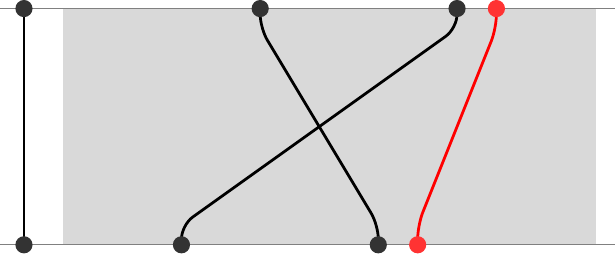}}\\
S_e\;\Big(
\;
\raisebox{-11pt}{\includegraphics[scale=.4]{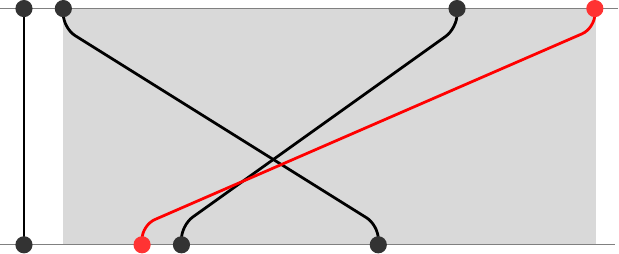}}
% \raisebox{-9pt}{\includegraphics[scale=1.75]{FigFol/Figure1_fig_PPlr3.pdf}}
\;
\Big)
&=
\raisebox{-11pt}{\includegraphics[scale=.4]{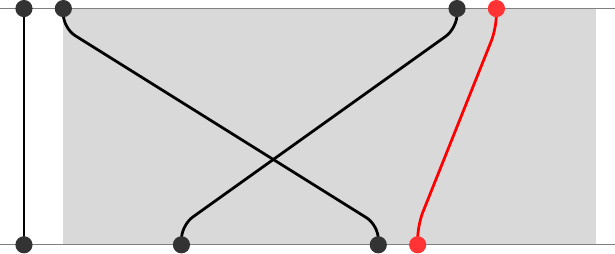}}
% \raisebox{-9pt}{\includegraphics[scale=1.75]{FigFol/Figure1_fig_PPlr4.pdf}}
\rule{0pt}{38pt}
&
S_e\;\Big(
\;
\raisebox{-11pt}{\includegraphics[scale=.4]{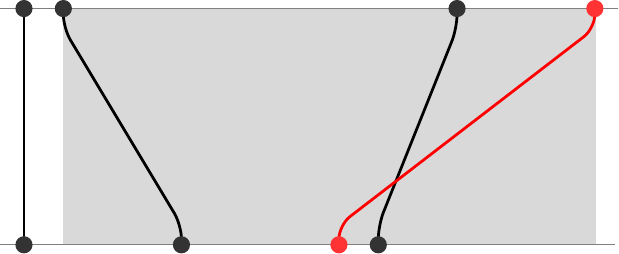}}
\;
\Big)
&=
\raisebox{-11pt}{\includegraphics[scale=.4]{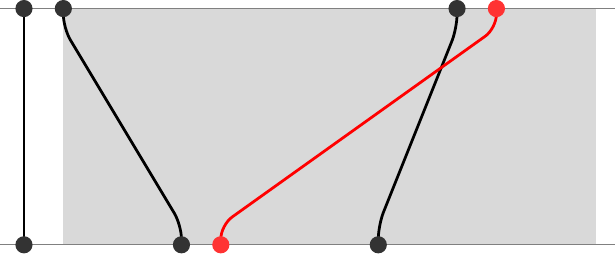}}\\
S_e\;\Big(
\;
\raisebox{-11pt}{\includegraphics[scale=.4]{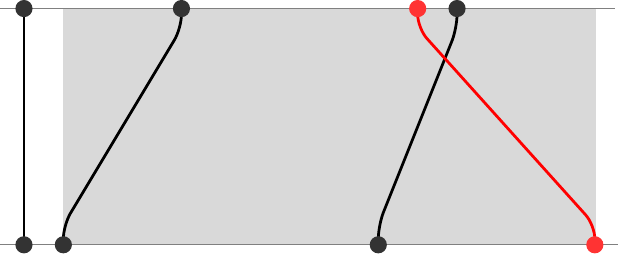}}
% \raisebox{-9pt}{\includegraphics[scale=1.75]{FigFol/Figure1_fig_PPlr3.pdf}}
\;
\Big)
&=
\raisebox{-11pt}{\includegraphics[scale=.4]{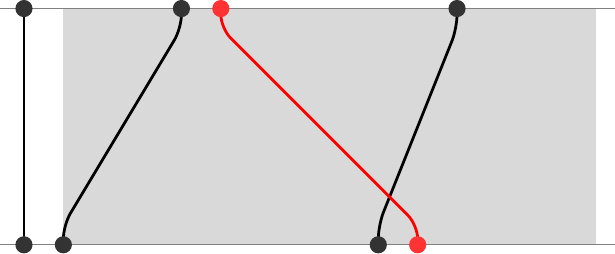}}
% \raisebox{-9pt}{\includegraphics[scale=1.75]{FigFol/Figure1_fig_PPlr4.pdf}}
\rule{0pt}{38pt}
&
S_e\;\Big(
\;
\raisebox{-11pt}{\includegraphics[scale=.4]{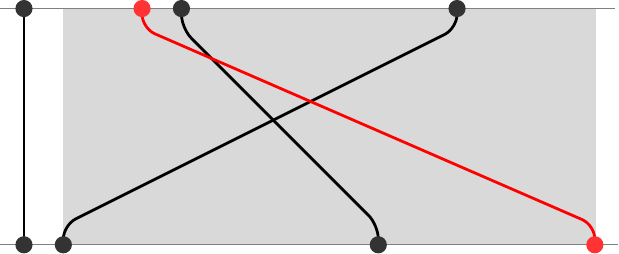}}
\;
\Big)
&=
\raisebox{-11pt}{\includegraphics[scale=.4]{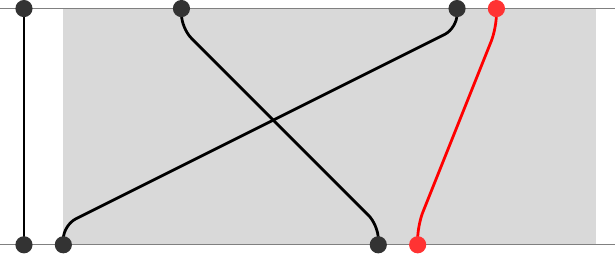}}\\
S_e\;\Big(
\;
\raisebox{-11pt}{\includegraphics[scale=.4]{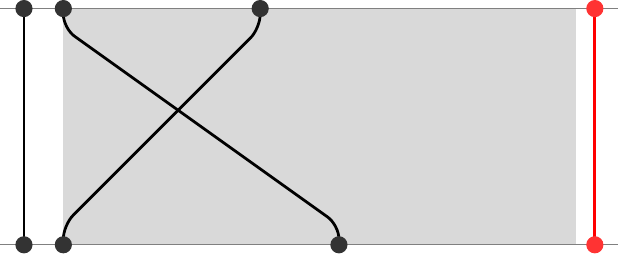}}
% \raisebox{-9pt}{\includegraphics[scale=1.75]{FigFol/Figure1_fig_PPlr3.pdf}}
\;
\Big)
&=
\raisebox{-11pt}{\includegraphics[scale=.4]{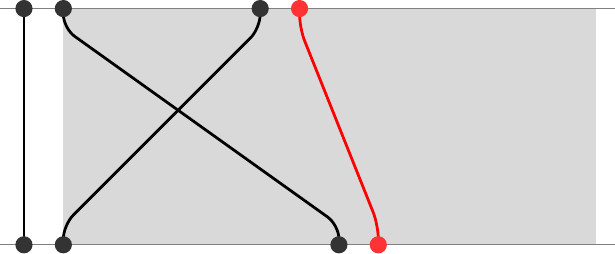}}
% \raisebox{-9pt}{\includegraphics[scale=1.75]{FigFol/Figure1_fig_PPlr4.pdf}}
\rule{0pt}{38pt}
&
S_e\;\Big(
\;
\raisebox{-11pt}{\includegraphics[scale=.4]{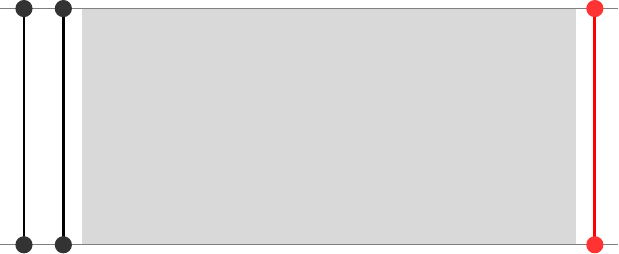}}
\;
\Big)
&=
\raisebox{-11pt}{\includegraphics[scale=.4]{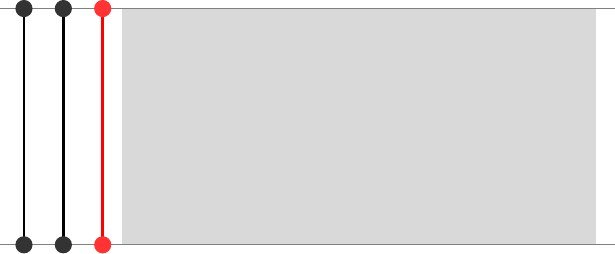}}\\
&
\rule{0pt}{38pt}
&
S_e\;\Big(
\;
\raisebox{-11pt}{\includegraphics[scale=.4]{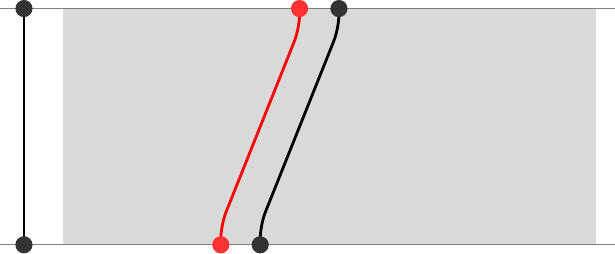}}
\;
\Big)
&=
\raisebox{-11pt}{\includegraphics[scale=.4]{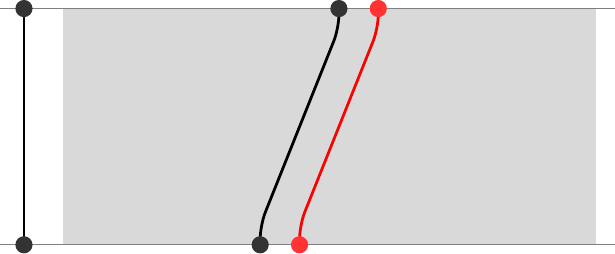}}
\end{align*}
\caption{\label{fig_Secase} All the possible cases for the $S_e$ operator (the edge $e$ is represented in red in the figure). Graphically $S_e$ makes the edge $e$ slide to the right along its two right-adjacent edges. To limit the number of cases we did not represent that when one of the endpoints of $e$ should end up to the position $n$ after the application of $S_e$ we shift it to position 2.}
\end{figure}
\begin{description}
\item[\phantom{$\matchs_n$}\gostrR{$\perms l_n$}\;:]\ The space of
  configuration is $\kS^{St}_n$, standard permutations of size
  $n$. There are $n+1$ generators, $L$ $L'$ as above, and one operator $S_e$ for each edge $e$ other than the edge $e_0=(1,1)$. Let $\s \in \kS^{St}_n$, $e=(i,\s(i))$ and define 
\[e_t=\begin{cases}(\s^{-1}(\s(i)+1),\s(i)+1)&\text{if } \s(i)<n \\
(\s^{-1}(2),2)& \text{if } \s(i)=n \end{cases}
\text{ and } 
e_b=\begin{cases}(i+1,\s(i+1))& \text{if } i<n \\
(2,\s(2))&\text{if } i=n \end{cases}\]
Graphically $e_t$, (respectively $e_b$) is the edge whose top endpoint (respectively bottom endpoint) is to the right of the top (respectively bottom) endpoint of $e$ (with the additionnal identification that the right of the rightmost bottom/top endpoint is the second bottom/top endpoint.)

Then $S_e(\s)$ is obtained from $\s$ in two steps: first we remove the edge $e$ from $\s$ and then we add an edge $e'$ whose bottom endpoint is to the right of the bottom endpoint of $e_t$ and top endpoint is to the right of the top endpoint of $e_t$. Recall that the right of the rightmost endpoint is the second endpoint so if for example $e_b=(i+1,n)$ then $e'=(j,2)$ for some $j$. See figure \ref{fig_Secase} for the possible cases depending on the values of $e,e_t$ and $e_b$. 
%% \end{description}
%% ppppppppppppppppppp
%% \begin{description}
\end{description}

\begin{remark}\label{rk_general} Clearly $L$ and $L'$ commutes with each other. The operators $(S_e)_e$ almost commute with $L$ and $L'$ indeed, for every $e$, we have:
\begin{align*}
LS_e(\s)=S_eL^{i(\s)}(\s) \quad & S_eL(\s)=L^{j(\s)}S_e(\s)\\
L'S_e(\s)=S_eL'^{i'(\s)}(\s)  \quad & S_eL'(\s)=L'^{j'(\s)}S_e(\s)
\end{align*}

 In particular this implies that when checking a property for the operators $(S_e)_e$, if we already know the property true for the operators $L$ and $L'$,  we can always place ourself in the case of figure \ref{fig_Secase} first line since all the other cases can be obtained from it by application of $L^i$ and $L'^j$ for some $i,j$. 

For example case second line left turns into case first line left by application of $L^i$ for some $i$. See below.
\[ 
\begin{tikzcd}
\raisebox{-10pt}{\includegraphics[scale=.4]{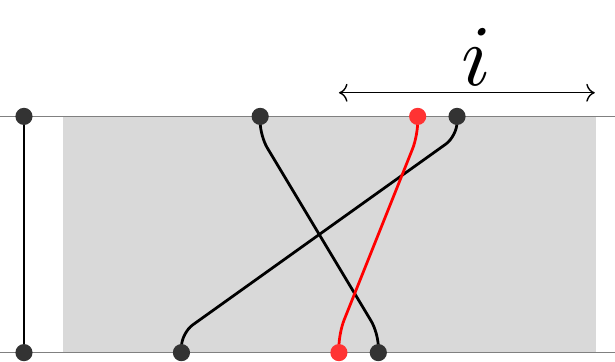}} \arrow{r}{L^i} \arrow[bend left=20]{rrr}{S_e} & \raisebox{-10pt}{\includegraphics[scale=.4]{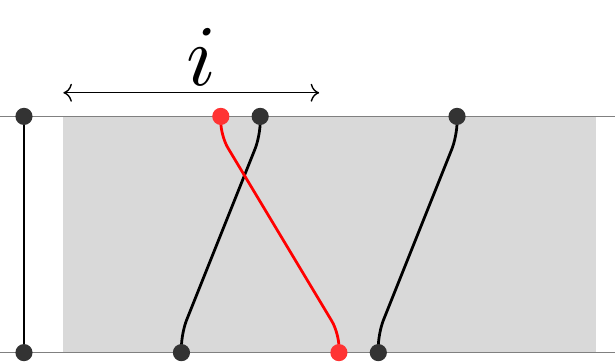}} \arrow{r}{S_e}& \raisebox{-10pt}{\includegraphics[scale=.4]{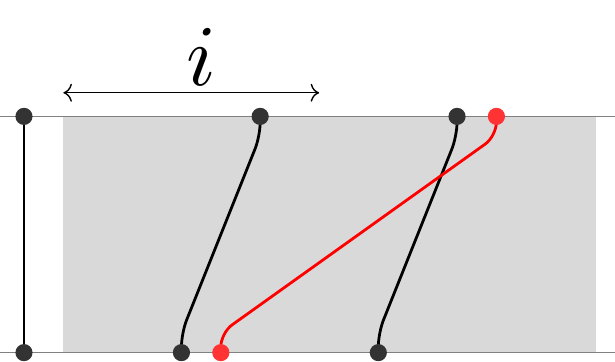}} \arrow{r}{L^{-i}} & \raisebox{-10pt}{\includegraphics[scale=.4]{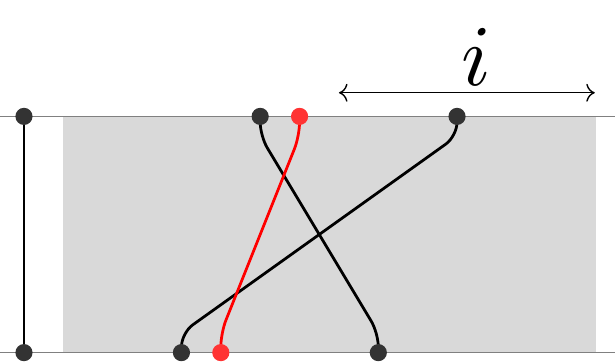}}
\end{tikzcd}
\]
 We will make use of this to verify that the sign invariant is indeed invariant for $S_e$ (see proposition \ref{pro_arf_inv}).
\end{remark}

\begin{remark}\label{rk_never_at_end}
Due to our definition of $S_e$ the edge $e$ will never have one endpoint at position $n$ in $S_e(\s)$, this will be important in the next section (cf text below proposition \ref{pro_Se_action}). 
\end{remark}
The main theorem of this first part of the paper is the following:

\begin{theorem}\label{thm_main}
The Rauzy classes of the Sliding dynamics are the restriction of the extended Rauzy classes (of the $\permsex$ dynamics) to the set of standard permutations.

Thus two permutations are connected for the $\perms l_n$ dynamics if they are both in the exceptional class or if neither are in the exceptional class and they have the same cycle invariant and sign invariant.
\end{theorem}

Let us first prove that the dynamics is invertible and keeps invariant the cycle invariant and sign invariant.
For this purpose, we reintroduce the notations developped in \cite{D18} on arcs and intervals adapted to the extended Rauzy classes.

\subsection{Consistent labelling}

In our previous article \cite{D18} we introduced the notion of consistent labelling (a labelling of the intervals between adjacent vertices (top or bottom) of the permutation or equivalently a labelling of the arcs added in the cycle invariant structure) for the non-extended Rauzy classes. We reproduce here the definition with some modifications to better support the fact that we are working with the extended Rauzy classes instead.
%Our notion of arcs departs slightly from the definition inherited from
%the combinatorial structure with arcs as we define both bottom and top
%arcs in our diagram representation of permutations. 

% \begin{definition}
Let $\s$ be a permutation of size $n$. The procedure to construct the
cycle invariant $(\lambda)$ (as described in Section
\ref{sssec.cycinv}) involves the introduction of $n-1$ top and bottom
arcs connecting adjacent top and bottom vertices. 

We number the top and bottom arcs from left to right, and refer to
them by their position: the bottom arc $\beta\in \{1,\ldots,n-1\}$ is
the $\beta$-th bottom arc, counting from the left. Likewise for the top arc $\alpha\in \{1,\ldots,n-1\}$.
% \begin{center}
\[
\begin{array}{cc}
\raisebox{10.9pt}{ \includegraphics[scale=.5]{figure/fig_n_arcs.pdf}}
&
\raisebox{11pt}{\includegraphics[scale=.5]{figure/fig_no_arc.pdf}}
\end{array}
\]
% \end{center}

By convention, the variables used to name the positions of the top
(bottom) arcs will be $\alpha$ (respectively $\beta$), in order not to
make confusion with other parts of the diagram (for which we will use
$i,j,\ldots$ or $x,y,\ldots$).
% \end{definition}

\begin{definition}\label{def_consecutive}
We say that two (bottom) arcs $\beta,\beta'$ are \emph{consecutive}
(in a cycle) if they are inside the same cycle and they are consecutive in the cyclic order induced by the cycle. This occurs in one of the three
graphical patterns:
\begin{center}\begin{tabular}{ccc}
    \includegraphics[scale=.5]{figure/fig_consecutive_bot_1.pdf} &
    \includegraphics[scale=.5]{figure/fig_consecutive_bot_2.pdf} &
    \includegraphics[scale=.5]{figure/fig_consecutive_bot_3.pdf} 
\end{tabular}
\end{center}
In formulas:
\[
\beta'=\s^{-1}(\s(\beta+1)+1) \text{ if } \s(\beta+1) <n \text{ and }
\beta'=\s^{-1}(\s(1)+1) \text{ if } \s(\beta+1) =n
\]
We define consecutive arcs for top arcs similarly:
\begin{center}\begin{tabular}{ccc}
    \includegraphics[scale=.5]{figure/fig_consecutive_top_1.pdf} &
    \includegraphics[scale=.5]{figure/fig_consecutive_top_2.pdf} &
    \includegraphics[scale=.5]{figure/fig_consecutive_top_3.pdf} 
\end{tabular}
\end{center}
In formulas:
\[
\alpha'=\s(\s^{-1}(\alpha)+1)+1 \text{ if } \s^{-1}(\alpha) <n \text{ and }
\alpha'=\s(\s^{-1}(1))+1 \text{ if } \s^{-1}(\alpha) =n
\]

\end{definition}

\begin{remark}
As we have seen above, when representing graphically the consecutive
arcs, we need three figures depending on the different cases (edges
crossing or not, and edges ending at a left corner of the
diagram). However, these cases are treated in a very similar way, and,
in the graphical explanation of our following properties, we shall
mostly draw consecutive arcs by representing the case of non-crossing
and non-corner edges, i.e.\ the left-most of the drawings above. It is
intended that the underlying reasonings remain valid for the other
cases.
\end{remark}

Next we define suitable alphabets used to label the top and bottom
arcs of a permutation.

\begin{notation}
For all $j$, let $\Sigma_{i,j}=\{b_{0,i,j},\ldots,b_{2i-2,i,j}\}$ and
$\Sigma'_{i,j}=\{t_{1,i,j},\ldots,t_{2i-1,i,j} \}$ be a pair of
alphabets which label the bottom arcs and the top arcs respectively of
a cycle of length $i$.
\end{notation}

\begin{figure}[bt]
\begin{center}\begin{tabular}{c} \includegraphics[scale=1]{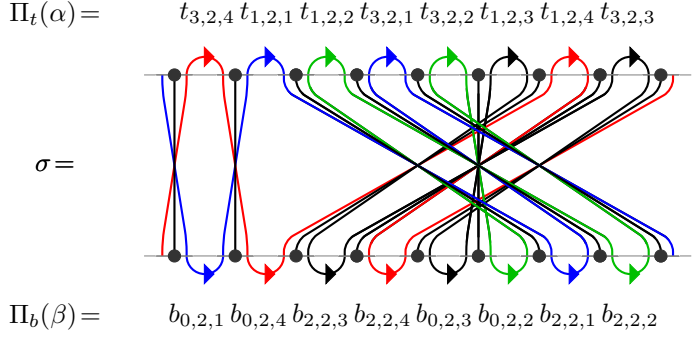}   \end{tabular}\end{center}
\caption{\label{fig_example_labelling} A consistent labellings $(\Pi_b,\Pi_t)$ of a permutation $\s$ with cycle invariant $(\{2,2,2,2\})$.}
\end{figure}

Finally, we can introduce our notion of consistent labelling (for the extended Rauzy classes).
% In essence, a consistent labelling labels the arcs of the
% permutation in a way  that respects the cycle invariant. 
\begin{definition}[Consistent labelling]
\label{def_consistent_lab}
Let $\s$ be a permutation with invariant
$(\lam=\{\lam_1^{m_1},\ldots,\lam_k^{m_k}\})$ and define a
\emph{consistent labelling} to be a pair $(\Pi_b,\Pi_t)$ of
bijections:
\[
\begin{array}{ccccc}
\Pi_b & : & \{1,\ldots,n-1\} & \to & \Sigma_b= 
 \bigcup_{i=1}^k \Big( \bigcup_{j=1}^{m_i}
  \Sigma_{\lam_i,j} \Big) 
\\
\Pi_t & : & \{1,\ldots,n-1\} &\to &\Sigma_t= 
\bigcup_{i=1}^k \Big( \bigcup_{j=1}^{m_i}
  \Sigma'_{\lam_i,j} \Big) 
\end{array}
\]
such that
\begin{enumerate}
\item Two arcs within the same cycle have labels within the same
  alphabet. Thus if $S_b=\{(\beta_k)_{1\leq k \leq \lam_i}\}$ and
  $S_t=\{(\alpha_k)_{1\leq k \leq \lam_i}\}$ are the sets of bottom
  (respectively top) arcs of a cycle of length  $\lam_i$, then
  $\Pi_b(S_b)=\Sigma_{\lam_i,j}$ and $\Pi_t(S_t)=\Sigma'_{\lam_i,j}$
  for some $1 \leq j\leq m_i$.
\item Two consecutive arcs of a cycle of length $\lam_i$ have labels
  with consecutive indices: if $\beta$ and $\beta'$ are consecutive,
  then $\Pi_b(\beta)=b_{k,\lam_i,j}$ for some $k<\lam_i$ and $j\leq
  m_i$ and $\Pi_b(\beta')=b_{k+2,\lam_i,j}$, where $k+2$ is intended
  modulo $2\lam_i$. Likewise for top arcs.\\[-7mm]
\begin{center}
\begin{tabular}{cc}
  \includegraphics[scale=.5]{figure/fig_consecutive_bot_lab.pdf} &
  \raisebox{24pt}{\includegraphics[scale=.5]{figure/fig_consecutive_top_lab.pdf}}
\end{tabular}
\end{center}
~\\[-12mm]
\item The bottom right arc $\beta$ and the top left arc $\alpha$ of an
  edge $i$ are labeled by consecutive indices:
\[
\text{ if } \beta=i,\alpha=\s(i), \text{ then }\, 
\Pi_t(\alpha)=t_{k,\lam_{\ell},j} \ \Leftrightarrow \ 
\Pi_b(\beta)=b_{k+1,\lam_\ell,j} 
\] \\[-12mm]
\begin{center}
\begin{tabular}{cc}
\includegraphics[scale=.5]{figure/fig_label_top_bot} &     \includegraphics[scale=.5]{figure/fig_consecutive_top_3_top_bot.pdf} 
\end{tabular}
\end{center}~\\[-14mm]
\end{enumerate}
\end{definition}
\noindent
Figure~\ref{fig_example_labelling} provides an example of consistent
labelling.

\begin{lemma}\label{rk_top_from_bot}
Let $\s$ be a permutation and 
$\Pi_b \,:\,\{1,\ldots,n-1\} \to \Sigma_b$ a labelling of bottom arcs
verifying property 1 and 2 of Definition
\ref{def_consistent_lab}. Then there exists a unique 
$\Pi_t \,:\,\{1,\ldots,n-1\} \to \Sigma_t$ such that $(\Pi_b,\Pi_t)$ is
a consistent labelling.
\end{lemma}

\begin{pf}
Let  $(\Pi_b,\Pi_t)$ a be a consistent labelling. Then by property 3
we must have 
$\Pi_t(\alpha)= t_{i,\lam_\ell,j}  \text{ if }
  \Pi_b(\s^{-1}(\alpha))=b_{i+1 \!\!\! \mod 2\lam_\ell,\lam_\ell,j}.$
This uniquely defines $\Pi_t$.  
\qed
\end{pf}

\noindent
This lemma implies that the data $(\s,\Pi_b)$ or $(\s,\Pi_t)$ are
sufficient to reconstruct $(\s,(\Pi_b,\Pi_t))$. Thus, occasionally, we
will consider just $(\s,\Pi_b)$ rather than $(\s,(\Pi_b,\Pi_t))$.

\subsection{Cycle invariant and edge addition\label{sec_edge_addition}}

This preliminary section study the change of the cycle invariant when inserting a few consecutives edges  in a permutation. 
%The results of this section will be used in sections \ref{sec_arf},\ref{sec_IIXshaped},\ref{sec_induction}.

\begin{figure}[bt!]
\begin{center}
\includegraphics{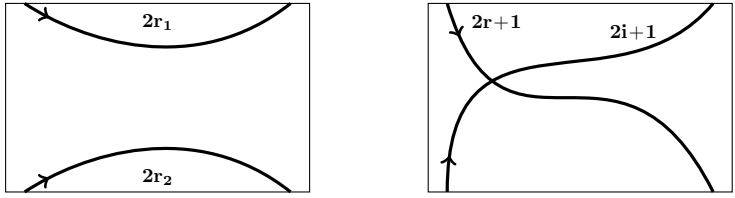}
\caption{\label{fig.typeXH}Left: a schematic representation of a
  permutation of type $H(r_1,r_2)$. The quantity $2r_1$ and $2r_2$ are the path lenght of the principal cycle. Right: a representation of a  permutation of type $X(r,j)$. The quantity $2r+1$ and $2i+2$ are the path lenght of the two principal cycles.}
  \end{center}
\end{figure}

First we reintroduce a notation from \cite{DS17} (it was also defined as $m|a$ and $m_l \cdot m_r$ in \cite{Del13}).
\begin{definition}
\label{def.XHtype}
A permutation $\s$ is \emph{of type $H$} if the top principal cycle and the bottom principal cycle are distinct, and \emph{of type $X$} otherwise.
\end{definition}
\noindent See also figure \ref{fig.typeXH}.

\begin{figure}[bt!]
\begin{center}
\begin{tabular}{cc}
\includegraphics[scale=.5]{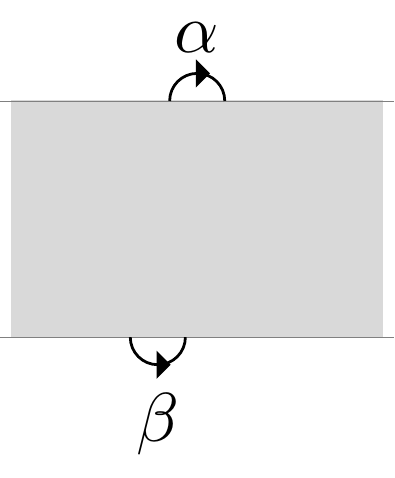}
&
\includegraphics[scale=.5]{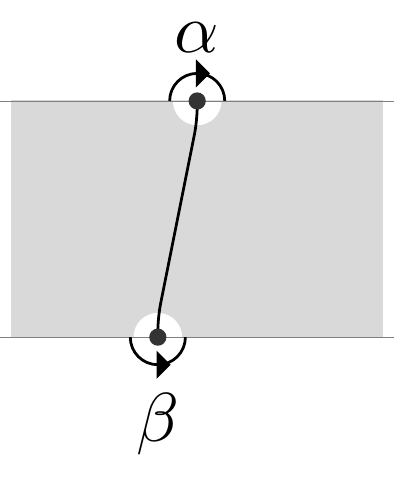}\\
$\s$ 
&
$\s|(1,\alpha,\beta)$
\end{tabular}
\end{center}
\caption{\label{fig_one_edge_example}The insertion of one edge within the arcs $\alpha$ and $\beta$.}
\end{figure}

\begin{notation}\label{not_add_edge}
let $\s$ be a permutation and let $\alpha$ be a top arc and $\beta$ be a bottom arc, we define $\s|(i,\alpha,\beta)$ to be the permutation obtained from $\s$ by inserting $i\in \mathbb{N}$ consecutive and parallel edges within $\alpha$ and $\beta$. (see figure \ref{fig_one_edge_example} for an example with $i=1$).

In the notation $\s|(i,\cdot,\cdot)$, we will refer to the arcs either by their position $\alpha,\beta$ or their labels $b,t \in \Sigma$ if a consistent labelling is defined. 
\end{notation}

\begin{definition}[(n-1,1)-coloring and reduction]\label{def_coloring}
Let $\s$ be a permutation, a coloring $c$ of $\s$ is a coloring of the edges of $\s$ into a black set of $n-1$ edges and one gray edge $e$. We call $\tau$ the \emph{reduction} of $\s$ if it is the restriction of $(\s,c)$ to the set of black edges.

Thus we have $\s=\tau|(1,\alpha,\beta)$ where $\alpha$ and $\beta$ are the positions of the top arc (respectively) bottom arc containing the gray edge of $(\s,c)$ in $\tau$.
\end{definition}

Note that we can define $\s=\tau|(1,\alpha,\beta)$ only if $e=(i,j)$ does not have $i=n$ or $j=n$ since (cf the previous section) we did not define the n-th (top or bottom) arc. However that will not be a problem since if $e=(i,j)$ with $i=n$ or $j=n$ then $L'(\s)$ or $L(\s)$ or $LL'(\s)$  (depending on whether $i=n$ or $j=n$ or both egal $n$) does not.

\begin{proposition}[One edge insertion into two cycles]\label{pro_one_edge}
Let $\tau$ be a permutation with cycle invariant $(\lambda)$. 

 Let $\tau|(1,\alpha,\beta)$ be the permutation resulting from the insertion of an edge within two arcs $\alpha$ and $\beta$ of two differents cycles of respective length $\ell$ and $\ell'$. Then the cycle invariant of  $\tau|(i,\alpha,\beta)$ is $\lambda\setminus\{\ell,\ell'\}\bigcup\{\ell+\ell'+1\})$.

\end{proposition}

\begin{pf}
See figure \ref{fig_add_one_edge}.
\begin{figure}[h!]
\begin{center}
\begin{tabular}{ll}
&$\tau, (\lambda\bigcup\{\ell,\ell'\},r) \qquad \qquad\qquad \tau|(1,\alpha,\beta), (\lambda\bigcup\{\ell\!+\!\ell'\!+\!1\},r)$\\[-2mm]
&\put(22,103){$\ell$}\put(22,103){$\ell$}
\put(65,90){$\ell'$}\put(65,90){$\ell'$}
\put(165,98){$\ell\!+\!\ell'\!+\!1$}\put(165,98){$\ell\!+\!\ell'\!+\!1$}
\put(30,43){$2r_1$}\put(30,43){$2r_1$}
\put(23,32){$2p_1$}\put(23,32){$2p_1$}
\put(74,55){$2r_2\!+\!1$}\put(74,55){$2r_2\!+\!1$}
\put(64.5,17){$2p_2\!+\!1$}\put(64.5,17){$2p_2\!+\!1$}
\includegraphics[scale=1]{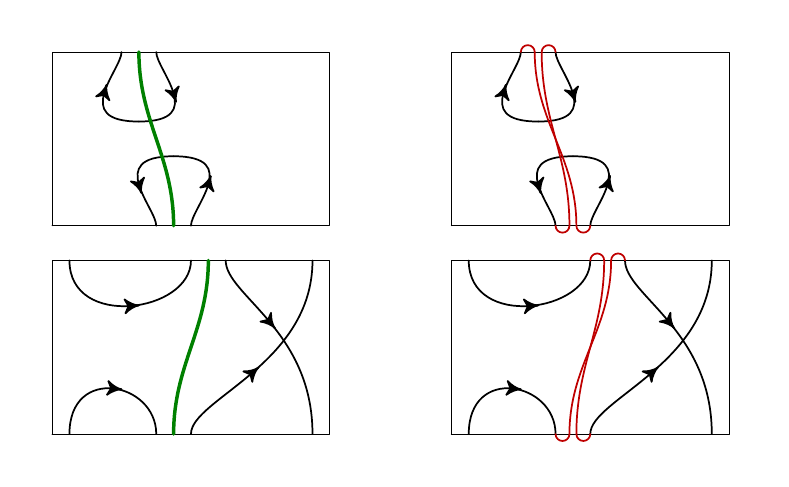}\\[-2mm]
&$\tau, (\lambda\bigcup\{p1\!+\!p2\},r_1\!+\!r_2) \qquad  \s|(1,\alpha,\beta), (\lambda,r_1\!+\!r_2\!+\!p_1\!+\!p_2\!+\!1)$\\[.5mm]
Type: & $\ X(r_1\!+\!r_2,p_1\!+\!p_2) \qquad\qquad  \quad H(r_1+p2+1,r_2+p_1+1)$
\end{tabular}
\end{center}
\caption{\label{fig_add_one_edge}The first line represents the case: Top arc : any cycle. Bottom arc: any cycle. The second line represents the case: Top arc : top principal cycle. Bottom arc: bottom principal cycle.}
\end{figure}
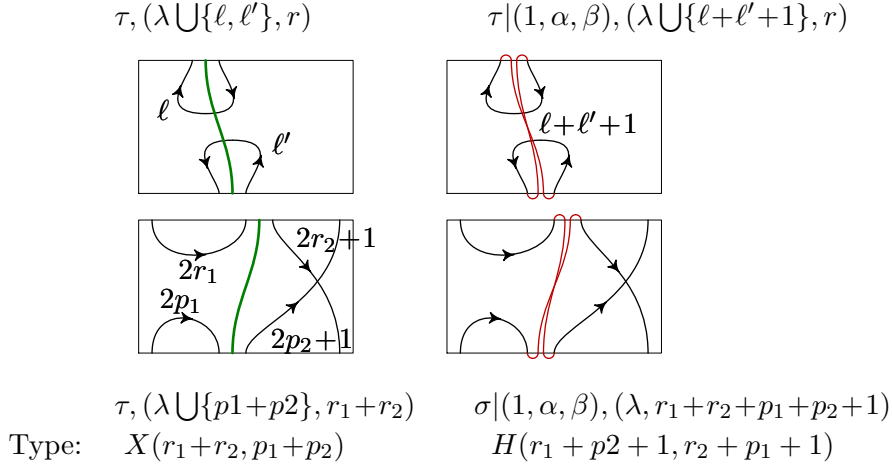

 Some cases are not represented in the figure.
The missing cases are: \begin{itemize}
 \item Top arc: Top principal cycle. Bottom arc: any cycle. 
 \item Top arc: bottom principal cycle. Bottom arc: any cycle. 
 \item Top arc: bottom principal cycle. Bottom arc: top principal cycle. 
 \item Top arc: any cycle. Bottom arc: top principal cycle. 
 \item Top arc: any cycle. Bottom arc: bottom principal cycle. 
\end{itemize} 
Their proof is nearly identical to the ones represented in the figure and are thus omitted. 
\qed
 \end{pf}
 
\begin{figure}[h!]
\begin{center}
\begin{tabular}{cc}
$\tau$ & $\s=\tau|(1,t_{3,2,1},b_{2,3,1})$\\
\includegraphics[scale=.5]{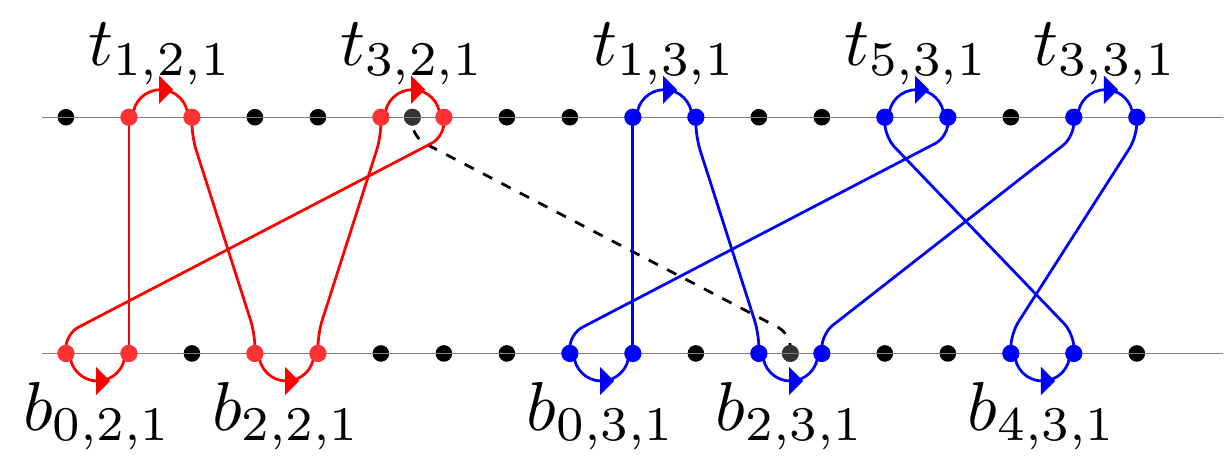} &
\raisebox{12pt}{\includegraphics[scale=.5]{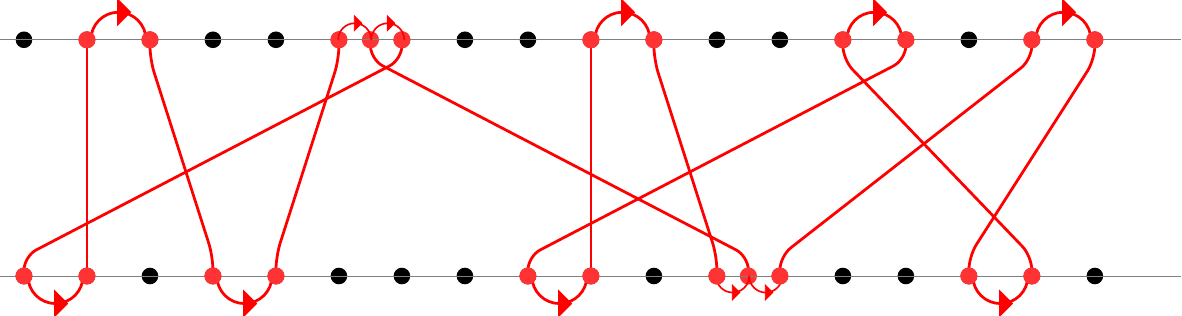}}
\end{tabular}
\end{center}
\caption{\label{figure_ex_appar_arc_merge} The two cycles of length 2 and 3 are merge into a cycle of length 6. The two arcs of $\tau$ containing the edge are broken into two arcs in $\s$ while the others are untouched.}
\end{figure}

The following lemma indicates the correspondance between the arcs of $\tau$ and the arcs of $\s$:

\begin{lemma}\label{lem_where_arc}
Let $\s$ and $\tau$ be as above and let $\Pi$ be a consistent labelling of $\tau$ such that $\s=\tau|(1,t_{2x+1,\ell,k},b_{2y,\ell,k'})$. Consider the arcs with labels $t_{2u+1,\ell,k}, t_{2u'+1,\ell',k'}, b_{2v,\ell,k}, b_{2v',\ell',k'}$ in $\tau$. Those arcs correspond to arcs in $\s$ (one to one except for the arcs $t_{2x+1,\lam_j,k},b_{2y,\lam_j,k}$ in $\tau$ which correspond to two arcs in $\s$ since we inserted the edge rigth within them) and all those arcs are part of the new cycle of length $\ell+\ell'+1$. 

Of course all the arcs part of the other cycles are unchanged (thus correspond one-to-one).
See figure \ref{figure_ex_appar_arc_merge}.
\end{lemma}

\begin{notation}
let us define the following quasimetric on $\mathbb{Z}_n:$ $q_n(x,y)=\begin{cases} y-x &\text{if } y\geq x\\ n-x+y &\text{if } x>y\end{cases}$. $q_n(x,y)$ is the smallest nonnegative integer such that $x+q_n(x,y)\equiv y \mod n$.
\end{notation}

\begin{proposition}[One edge insertion into one cycle]\label{pro_one_edge_two}
Let $\tau$ be a permutation with cycle invariant $(\lambda)$ and let $\Pi=(\Pi_t,\Pi_b)$ be a consistent labelling.

 Let $\tau|(1,t_{2x+1,\lam_j,k},b_{2y,\lam_j,k})$ be the permutation resulting from the insertion of an edge within two arcs of the same cycle of length $\lam_j$. Then the cycle invariant of $\tau|(1,t_{2x+1,\lam_j,k},b_{2y,\lam_j,k})$ is \[\lambda\setminus\{\lam_j\}\bigcup\{\frac{q_{2\lam_j}(2x+1,2y)+1}{2},\frac{q_{2\lam_j}(2y,2x+1)+1}{2}\}\].

\end{proposition}

\begin{pf}
See figure \ref{fig_add_one_edge_one_cycle}. $q_{2\lam_j}(2x\!+\!1,2y)$ represente the path length between $t_{2x+1,\lam_j,k}$ and $b_{2y,\lam_j,k}$ and $q_{2\lam_j}(2y,2x\!+\!1)$ the path length between $b_{2y,\lam_j,k}$  and $t_{2x+1,\lam_j,k}$. Thus when inserting the edge within $t_{2x+1,\lam_j,k}$ and $b_{2y,\lam_j,k}$ the cycle of length $\lam_j$ is broken into two cycles, one of length $\frac{q_{2\lam_j}(2x+1,2y)+1}{2}$ and the other of length $\frac{q_{2\lam_j}(2y,2x+1)+1}{2}$.
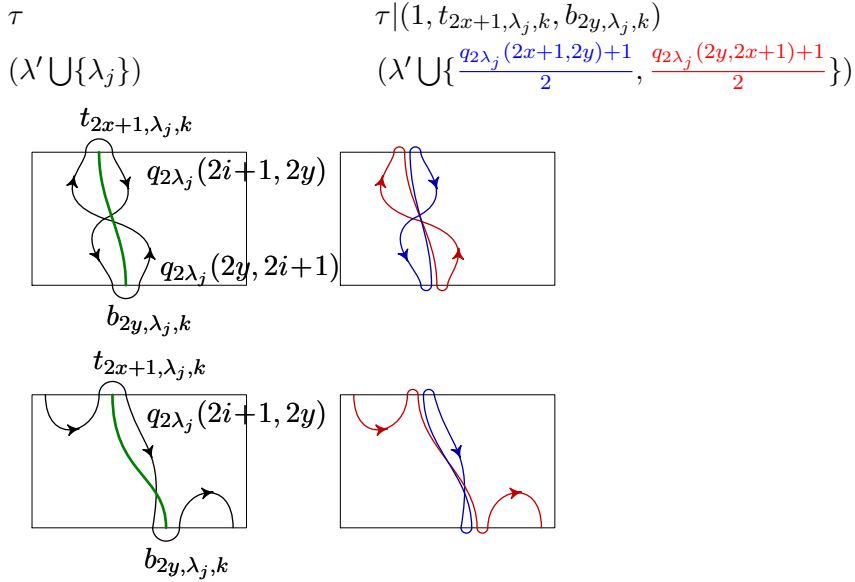
\begin{figure}[h!]
\begin{center}
\begin{tabular}{ll}
&$\tau \qquad\qquad\qquad\qquad\qquad\qquad \tau|(1,t_{2x+1,\lam_j,k},b_{2y,\lam_j,k})$\\[1mm]
&$ (\lambda' \bigcup \{\lam_j\})\qquad\qquad\qquad\qquad (\lambda'\bigcup\{{\color{blue}\frac{q_{2\lam_j}(2x+1,2y)+1}{2}},{\color{red}\frac{q_{2\lam_j}(2y,2x+1)+1}{2}} \})$\\[2.5mm]
&
\put(33,75){$t_{2x+1,\lam_j,k}$}\put(33,75){$t_{2x+1,\lam_j,k}$}
\put(58,55){$q_{2\lam_j}(2i\!+\!1,2y)$}\put(58,55){$q_{2\lam_j}(2i\!+\!1,2y)$}
\put(63,19){$q_{2\lam_j}(2y,2i\!+\!1)$}\put(63,19){$q_{2\lam_j}(2y,2i\!+\!1)$}
\put(43,0){$b_{2y,\lam_j,k}$}\put(43,0){$b_{2y,\lam_j,k}$}
\includegraphics[scale=1]{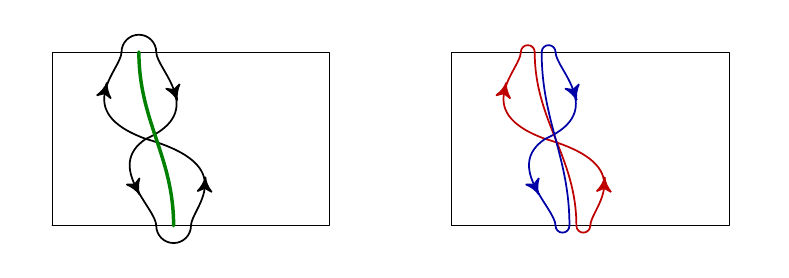}\\[2mm]
&\put(38,75){$t_{2x+1,\lam_j,k}$}\put(38,75){$t_{2x+1,\lam_j,k}$}
\put(58,55){$q_{2\lam_j}(2i\!+\!1,2y)$}\put(58,55){$q_{2\lam_j}(2i\!+\!1,2y)$}
\put(57,0){$b_{2y,\lam_j,k}$}\put(57,0){$b_{2y,\lam_j,k}$}
\includegraphics[scale=1]{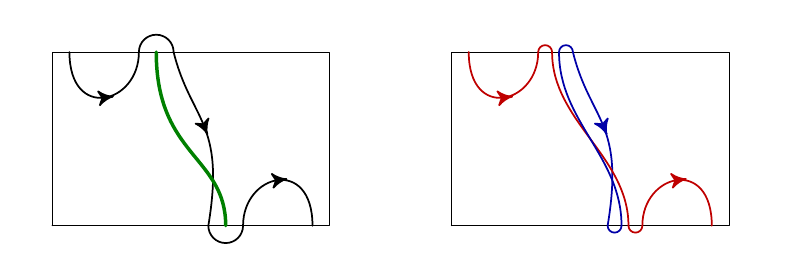}
\end{tabular}
\end{center}
\caption{\label{fig_add_one_edge_one_cycle}The first line represents the case where the edge is inserted in a cycle that is not a principal cycle. The second line represents the case where the edge is inserted in the top principal cycle and the top arc appears before the bottom arc is the left to right path order of the principal cycle}
\end{figure}

 Some cases are not represented in the figure.
The missing cases are: \begin{itemize}
 \item The edge is inserted in the top principal cycle and the bottom arc appear before the top arc is the right to left path order
 \item The edge is inserted in the bottom principal cycle and the bottom arc appear before the top arc is the right to left path order
 \item The edge is inserted in the top principal cycle and the top arc appear before the bottom arc is the right to left path order.
\end{itemize} 
Their proof is nearly identical to the ones represented in the figure and are thus omitted. 
\qed
 \end{pf}

The following lemma give some precision on which arcs of the cycle $\lam_j$ goes to the cycle of length $\frac{q_{2\lam_j}(2x+1,2y)+1}{2}$ and $\frac{q_{2\lam_j}(2y,2x+1)+1}{2}$ respectively. 

\begin{figure}[h!]
\begin{center}
\begin{tabular}{cc}
$\tau$ & $\s=\tau|(1,t_{3,5,1},b_{6,5,1})$\\
\includegraphics[scale=.5]{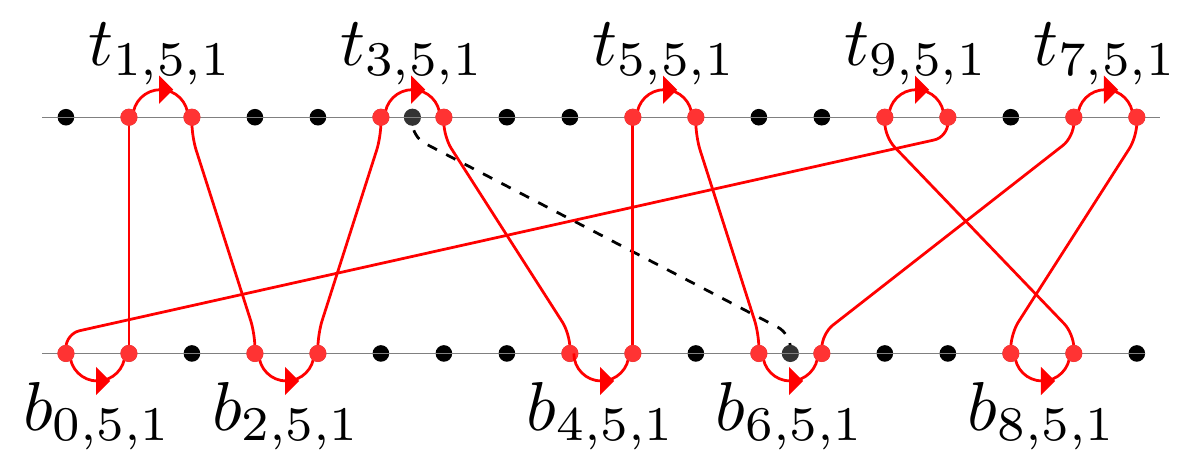} &
\raisebox{12pt}{\includegraphics[scale=.5]{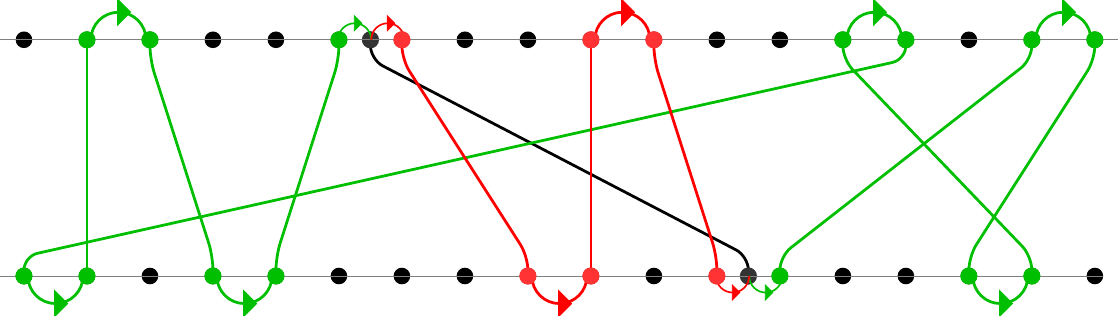}}
\end{tabular}
\end{center}
\caption{\label{figure_ex_appar_arc} The cycle of length 5 is broken into two cycles of length 2 and 4. The two arcs of $\tau$ containing the edge are broken into two arcs in $\s$ and the others are distributed between the two new cycles according to table \ref{table_affi}.}
\end{figure}

\begin{lemma}\label{lem_where_arc2}
Let $\tau$ and $\Pi$ and  $\s=\tau|(1,t_{2x+1,\lam_j,k},b_{2y,\lam_j,k}) $ be as above. Consider the arcs with labels $t_{2u+1,\lam_j,k},b_{2v,\lam_j,k}$ in $\tau$. Those arcs correspond to arcs in $\s$ (one to one except for the arcs $t_{2x+1,\lam_j,k},b_{2y,\lam_j,k}$ in $\tau$ which correspond to two arcs in $\s$ since we inserted the edge rigth within them) and are part of either of the two new cycles.

More precisely we have (the intervals are taken modulo $2\lam_j$ and are oriented):
\begin{table}[h!]
\begin{tabularx}{\textwidth}{|l|X|}
\hline
& arc of $\s$ corresponding to $t_{i,\lam_j,k}$ or $b_{i,\lam_j,k}$ \\ \hline
$i\in ]2x+1,2y[$ & part of the cycle of length $\frac{q_{2\lam_j}(2x+1,2y)+1}{2}$ \\ \hline
$i\in ]2y,2x+1[$ & part of the cycle of length $\frac{q_{2\lam_j}(2y,2x)+1}{2}$ \\ \hline
$i=2x+1$ & two top arcs: the one adjacent left to the edge is part of $\frac{q_{2\lam_j}(2y,2x+1)+1}{2}$ and the other is part of   $\frac{q_{2\lam_j}(2x+1,2y)+1}{2}$ \\ \hline
$i=2y$ & two bottom arcs: the one adjacent left to the edge is part of $\frac{q_{2\lam_j}(2x+1,2y)+1}{2}$ and the other is part of $\frac{q_{2\lam_j}(2y,2x+1)+1}{2}$ \\ \hline
\end{tabularx}
\caption{\label{table_affi}Affiliation of the arcs of $\s$ corresponding to the arcs of the cycle $\lam_j$ of $\tau$.}
\end{table}

Of course all the arcs part of the other cycles are unchanged (thus correspond one-to-one).
See figure \ref{figure_ex_appar_arc}.
\end{lemma}

The two propositions put together tell us that adding an edge to a permutation will either break a cycle into two smaller cycles or merge two cycles into a larger one.

For the sake of convenience we describe the effect of the inverse operation: removing a edge from a permutation.

\begin{proposition}\label{pro_top_adj_arc}
Let $\s$ be a permutation with cycle invariant ($\lam$), $e$ an edge and $\tau$ the permutation with $e$ removed and cycle invariant $(\lam')$. 

If the two arcs (top left and top right) adjacent to the edge $e$ are part of the same cycle of length $\ell$ then $\lam'=\lam\setminus \ell \bigcup \{\ell_1,\ell_2\}$ with $\ell_1+\ell_2+1=\ell$.

More precisely, let $\Pi$ be a consistent labelling of $\s$ and let $t_{1,\ell,1}$, $t_{2x+1,\ell,1}$ be the top left (respectively top right) adjacent arc of $e$, Then $\lam'=\lam\setminus \ell \bigcup \{x-1,\ell-x\}$ and the arcs of  $\ell$ with indices $2,\ldots,2x-1$ correspond to the top/bottom arcs of the cycle of length $x-1$ of $\tau$.\\

\begin{center}
\begin{tabular}{ll}
&$\quad \s,\Pi \ \ \qquad\qquad\qquad\qquad\qquad \tau$\\[1mm]
&$\quad  (\lambda)\quad\qquad \qquad\qquad\qquad\qquad (\lam\setminus \ell \cup \{x-1,\ell-x\})$\\[2.5mm]
&\put(20,71){$t_{1,\ell,1}$}\put(20,71){$t_{1,\ell,1}$}
\put(43,71){$t_{2x\!+\!1,\ell,1}$}\put(43,71){$t_{2x\!+\!1,\ell,1}$}
\put(180,29){$x\!-\!1$}\put(180,29){$x\!-\!1$}
\put(170,49){$\ell\!-\!x$}\put(170,49){$\ell\!-\!x$}
\includegraphics[scale=1]{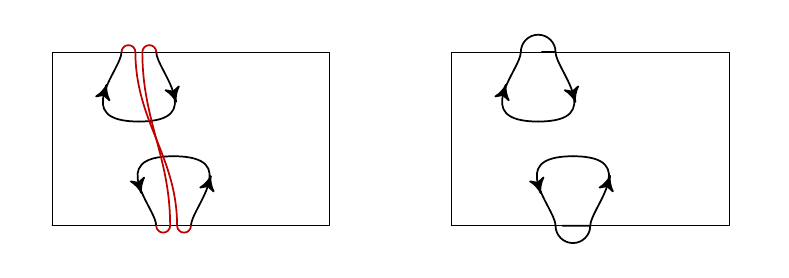}
\end{tabular}
\end{center}

If the two arcs (top left and top right) adjacent to the edge $e$ are part of the two different cycles of length $\ell$ and $\ell'$ respectively then $\lam'=\lam\setminus \{\ell,\ell'\} \bigcup \{\ell+\ell'-1\}$ with $\ell_1+\ell_2+1=\ell$
\end{proposition}

\subsection{Invertibility and invariance of the $\perms l_n$ dynamics}

\begin{proposition}\label{pro_Se_action}
 Let $\s \in \kS^{St}_n$, $e=(i,j)\neq (1,1)$ with neither $i=n$ nor $j=n$.
Let $\tau$ to be the reduction of $(\s,c)$ where the edge $e$ is grayed. Let $(\Pi_b,\Pi_t)$ be consistent labelling and let $t_{2x+1,\ell,k}$ and $b_{2y,\ell',k'}$ be the labels of the arcs containing $e$ in $\tau$ (i.e. $\s=\tau|(1, t_{2x+1,\ell,k},b_{2y,\ell',k'})$). We can have $\ell=\ell'$ and $k=k'$.

Then $S_e(\s)=\tau|(1,t_{2y+1 \!\!\! \mod 2j' ,j',k'},b_{2i+2 \!\!\! \mod 2j,j,k} ))$.
\end{proposition}

\begin{pf}
By definition of $S_e$ and a consistent labelling, it is clear that $S_e(\s)=\tau|(1,t_{2y+1 \!\!\! \mod 2ell' ,\ell',k'},b_{2x+2 \!\!\! \mod 2\ell,\ell,k} ))$.\qed
\end{pf}

\begin{figure}[tb!]
\begin{center}
\begin{tabular}{c}
\includegraphics[scale=1]{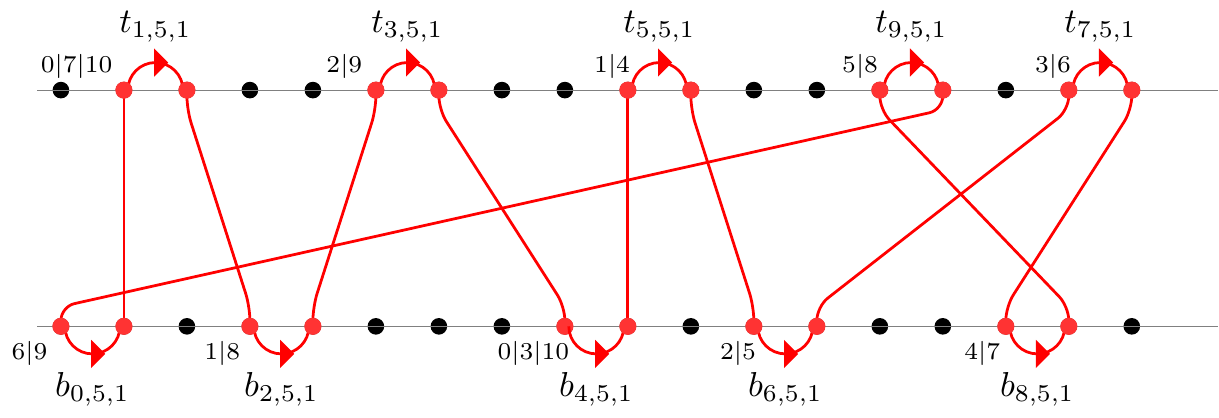}\\
 \includegraphics[scale=1]{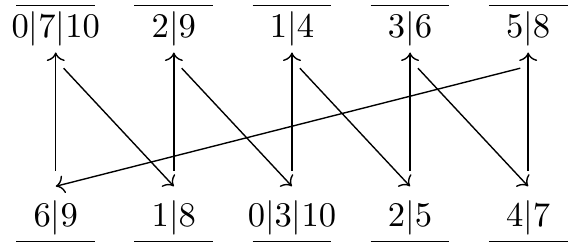}
\end{tabular}
\end{center}
\caption[caption]{\label{fig_one_cycle_apply_dyn}First line: Let $(\tau,\Pi)$ be the represented permutation and $\s=\tau|(1,t_{1,5,1},b_{4,5,1})$. The lists $x|y|\ldots$ left to the arcs represent the position of the two endpoints of the the gray edge of $S_e^i(\s)$ for all $i$.\\ For example, $S_e^0(\s)=\s$ so the edge is within  $t_{1,5,1}$ and $b_{4,5,1}$ corresponding to the pair (0,0), for $S_e(\s)$ the edge is within $b_{2,5,1}$ and $t_{5,5,1}$ corresponding to the pair $(1,1)$ etc... It is clear that after $LCM(2*5,2*5)=10$ iterations we have $S_e^{10}(\s)=\s$, thus the operators are invertibles.\\[1mm]
Second line: We represent  the cycle of length 5 in a schematic way. The numbers represent again the position of the gray edge after $i$ iterations of $S_e$. Note that we follow the natural cyclic order of the cycle rather than order derivated from the position of the arcs (top and bottom) in the permutation.\\ For example, in the permutation, the arc $t_{9,5,1}$ is before $t_{7,5,1}$ in position but not in the cyclic order, thus in the schematic representation the numbers of the fourth top arc correspond to that of $t_{7,5,1}$.}
\end{figure}

\begin{figure}[tb!!]
\begin{center}
\begin{tabular}{c}
\includegraphics[scale=1]{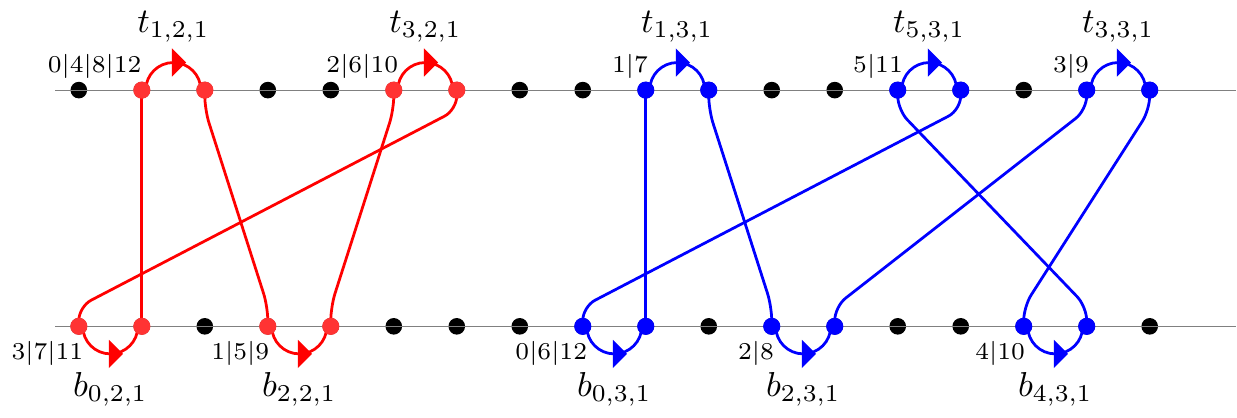}\\
\includegraphics[scale=1]{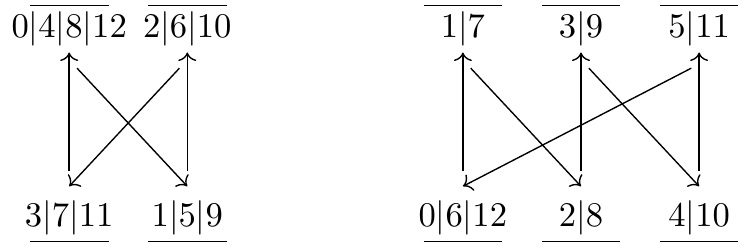}
\end{tabular}
\end{center}
\caption[caption]{\label{fig_two_cycle_apply_dyn}First line: Let $(\tau,\Pi)$ be the represented permutation and $\s=\tau|(1,t_{1,5,1},b_{0,3,1})$. As in figure \ref{fig_one_cycle_apply_dyn}, the lists of numbers represent the position of the gray edge after $i$ iterations of $S_e$. Again, we note that after $LCM(2*3,2*2)=12$ we have $S_e^{12}(\s)=\s$, thus the operators are invertibles.\\
Second line: The schematic representation of the two cycles.}
\end{figure}

In other words, the operator $S_e$ move the edge $e$ along the cycle (or the pair of cycles) in $\tau$ containing it. See figure \ref{fig_one_cycle_apply_dyn} and \ref{fig_two_cycle_apply_dyn} for an illustration of the proposition in the case $\ell=5, \ell'=5,k=k'=1$ and the case $\ell=2,\ell'=3,k=k'=1$.

As we noted below definition \ref{def_coloring}, the notation $\tau|(1,t,b)$ does not make sense if the edge $e=(i,j)$ has $i=n$ or $j=n$, thus to apply the proposition we need the technical condition that this is not the case on $\s$. Thankfully we already know that it will not be the case for  $S_e(\s)$ (cf remark \ref{rk_never_at_end}) thus we can apply the lemma repeatedly and so 
\begin{align}
S_e^i(\s)&=\tau|(1,t_{2y+i \!\!\! \mod 2ell' ,\ell',k'},b_{2x+1+i \!\!\! \mod 2\ell,\ell,k} )) \quad\text{ if $i$ odd}\\
S_e^i(\s)&=\tau|(1,t_{2x+1+i \!\!\! \mod 2ell ,\ell,k},b_{2y+i \!\!\! \mod 2\ell',\ell',k'} )) \quad \text{ if $i$ even}
\end{align}

\begin{corollary}
The $\perms l_n$ dynamics is invertible and leaves invariant the cycle invariant.
\end{corollary}

\begin{pf}
The statements must be proven for the operators $S_e$ only as it was already done in section \ref{sssec.cycinv} (cf article \cite{DS17}) for $L$ and $L'$.

Let $\s \in \kS^{St}_n$ and $e=(i,\s(i))$ be an edge, we can suppose $i\neq n$ and $\s(i)\neq n$ (if not choose $\s'=L'\s$ or $L\s$ or $LL'\s$ depending on whether $i=n$, $\s(i)=n$ or both).

Let $\tau$ be as in proposition \ref{pro_Se_action}, then \[\s=\tau|(1,t_{2x+1,j,k},b_{2y,j',k'}) \text{ and }  S_e(\s)=\tau|(1,t_{2y+1 \!\!\! \mod 2j' ,j',k'},b_{2i+2 \!\!\! \mod 2j,j,k} )).\]
Thus by proposition \ref{pro_one_edge} or \ref{pro_one_edge_two}, $\s$ and $S_e(\s)$ have the same cycle invariant since the edge is added within two cycles of the same length or within one cycle and such that $q_{2j}(2x+1,2y)=q_{2j}(2i+2,2y+1)$ and $q_{2j}(2y,2x+1)=q_{2j}(2y+1,2i+2)$. Moreover $S_e$ is clearly invertible since $S_e^{\lcm(2j,2j')}(\s)=\s$. (see figure \ref{fig_one_cycle_apply_dyn} and \ref{fig_two_cycle_apply_dyn} for an example).
\qed
\end{pf}

\begin{proposition}\label{pro_arf_inv}
We have:\begin{align*}
\Abar\;\Big(
\;
\raisebox{-11pt}{\includegraphics[scale=.4]{../FigFol/fig_sliding_dyn_1.pdf}}
\;
\Big)
&=\Abar\;\Big(
\;
\raisebox{-11pt}{\includegraphics[scale=.4]{../FigFol/fig_sliding_dyn_2.pdf}}
\;
\Big)
\end{align*}
\end{proposition}
\begin{pf}
The proposition enters the framework of theorem 47 of \cite{D18} thus it can be proven automatically.
\qed
\end{pf}

Thus the sign is also invariant since, by remark \ref{rk_general}, if the edge $e$ of $\s$ is not in this configuration then it is for some $L^iL'^j(\s)$ and the sign is invariant by both operator $L$ and $L'$.

%%%%%%%%%%%%%%%%%%%%%%%%%%%%%%%%%%%%%%%%%%%%%%%%%%%%%%%
\section{Proof overview}
\label{ssec_pf_overview}
\label{ch.rauzy2}

In this section we present a proof of the classifcation of the
Rauzy classes of the dynamics $\cS l_n$ by applying the labelling
method.

However we will not need the full extent of the labelling method due to the particularity of the dynamics. 
Let us first start by recalling the labelling method (the full details can de found in \cite{D18} section 2) and then we will explain what we need of it and how we organise the proof in the case of the $\cS l_n$ dynamics.

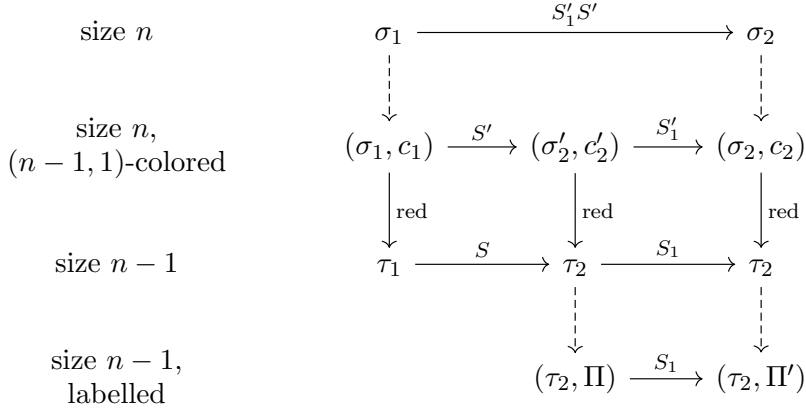
\begin{figure}[tb]
\begin{center}
$\begin{tikzcd}
\textrm{size $n$} 
&
& \s_1 \arrow[d,dashrightarrow]
\arrow{rr}{S'_1 S'} 
&& \s_2  \arrow[d,dashrightarrow]
& 
\\
\makebox[0pt][c]{\raisebox{6.5pt}{\textrm{size $n$,}}}%
\makebox[0pt][c]{\raisebox{-6.5pt}{\textrm{$(n-1,1)$-colored}}}
% \makebox[0pt][c]{\raisebox{-6.5pt}{\textrm{$(n-r,r)$-colored}}}
&
& (\s_1,c_1) \arrow{d}{\textrm{red}} \arrow{r}{S'} 
& (\s'_2,c'_2)  \arrow{d}{\textrm{red}} \arrow{r}{S'_1} 
& (\s_2,c_2) \arrow{d}{\textrm{red}} 
% & x_2  \arrow[swap]{l}{S_2} 
\\
\textrm{size $n-1$}
&
& \tau_1 \arrow{r}{S} 
& \tau_2 \arrow{r}{S_1} 
\arrow[d,dashrightarrow]
& \tau_2 
\arrow[d,dashrightarrow]
\\
\makebox[0pt][c]{\raisebox{6.5pt}{\textrm{size $n-1$,}}}%
\makebox[0pt][c]{\raisebox{-6.5pt}{\textrm{labelled}}}
&
& 
& (\tau_2,\Pi) \arrow{r}{S_1}&
(\tau_2,\Pi')
\\
\end{tikzcd}$
\end{center}
\caption{\label{fig_labelling_method} Outline of the proof of
  connectivity between $\s_1$ and $\s_2$, using the labelling
  method. The sequence $S$ sends $\tau_1$ to $\tau_2$, however the intervals
  containing the gray edge of $\s'_2$ may not be at their correct
  place, in order to match with those of $\s_2$. The sequence $S_1$
  corrects for this. Thus the boosted sequence $S'_1 S'$ sends $\s_1$
  to $\s_2$.}
\end{figure}

\begin{figure}[tb!]
\begin{center}
\includegraphics[scale=1]{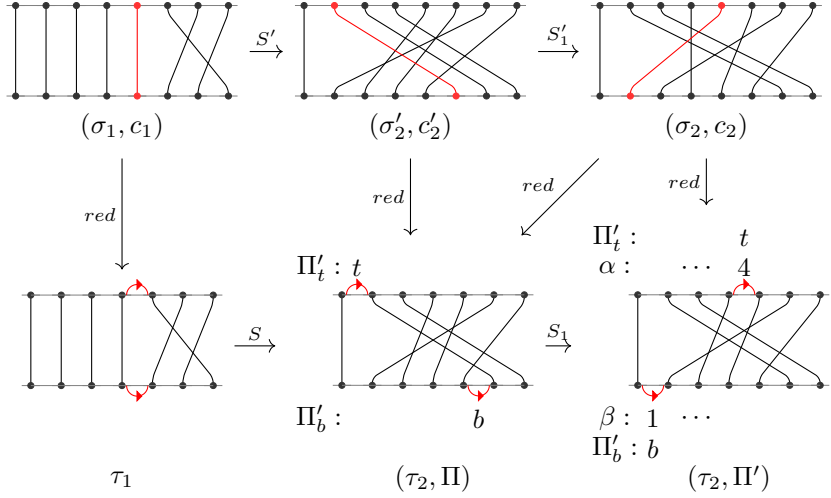}
\end{center}
\caption{\label{fig_labelling_meth_ex}An example of the labelling method. $\s_1$ and $\s_2$ have invariant $(\{7\},+)$ and $\tau_1$ and $\tau_2$ have invariant $(\{2,4\},0)$. In applying the 2-point monodromy, we must find a sequence $S_1$ such that the label $t$ is sent to the position 4 and the label $b$ to the position 1.}
\end{figure}

\paragraph*{The labelling method}$\!\!\!\!\!\!$ is a procedure to prove by induction that two given permutations $\s_1,\s_2$ with the same invariant are connected. 
It proceeds more or less as follows: 
\begin{enumerate}
\item We choose a two coloring $c_1,c_2$ such that $(\s_1,c_1)$, $(\s_2,c_2)$ have both one gray edge and their reduction $\tau_1,\tau_2$ have the same invariant.
\item By induction $\tau_1$ and $\tau_2$ are connected, let $S$ be such that $\tau_2=S(\tau_1).$
\item We now lift the sequence $S$ into a sequence $S'$ such that $(\s_2',c_2')=S'(\s_1,c_1)$ and the reduction of $(\s_2',c_2')$ is $\tau_2$. That is to say the sequence $S'$ complete the following diagram: $\begin{tikzcd}[column sep=15pt] 
(\s_1,c_1) \arrow{d}[swap]{red} \\
\tau_1 \arrow{r}{S}& \tau_2
\end{tikzcd}$ into a commutative square: $\begin{tikzcd}[column sep=15pt] 
(\s_1,c_1) \arrow{d}[swap]{red} \arrow{r}{S'} & (\s_2',c_2') \arrow{d}[swap]{red} \\
\tau_1 \arrow{r}{S}& \tau_2
\end{tikzcd}$.
We call such a sequence a boosted sequence and we show at the beginning of the labelling method that such boosted sequence always exists (we call this the boosted dynamics).
\item We now have $(\s_2,c_2)$ and $(\s_2',c_2')$ that are equal on the set of black edges and thus only differ on the position of the grey edge. The question becomes can we find a sequence that only moves the gray edge of  $(\s_2',c_2')$ to the gray edge of $(\s_2,c_2)$ ?
\item To answer such question we define a labelling $\Pi$ of the intervals between the pair of adjacent vertices (top and bottom) of $\tau_2$ with the following property: if the grey edge of $(\s_2',c_2')$ is within the intervals with label $t$ and $b$ in $(\tau_2,\Pi)$ and if $(\tau_3,\Pi')=S_1(\tau_2,\Pi)$ then the grey edge of $(\s_3,c_3)=S_1'(\s_2',c_2')$ is within the intervals with labels $t$ and $b$ in $(\tau_3,\Pi')$.

In other word the gray edge will follow the labels of the intervals that contains it.

\item With this labelling $\Pi$ we need to resolve the following problem: given two label $t$ and $b$ is there a loop $S_1$ of $\tau_2$ such that $(\tau_2,\Pi')=S_1(\tau_2,\Pi)$ verifies $\Pi'^{-1}(t)=\alpha$ and $\Pi'^{-1}(b)=\beta$ for any $\alpha, \beta$. We call this the two-point monodromy problem.

\item Then we choose $\alpha$ and $\beta$ to be the position of the intervals containing the gray edge of $(\s_2,c_2)$ and if the sequence $S_1$ exists we have $S_1'(\s_2',c_2')=(\s_2,c_2)$. 

Thus we have $\s_2=S_1'S'(\s_1)$. See figure \ref{fig_labelling_method} and \ref{fig_labelling_meth_ex}.
\end{enumerate}

In the case of this article we will stop at step 4 (thus we will not need to define a labelling that tracks the gray edge and solve the 2-point monodromy problem). The reason being that the dynamics itself will already answer the problem. 

Indeed at step 4 we have two configurations $\s_2$ and $\s_2'$ that only differ from one edge and we need to find a sequence of operators that sends the gray edge of $(\s_2',c_2')$ to the position of the gray edge of $(\s_2,c_2)$. In a general dynamics, this is a difficult problem since the operators will modify $\s_2'$ in a non-trivial way and thus both the gray edge $e$ and the remaining black edges will move. 

However this is not the case with our current dynamics, indeed the operator $S_e$ only moves the gray edge. Moreover, it moves the edge in an explicit way in terms of the reduced permutation $\tau_2$. As we have seen in proposition \ref{pro_Se_action} the edge moves along the cycle (or the pair of cycles) of $\tau_2$ that contains it. Thus we can hope that $\s_2=S_e^i(\s_2')$ for some $i$, and the remaining of our discussion will identify a sufficient condition for this to be the case.

We have the following proposition: 
\begin{proposition}\label{pro_connect_1}
Let $(\s,c)$ and $(\s',c')$ be two configurations with same invariants $(\lam,s)$, one gray edge and same reduction $\tau$. Let $\Pi$ be a consistent labelling of $\tau$ and let $(\lam',s')$ be its invariant. 

 If $\s=\tau|(1,t_{x,\lam'_i,k},b_{y, \lam'_i,k})$ and $\s'=\tau|(1,t_{x',\lam'_i,k},b_{y',\lam'_i,k})$ then $\s'=S_e^m(\s)$ for some $m$.
\end{proposition}

In other words, if the gray edge of $(\s_2',c_2')$ and $(\s_2,c_2)$ are within the same cycle of $\tau_2$, we have $\s_2=S_e^i(\s_2')$. Moreover the cycle of $\tau$ is broken into two smaller cycles in $\s_2$ and $\s_2'$ by proposition \ref{pro_one_edge_two}. 

and the proposition:
\begin{proposition}\label{pro_connect_2}
Let $(\s,c)$ and $(\s',c')$ be two configurations with same invariants $(\lam,s\neq 0)$, one gray edge and same reduction $\tau$. Let $\Pi$ be a consistent labelling of $\tau$ and let $(\lam',s')$ be its invariant. Suppose $\lam'$ contains two even cycles of length $(2,2p)$.

 If $\s=\tau|(1,t_{x,\ell_1,1},b_{y,\ell_2,1})$ and $\s'=\tau|(1,t_{x',\ell_i,1},b_{y',\ell_{1-i},1})$ with $(\ell_1,\ell_2)=(2,2p)$ or $(\ell_1,\ell_2)=(2p,2)$ and $i\in \{0,1\}$ then $\s'=S_e^m(\s)$ for some $m$.
\end{proposition}

In other words if the gray edge of $(\s_2',c_2')$ and $(\s_2,c_2)$ has one endpoint within the cycle of length two and the other within the other cycle of even length of $\tau_2$, we have $\s_2=S_e^i(\s_2')$. Moreover the two even cycles of $\tau_2$ are merged into a cycle of length 2p+3 in $\s_2$ and $\s_2'$  by proposition \ref{pro_one_edge} (This is consistent with the fact that $s\neq 0$ since by the characterisation theorem the sign is non zero if and only if the cycle invariant does not contain even cycles).  \\

\noindent These criteria will be our guideline to organise the proof by induction:

By the classification theorem, we know that, for every even $n$, Rauzy classes with cycle invariant a unique cycle (i.e. $\lam=\{n-1\}$) exists. 
By induction we suppose that for a given even $n$ the Rauzy classes with cycle invariant a unique cycle are classified. (as a  reminder there are three such classes: the exceptional class and two Rauzy classes with sign + and - respectively). 

Then we apply the labelling method (up to step 4) and proposition \ref{pro_connect_1} to classify classes with cycle invariant two cycles (i.e $\lam=\{i,n-i\}$ for $1 < i< n$) then three then four etc... The processus is finite since every time we cut a cycle into two smaller cycles (since cycle of length 1 are not allowed refer to the end of section \ref{sssec.cycinv}).

Finally since the Rauzy classes with cycle invariant $\lam=\{2,n-2\}$ have been classified we apply the labelling method (up to step 4) and proposition \ref{pro_connect_2} to classify the classes with size $n+2$ and cycle invariant $\lam'=\{n+1\}$. 

Clearly this induction scheme cover all cases.

\begin{remark}
This induction scheme is very close to the original of Kontsevich-Zorich in \cite{KZ03}. Indeed, in their paper, the Rauzy classes with cycle invariant $\lam=\{n-1\}$ are called the minimal statum for genus $2g=n$. Moreover the first step of our induction repeatedly breaks up a singularity at fixed genus and the second step increases the genus.
\end{remark}

A last hurdle remains to be solved before laying out the full structure of the proof. In order to apply propositions \ref{pro_connect_1} and \ref{pro_connect_2}, both $\s_2$ and $\s_2'$ must have a particular structure i.e. the gray edge must be within the single cycle of a given length or within a cycle of length 2 and another even cycle in $\tau_2$.

More precisely we have the two following statements:

\begin{proposition}\label{pro_structure_1}
Let $\s\in \kS^{st}_n$ be a configuration with n even and invariant $(\lam=\{n-1\},s)$ then there exists $\s'\sim \s$ and a coloring $c'$ of $\s'$ with one gray edge with the following property:

The gray edge has one endpoint within a cycle of length $2$ and the other within a cycle of length $2p>2$ of the reduction $\tau$ and there are no other cycles in $\tau$.  Moreover $\tau$ has invariant $(\{2,2p=n-4\},0)$.
\end{proposition}
and

\begin{proposition}\label{pro_structure_2}
Let $\s\in \kS^{st}_n$ be a configuration with invariant $(\lam\neq\{n-1\},s)$ then there exists $\s'\sim \s$ and a coloring $c'$ of $\s'$ with one gray edge with the following property:

The gray edge is within a cycle of length $k$ of the reduction $\tau$ and there are no other cycle of length $k$ in $\tau$. Moreover $\tau$ has invariant $(\lam \setminus \{\ell,\ell'\}\cup\{k=\ell+\ell'-1\},s')$ where $\ell$ is the largest cycle of $\lam$ and $s'$ and $\ell'$ depend only on $(\lam,s)$.
\end{proposition}

\begin{remark}
We will apply those two propositions before starting the labelling method. i.e have two permutations $\pi_1$ and $\pi_2$ and we must prove that $\pi_1 \sim \pi_2$. We apply the above proposition and end up with $(\s_1c_1)$ and $(\s_2,c_2)$. Then we show that if $(\s_1,c_1)$ has the property so does $(\s_2',c_2)$.
\end{remark}
$\s_1$ and $\s_2$ are normal forms in the framework of the labelling method.

 Thus we can organise this proof as follows 
\begin{itemize}
\item{Section 4.1:} We define the boosted dynamics.
\item{Sections 4.2:} We prove proposition \ref{pro_connect_1} and \ref{pro_connect_2}
\item{Section 4.3:} We construct permutations of some particular forms for every given $(\lam,s)$ that will allow us to prove proposition \ref{pro_structure_1} and \ref{pro_structure_2} during the induction. We also demonstrate a few technical statements that will be useful for their proof.
\item{Section 8:} We proceed with the induction.
\end{itemize}

\section{Preparing the induction}

\subsection{Boosted dynamics}
In this section we define the boosted dynamics for the operators $L,L'$ and $S_e$.

Let $(\s,c)$ be a permutation with one gray edge $e'$ and let $\tau$ be the reduction. Let $S$ be a sequence for $\tau$ and $B(S)$ be the boosted sequence for $(\s,c)$. We define for $(\s,c)$, $B(L),B(L')$ and $B(S_e)$ for every $e$ of $\tau$ and then for a sequence $S=O_k,\ldots O_1$ where $O_i\in \{l,L',(S_e)_e\}$ $B(S)=B(O_k)\ldots B(O_1)$. 

We have \[B(L)(\s,c)=\begin{cases}L^2 & \text{If $e'=(\s^{-1}(n),n)$ in $(\s,c)$ or $e'=(L(\s)^{-1}(n),n)$ in $L(\s,c)$}\\
L & \text{Otherwise.}
\end{cases}\]
and
\[B(L')(\s,c)=\begin{cases}L^{\prime 2} &\text{If $e'=(n,\s(n))$ in $(\s,c)$ or $e'=(n,L'(n))$ in $L'(\s,c)$}\\
L' &\text{Otherwise.}
\end{cases}\]

This definition garantees that the gray edge will never be have an endpoint in the top or bottom right corner (recall that we do not want the gray edge to be in the right corner due to our definition of a consistent labelling).

\[B(S_e)(\s,c)=\begin{cases}S_e S_{e'} &\text{If $e'=(i',j')$ has one endpoint adjacent to the right of $e=(i,j)$.}\\
&\text{  i.e $i'=i+1$ or $i=n$ and $i'=2$ or $j'=j+1$ or $j=n$ and $j'=2$.}\\
 S_e &\text{ Otherwise.}
\end{cases}\]

The idea of this definition is the following: if the edge $e'$ is not adjacent to the right of $e$ then clearly the reduction of $S_e(\s,c)$ is $S_e(\tau)$. Thus if the edge $e'$ adjacent to the right of $e$ we first move it away with $S_{e'}$ and then we apply $S_e$.

\begin{figure}[h!]
\begin{center}
\includegraphics[scale=.5]{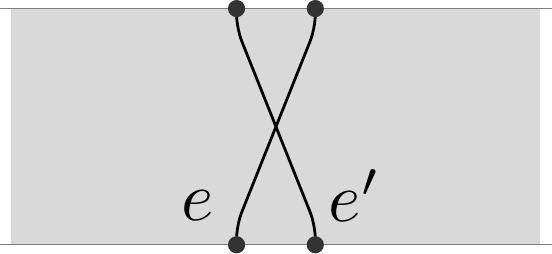}
\end{center}
\caption{The only case in which an edge $e'$ can still be right adjacent to another edge $e$ after application of $S_{e'}$.}
\end{figure}

Thus we only need to justify that in $S_{e'}(\s,c)$ the edge $e'$ is not adjacent to the right of $e$. Let us say wlog that in $(\s,c)$ the edge $e'$ has it bottom endpoint adjacent to the right of $e$ then it is clear that the top endpoint of $e'$ cannot be adjacent to the right of $e$ in $S_{e'}(\s,c)$ since this could only happens if its the bottom endpoint was adjacent to the left of $e$ in $(\s,c)$. Now if the bottom endpoint of $e'$ were adjacent to the right of $e$ in $S_{e'}(\s,c)$ it would mean that the top endpoint of $e'$ is adjacent to the left of $e$ in $(\s,c)$ but then we would have $e=(i,j)$ and $e'=(i+1,j-1)$ in $(\s,c)$ which in turn implies that $\s$ has a cycle of length 1 (cf end of section \ref{sssec.cycinv}). Since we banned this case, this cannot happen.

\subsection{Proof of proposition \ref{pro_connect_1} and \ref{pro_connect_2}}
Let us prove the two propositions, introduced in the proof overview, that allows us to to connect two permutations $(\s,c)$ and $(\s',c)$ having the same reduction $\tau$ by an operator $S_e^i$.

\begin{proof}[Proof of proposition \ref{pro_connect_1}]
Let $(\s,c)$ and $(\s',c')$ two permutations with the same invariant $(\lam,s)$ one gray edge and reduction $\tau$ with invariant $(\lam',s)$ such that $\s=\tau|(1,t_{2x+1,\lam'_i,k},b_{2y, \lam'_i,k})$ and $\s'=\tau|(1,t_{2x'+1,\lam'_i,k},b_{2y',\lam'_i,k})$. 

Then by proposition \ref{pro_one_edge_two} we must have 
\begin{align*}
\lam&=\lam'\setminus \{\lam'_i\} \bigcup\{ \frac{q_{2\lam'_i}(2x+1,2y)+1}{2},\frac{q_{2\lam'_i}(2y,2x+1)+1}{2}\}\\
&=\lam'\setminus \{\lam'_i\} \bigcup\{ \frac{q_{2\lam'_i}(2x'+1,2y')+1}{2},\frac{q_{2\lam'_i}(2y',2x'+1)+1}{2}\} 
\end{align*}
thus \[\left(q_{2\lam'_i}(2x+1,2y), q_{2\lam'_i}(2y,2x+1)\right)=\begin{cases}
\left(q_{2\lam'_i}(2x'+1,2y'), q_{2\lam'_i}(2y',2x'+1)\right)\\
\left(q_{2\lam'_i}(2y',2x'+1), q_{2\lam'_i}(2x'+1,2y')\right) \text{otherwise.}
\end{cases}\]
In the first case we have $S_e^{q_{2\lam'_i}(2x+1,2x'+1)}(\s)=\s'$. Indeed by proposition \ref{pro_Se_action} 
\[
S_e^{q_{2\lam'_i}(2x+1,2x'+1)}(\tau|(1,t_{2x+1,\lam'_i,k},b_{2y, \lam'_i,k}))=\tau|(1,t_{2x+1+q_{2\lam'_i}(2x+1,2x'+1)\!\!\mod 2\lam'_i,\lam'_i,k},b_{2y+q_{2\lam'_i}(2x+1,2x'+1)\!\!\mod 2\lam'_i, \lam'_i,k})
\]
 and we have
\begin{align*}
2x+1+q_{2\lam'_i}(2x+1,2x'+1)\mod 2\lam'_i&=2x'+1\text{$\quad$by definition of $q_{2\lam'_i}(2x+1,2x'+1)$}\\
2y+q_{2\lam'_i}(2x+1,2x'+1)\mod 2\lam'_i 
&= 2x+1 + q_{2\lam'_i}(2x+1,2y) +q_{2\lam'_i}(2x+1,2x'+1)\\
& \ \mod 2\lam'_i  \text{$\quad$by definition of $q_{2\lam'_i}(2x+1,2y)$}\\
& =2x+1+q_{2\lam'_i}(2x+1,2x'+1) + q_{2\lam'_i}(2x'+1,2y') \\
& \ \mod 2\lam'_i\\
& \text{$\quad$$\quad$since $q_{2\lam'_i}(2x+1,2y) =q_{2\lam'_i}(2x'+1,2y')$}\\
&=2x'+1 + q_{2\lam'_i}(2x'+1,2y') \mod 2\lam'_i \\
&=2y'.
\end{align*}

Likewise in the second case we have $S_e^{q_{2\lam'_i}(2x+1,2y')}(\s)=\s'$.
Indeed 
\begin{align*}
2x+1+q_{2\lam'_i}(2x+1,2y')\mod 2\lam'_i& =2y'\\
2y+q_{2\lam'_i}(2x+1,2y')\mod 2\lam'_i &= 2x+1 + q_{2\lam'_i}(2x+1,2y) +q_{2\lam'_i}(2x+1,2y')\mod 2\lam'_i \\
&=2y' + q_{2\lam'_i}(2x+1,2y)   \mod 2\lam'_i \\
&=2y'+q_{2\lam'_i}(2y',2x'+1) \mod 2\lam'_i \\
& \text{$\quad$since $q_{2\lam'_i}(2x+1,2y) =q_{2\lam'_i}(2y',2x'+1)$}\\
&=2x'+1.
\end{align*}

  \end{proof}
before starting the proof of the second proposition, let us recall a proposition from \cite{D18}
\begin{proposition}\label{pro_opposite_sign}
Let $\tau$ be a permutation with two even cycles $2\ell$ and $2\ell'$ and let $\Pi$ be a consistent labelling of $\tau$. 
If $\s= \tau|(1,t_{2x,2\ell,1},b_{2y,2\ell',1})$ and $\s'= \tau|(1,t_{2x,2\ell,i},b_{2y+2,2\ell',j})$ then $s(\s)=-s(\s').$
\end{proposition}

This corresponds (with some minor notational changes) to proposition 53 proven with thanks to the identities of proposition 51 in section 8.3.

\begin{proof}[Proof of proposition \ref{pro_connect_2}]
Let $(\s,c)$ and $(\s',c')$ two permutations with the same invariant $(\lam,s\neq0)$ one gray edge and reduction $\tau$ with invariant $(\lam',s)$ such that $\s=\tau|(1,t_{2x+1,\ell_1,1},b_{2y, \ell_2,1})$ and $\s'=\tau|(1,t_{2x'+1,\ell_i,1},b_{2y',\ell_{1-i},1})$ and $(\ell_1,\ell_2)=(2,2p)$ or $(2p,2)$ and $i\in \{0,1\}$. 

Wlog we can suppose 
\[\s=\tau|(1,t_{2x+1,2p,1},b_{2y, 2,1}) \text{ and } \s'=\tau|(1,t_{2x'+1,2p,1},b_{2y',2,1}),\]
 indeed if $\s=\tau|(1,t_{2x+1,2,1},b_{2y, 2p,1})$ then $S_e(\s)=\tau|(1,t_{2y+1 \mod 4p,2p,1},b_{2x+2\mod 4, 2,1})$, by proposition \ref{pro_Se_action} thus we can always exchange $\s$ with $S_e(\s)$ and $\s'$ with $S_e(\s')$.

Then by proposition \ref{pro_Se_action} we have $S_e^{q_{2p}(2x+1,2x'+1)}(\s)=\s'$. Since \[S_e^{q_{2p}(2x+1,2x'+1)}(\s)=\tau|(1,t_{2x+1+q_{2p}(2x+1,2x'+1)\!\! \mod 4p ,2p,1},b_{2y+q_{2p}(2x+1,2x'+1)\!\! \mod 4,2,1})\] and 
\[
2x+1+q_{2p}(2x+1,2x'+1)\!\!\mod 4p=2x'+1
\]
as well as 
\[
2y+q_{2p}(2x+1,2x'+1) \!\!\mod 4 =\begin{cases} 2y'\\
2y' +2 \mod 4 
\end{cases} \text{$\quad$since $q_{2p}(2x+1,2x'+1)$ is even.} 
\]
Now if $2y+q_{2p}(2x+1,2x'+1) \mod 4=2y'+2 \mod 4$ then $s(S_e^{q_{2p}(2x+1,2x'+1)}(\s))=-s(\s')$ by proposition \ref{pro_opposite_sign} but $s(S_e^{q_{2p}(2x+1,2x'+1)}(\s))=s(\s)$ since the operator $S_e$ leaves the sign invariant and $s(\s)=s(\s')\neq 0$ thus we must have $2y+q_{2p}(2x+1,2x'+1) \mod 4 = 2y'$ and $S_e^{q_{2p}(2x+1,2x'+1)}(\s)=\s'$.
  \end{proof}

\subsection{Technical statements for the proof of propositions \ref{pro_structure_1} and \ref{pro_structure_2} \label{section_tech}}
This section contains all technical propositions that will be needed to prove the two propositions \ref{pro_structure_1} and \ref{pro_structure_2} during the induction. 

We have the following proposition (it corresponds to the proposition 39 page 43 of \cite{DS17} adapted to the extended Rauzy dynamics)

\begin{proposition}[Properties of the standard family]
\label{trivialbutimp}
Let $\s$ be a standard permutation, and 
$S=\{\s^{(i)}\}_i = \{L^i(\s)\}_i$ its L-standard family. The latter has
the following properties:
\begin{enumerate}
  \item Every $\tau \in S$ has $\tau(1)=1$;
  \item The $n-1$ elements of $S$ are all distinct;
  \item Let $m_i$ be the multiplicity of the integer $i$ in $\lambda$
    (i.e.\ the number of cycles of length $i$) and $r$ be the top principal cycle of $\s$. 
There are $i\, m_i$ permutations of $S$ which are of type
    $X(r,i)$, and $1$ permutation of type $H(r-j+1,j)$, for each $1
    \leq j\leq r$.\footnote{Note that, as $\sum_i i\, m_i=n-1$ by
      the dimension formula (\ref{eq.size_inv_cycle}), this list
      exhausts all the permutations of the family.}
  \item Among the permutations of type $X(r,i)$ there is at least one
    $\tau$ with $\tau^{-1}(2) <\tau^{-1}(n)$.
 \end{enumerate}
Likewise define $S'=\{\s^{(i)}\}_i = \{L'^i(\s)\}_i$ its L'-standard family. We have 
\begin{enumerate}
  \item Every $\tau \in S'$ has $\tau(1)=1$;
  \item The $n-1$ elements of $S'$ are all distinct;
  \item Let $m_i$ be the multiplicity of the integer $i$ in $\lambda$
    (i.e.\ the number of cycles of length $i$) and $r$ be the bottom principal cycle of $\s$. 
There are $i\, m_i$ permutations of $S$ which are of type
    $X(i,r)$, and $1$ permutation of type $H(r-j+1,j)$, for each $1
    \leq j\leq r$.
 \end{enumerate}
\end{proposition}

\begin{lemma}\label{lem_top_adja_arc}
Let $\s$ be a permutation with cycle invariant $\lam$ and at least two cycles. For every $\ell\in \lam$ there is an edge $e$ in $\s$ or $L(\s)$ and $\ell' \in \lam$ such that the top left (respectively right) adjacent arc of $e$ is part of the cycle of length $\ell$ (respectively $\ell'$).
\end{lemma}

\begin{pf}
let $\alpha_1,\ldots,\alpha_\ell$ be the positions of the top arcs part of a cycle of length $\ell$, since there are more than one cycle there must be a top arc not part of this cycle, in particular either the positions $\alpha_1+1,\ldots,\alpha_\ell+1$ or $\alpha_1-1,\ldots,\alpha_\ell-1$ must contains one such arc. In the first case, let $\alpha_i+1$ be such arc, then the edge $e$ adjacent to $\alpha_i$ and $\alpha_i+1$ verifies the lemma.
If this is not the case then we must have $\alpha_{i+1}=\alpha_i+1$ and $\alpha_\ell=n-1$, then in $L(\s)$ the first top arc is part of the cycle of length $\ell$ and the second top arc is not thus the edge $e$ adjacent to those two top arcs verifies the lemma. 
\begin{center}\begin{tabular}{cc}
\includegraphics[scale=.5]{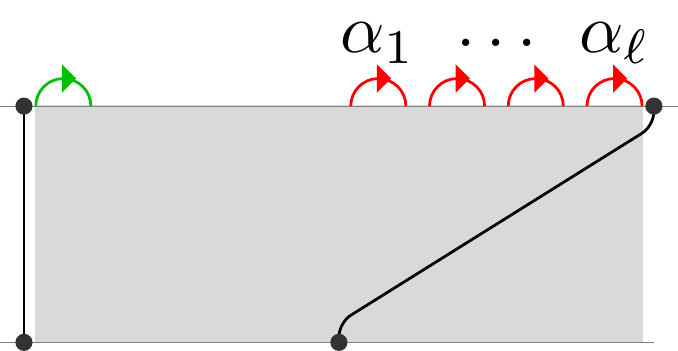} &
\includegraphics[scale=.5]{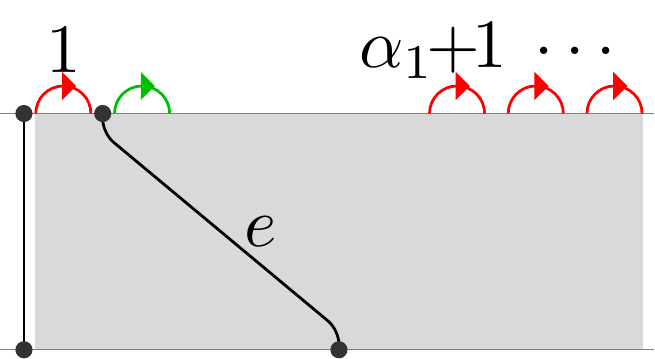}\\
$\s$ & $L(\s)$
\end{tabular}\end{center}
\qed
\end{pf}

Note that the edge $e$ is never (1,1) since it only adjacent to a right arc.

This lemma will be key to prove proposition \ref{pro_structure_1}.
\begin{lemma}\label{lem_choice_left_adj}
Let $C^{st}$ be the non-exceptional class (restricted to the standard permutations) of size $n>6$ with invariants $(\lam=\{\ell,\ell'\},s)$, $\ell>\ell'$. Let $\tau$ be a standard permutation with invariant $(\lam'=\{2,\ell-3,\ell'\},0)$ and $\Pi$ a consistent labelling. Then either $\s=\tau|(1,t_{1,\ell-3,1},b_{0,2,1})$ or $\s'=\tau|(1,t_{1,\ell-3,1},b_{2,2,1})$ is in $C^{st}$.
\end{lemma}
\begin{pf}
First we note that by proposition \ref{pro_one_edge}, $\s$ and $\s'$ have cycle invariant $\lam=\{\ell,\ell'\}$ moreover by proposition \ref{pro_opposite_sign} we have $s(\s)=-s(\s)$. Thus either $s=0$ in which case $\ell$ and $\ell'$ are even and $\s$ and $\s'$ are in $C^{st}$ by the classification theorem \ref{thm.Main_theorem_2} for the extended Rauzy dynamics, or $s\neq 0$ and either $s(\s)$ or $s(\s')$ egal $s$. Let us say $s(\s)=s$ then $\s$ is in $C^{st}$ by the classification theorem \ref{thm.Main_theorem_2}. \qed
\end{pf}

The two next lemmas allow us to ignore the problematic existence of the exceptional class in the induction by showing that when removing an edge we can always avoid falling into the exceptional class.

\begin{lemma}\label{lem_id_1}
let $\s$ be a permutation of size $n$ with cycle invariant $\lam=\{\ell,n-\ell-1\}$ not in the exceptional class, there exists $\s'$ connected to $\s$ and a coloring $c'$ where one edge $e\neq (1,1)$ is grayed such that the reduction $\tau$ of $(\s',c')$ is not in the exceptional class and has cycle invariant $\lam'=\{n-2\}$.
\end{lemma}

\begin{pf}
By lemma \ref{lem_top_adja_arc} there exists $e$ such that the two top adjacent arcs are part of $\ell$ and $n-\ell-1$ respectively so the reduction $\tau$ of $(\s,c)$ where $e$ is grayed has cycle invariant  $\lam'=\{n-2\}$ by proposition \ref{pro_top_adj_arc}.
If $\tau$ is not in the exceptional class then we are done. If it is then $L^i(\tau)=id_{n-1}$ let $\s'=B(L^I)(\s)$ where $B(L^i)$ is the boosted sequence then $\s'$ has the form:
\begin{center}
\begin{tabular}{ccc}
\raisebox{-15pt}{\includegraphics[scale=.5]{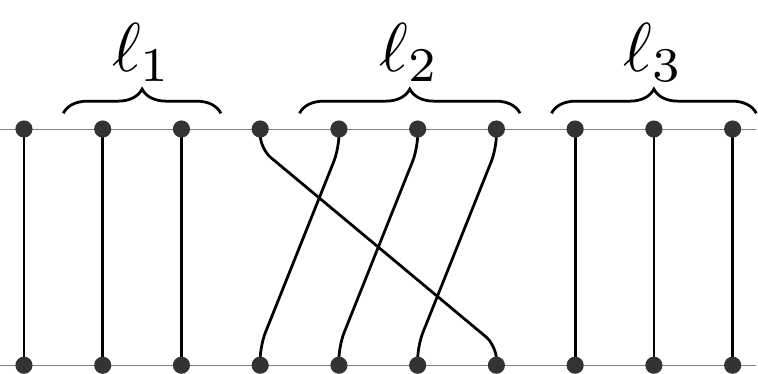}} & or &
\raisebox{-15pt}{\includegraphics[scale=.5]{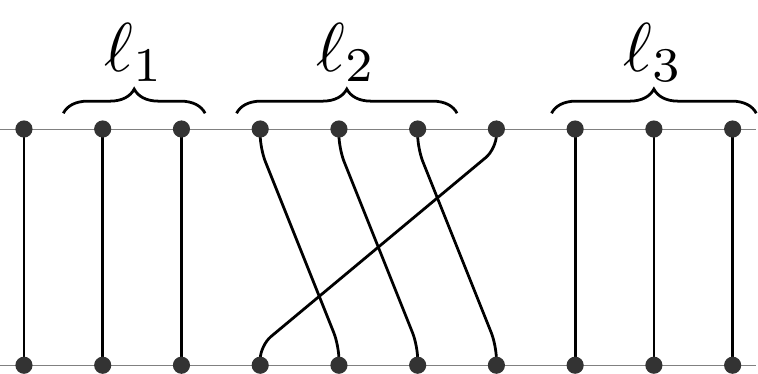}}
\end{tabular}
\end{center}
with $\ell_1+\ell_3\geq1$ since otherwise $L'(\s')=id_{n-1}$ or $L(\s')=id_{n-1}$ respectively, moreover $\ell_2\geq 2$ otherwise there is a cycle of length 1. Now it is clear by inspection of the cycle invariant that we can choose another edge $e'$ such that the reduction $\tau'$ has invariant $\lam=\{n-2\}$ and is not in an exceptional class since $e$ remains (and $\ell_2>1$ so we are not left with the identity).
\end{pf}

\begin{lemma}\label{lem_id_2}
let $\s$ be a permutation with cycle invariant $\lam=\{n-1\}$ not in the exceptional class, there exists an coloring $c$ where one edge $e\neq (1,1)$ is grayed such that the reduction $\tau$ of $(\s,c)$ is not in the exceptional class.
\end{lemma}

\begin{pf}
Almost identical to the one above with the simplification that any edge has its two top adjacent arcs part of the same cycle since there is only one cycle.
\end{pf}

The following definition and two propositions will be essential to prove proposition \ref{pro_structure_2}

\begin{definition}[$T_i$ operator and $TS_i$ operators]
Let $\s$ be a permutation and $e=(i,\s(i))$ an edge then $T_i(\s) $is the following permutation.
Likewise $TS_i$ is defined by a symmetry of $T_i$ and $TS_i(S)$ is the following permutation.
\begin{center}\begin{tabular}{ccc}
\raisebox{5pt}{\includegraphics[scale=.5]{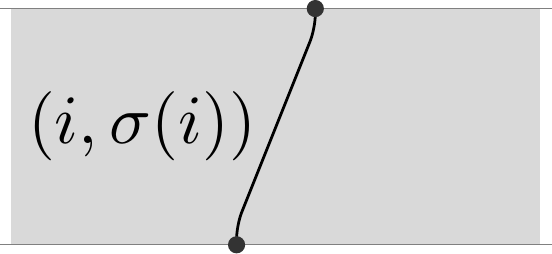}} &
\includegraphics[scale=.5]{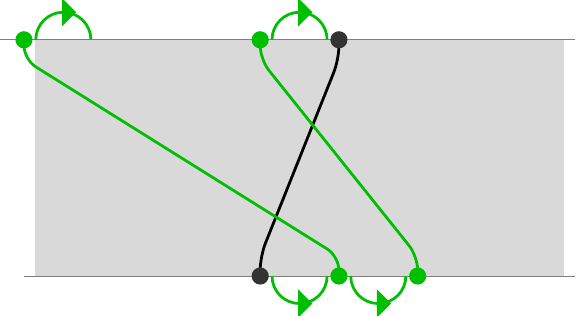}&
\includegraphics[scale=.5]{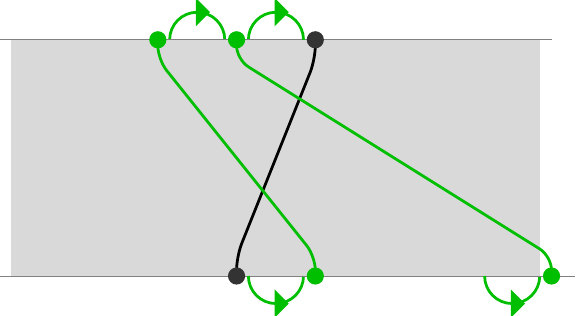}\\
$\s$ & $T_i(\s)$ & $TS_i(\s)$
\end{tabular}\end{center}
\end{definition}
$T_i$ and $TS_i$ have a clear action on the invariant: they do not change the sign and increase by two the length of the top principal cycle (either at the beginning for $T_i$ or at the end for $TS_i$). 
\begin{proposition}\label{pro_T_action}
Let $\s$ be a permutation with invariant $(\lam,s)$ and top principal cycle $r$ and type $X(r,i)$. Then $TS_i(\s)$ and $T_i(\s)$ have invariant $(\lam\setminus\{r\}\bigcup\{r+2\},s)$ and type $X(r+2,i)$.
\end{proposition}
For a proof, the type part is obvious and the invariant part is done in lemma 5.4 of \cite{DS17} (this concerns the case $T_1$, then we obtain the case $T_i$ by applying corollary 3.14 since $T$ is a square constructor). For $TS_i$ it suffices to note that the operator is symmetric to $T_i$ and the invariants are invariant by this symmetry).

\begin{proposition}\label{pro_T_end_perm}
Let $C$ be a class (restricted to standard permutations) with invariant $(\lam,s)$, for every $\ell,\ell'\in \lam$ with $\ell>2$ there exist $\s \in C$ and a consistent labelling $\Pi$ of $\s$ with the following form:
\[(\s,\Pi)= \raisebox{-35pt}{\includegraphics[scale=.5]{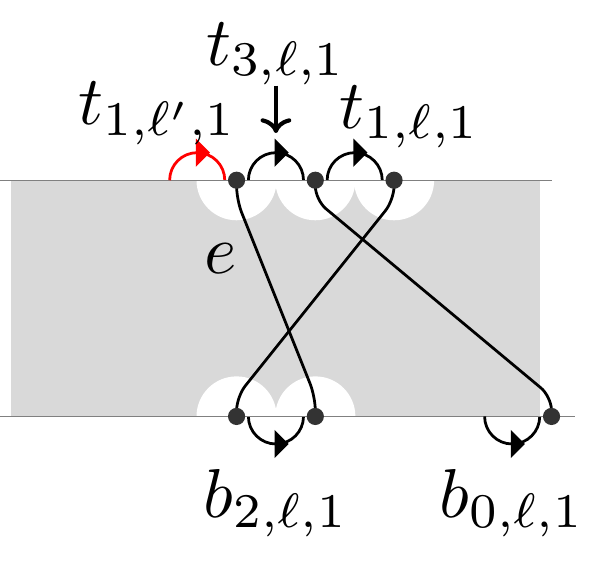}}\]
\end{proposition}
\begin{pf}
Let $C'$ be the class with invariant $(\lam \setminus \{\ell\} \bigcup\{\ell-2\},s)$, let $\tau$ be a standard permutation of type $X(\ell-2,i)$ for some $i$ (such permutation exists by the proposition \ref{trivialbutimp}: start with any standard permutation and consider the $L'$-standard family). By lemma \ref{lem_choice_left_adj} either $\tau$ or $L(\tau)$ (let us say $\tau$) has an edge $e=(i,\tau(i))$ such that its top left adjacent arc is part of a cycle of length $m'$. Then $TS_i(\tau)$ verifies the conditions of the proposition and is in $C$ by the classification theorem \ref{thm.Main_theorem_2} and proposition \ref{pro_T_action}
\qed
\end{pf}

The following lemma is key to get the correct sign invariant in proposition \ref{pro_structure_2}. 

\begin{lemma}\label{lem_sign_control}
Let $\s$ be a permutation with invariant $(\lam,s)$, let c be a coloring with one gray edge $e=(i,\s(i))$ such that the reduction $\tau$ has invariant $(\lam'=\lam\setminus\{\ell,\ell'\}\cup\{\ell+\ell'-1\},s')$. 
Then:

If $\lam'$ contains even cycles then the sign of $\tau$ is s'=0. 

If $\lam'$ does not contains even cycles then the sign of $\tau$ is $s'=s$.
\end{lemma}

\begin{pf}
Clearly if $\lam'$ contains even cycles we know by the classification theorem \ref{thm.Main_theorem_2} that the sign must be 0. 

If $\lam'$ does not contains cycles of length 2, then since $\lam=\lam'\cup\{\ell,\ell'\}\setminus \{\ell+\ell'-1\}$ and that there are always an even number of even cycles (cf theorem \ref{thm.Main_theorem_2}) $\lam$ must not contains even cycles either and its sign $s$ is non zero.

Let us consider $(\s',c')=L'^{n-i}(\s,c)$ then in $\s'$, $e=(n,\s'(n)=\s(i))$ and the reduction of $(\s',c')$ is $\tau'=L^{n-i}(\tau)$ thus $\tau'$ has also invariant $(\lam',s')$.

\begin{center}
\begin{tikzcd}
\raisebox{-12pt}{\includegraphics[scale=.5]{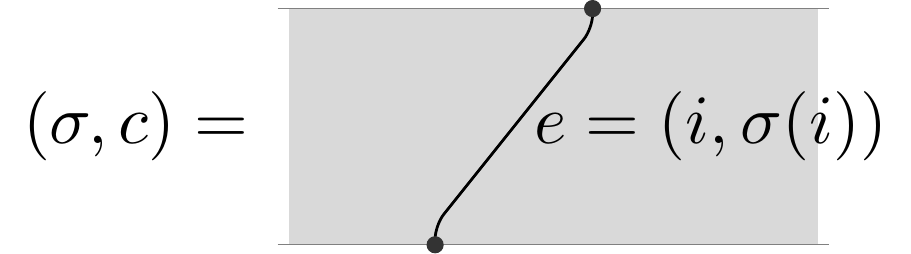}} \arrow[d, "red"] \arrow[r, "L^{n-i}"]&\raisebox{-12pt}{\includegraphics[scale=.5]{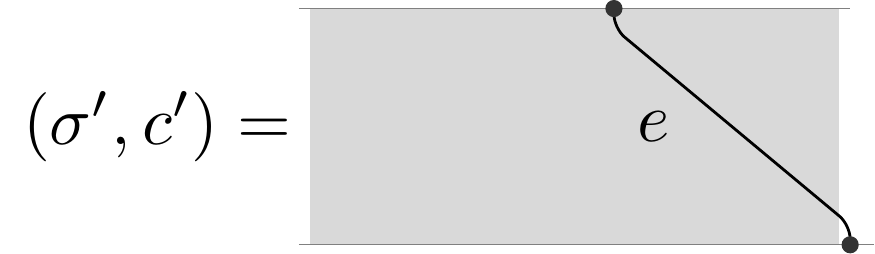}}\arrow[d, "red"]\\
\raisebox{-12pt}{\includegraphics[scale=.5]{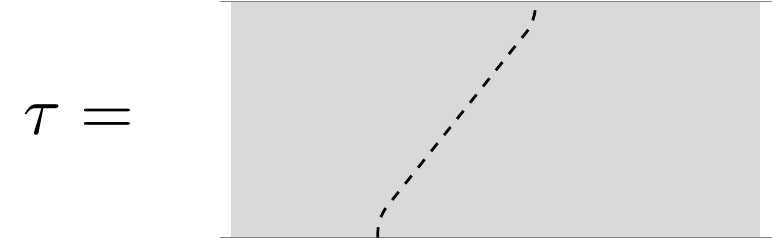}}\arrow[r, "L^{n-i}"] &\raisebox{-12pt}{\includegraphics[scale=.5]{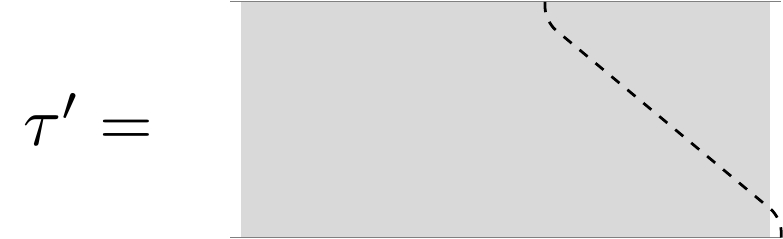}}
\end{tikzcd}\end{center}

Proposition 55 of \cite{D18} tell us that if $\s_1$ has invariant $(\lam,s)$, $\tau_1$ has invariant $(\lam',s')$ and is the reduction of $(\s_1,c_1)$ where the gray edge is $e_1=(\s_1^{-1}(1),1)$ then $s'=0$ or $s'=s$. (Indeed it says that if both endpoint of the edge $e_1$ are inserted within the top principal cycle of $\tau_1$ and its top endpoint is at position 1 then we have the statement. Thus by proposition \ref{pro_one_edge_two} this is equivalent to the fact that $\lam=\lam'\cup\{\ell,\ell'\}\setminus\{\ell+\ell'-1\}$ since we cut the top principal cycle (of length $\ell+\ell'-1$) of $\tau_1$ in two.)

Since the both the cycle invariant and the sign are invariant by the rotation of 180 degrees taking $\s_1=rot_{\pi}(\s')$ and $\tau_1=rot_\pi(\tau')$ implies that $s'=0$ or $s'=s$ since $s'\neq 0$ we have $s'=s$.

\begin{center}
\begin{tikzcd}
\raisebox{-12pt}{\includegraphics[scale=.5]{figure/fig_Lni_4.pdf}} \arrow[r, "rot_\pi"]&\raisebox{-12pt}{\includegraphics[scale=.5]{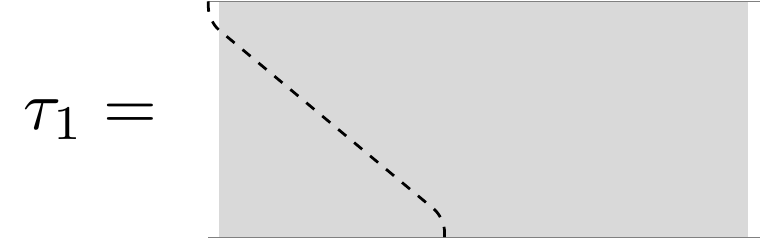}}\\[-20pt]
\qquad \quad(\lam',s') &\qquad\quad (\lam',s')
\end{tikzcd}\end{center}
\qed
   \end{pf}

\section{The induction\label{sec_induction}}

\noindent Let us proceed with the induction. The classification theorem is true at small size $<10$ (as verified by Boissy in \cite{BoiCM}).

\noindent \textbf{Inductive case:} Let $g(\s)=n-\ell$ where $n$ is the size of $\s$ and $\ell$ the number of cycle, $g$ is always odd (and is related to the genus of the translation surface associated to the permutation\footnote{more precisely the genus is $(n-\ell+1)/2$}). 

The induction is on both $g$ and $n$. More precisely we suppose the classification theorem true for all $g'<g$ and for $g$ it is true for all $g-1 \leq n'< n$ then we prove that it is true at size $n$. Finally we prove that it is true for $g+2$ and $n=g+1$.

Indeed recall the proof overview, the induction is in two steps: we suppose the theorem true for $n$ even with $\lam=\{n-1\}$, then we prove the theorem for $n+1$ and two cycles $n+2$ and three cycles etc. Clearly this correspond to $g=n+1$ fixed and $n$ increasing. Finally we prove the theorem for $n+2$ and $\lam=\{n+1\}$ this corresponds to $g+2$ and $n'=g+1=n+2$.

Let us first prove the proposition \ref{pro_structure_1} and \ref{pro_structure_2}.

\begin{figure}[tb!]
\begin{center}
\begin{tikzcd}[column sep=.7em,every label/.append style = {font = \normalsize}]
\defineframenode{$(\s,c)$\\ $(\lam=\{n-1\},s)$\\edge: $e$} \arrow[d,"red"] \arrow[r,"S'"] 
&\defineframenode{$(\s',c')$\\ $(\lam,s)$\\edge: $e$}   \arrow[d,"red"] \arrow[r,"S_e^{2x}"] &\defineframenode{$(\s'',c_{e})$\\ $(\lam,s)$\\edge: $e$}   \arrow[ld,"red"]  \arrow[<->]{r}{change}[swap]{coloring}  
&\defineframenode{$(\s'',c_{e'})$\\ $(\lam,s)$\\edge: $e'$}   \arrow[d,"red"]\\
\defineframenode{$\tau$\\ $(\lam'=\{\ell,n\!-\!\ell\!-\!2\},s')$}  \arrow[r,"S"] 
&\defineframenode{$\tau'$\\ $(\lam',s')$}  \arrow[dr,"\text{remove $e'$}"] 
&
&\defineframenode{$\tau''$\\ $(\{2,n\!-\!4\},0)$}\arrow[dl,swap,"\text{remove $e$}"] \\
&&\defineframenode{$\pi$\\ $(\lam''=\{2,\ell\!-\!3,n\!-\!\ell\!-\!2\},0)$} &
\end{tikzcd}
\end{center}
\caption[caption]{\label{fig_2_nminus2_sequence} The configurations involved in the proof of proposition \ref{pro_structure_1}.}
\end{figure}
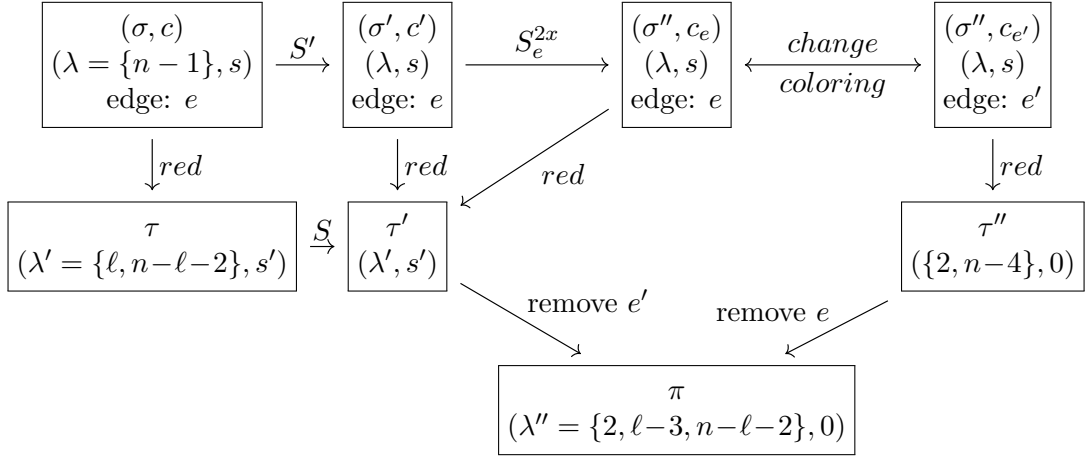

\begin{figure}[tb!]
\begin{center}
\begin{tikzcd}[column sep=-2em,every label/.append style = {font = \normalsize}]
&\rule{-30pt}{0pt} \fbox{\includegraphics[scale=.5]{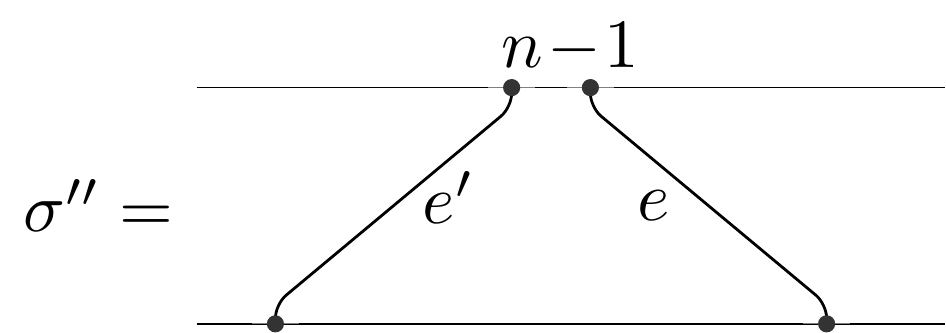} } \arrow[swap,pos=.68]{dl}{\text{red for $(\s'',c_e)$}} \arrow[pos=0.78]{dr}{\text{red for $(\s'',c_{e'})$}}&[-25pt]\\[5pt]
 \fbox{\includegraphics[scale=.5]{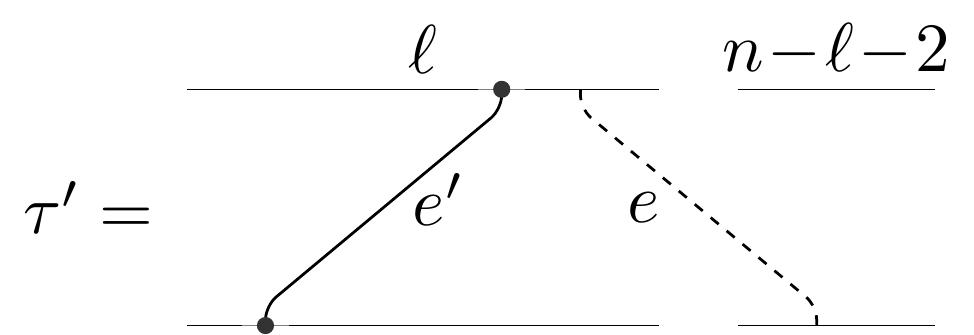}} \arrow{dr} && \fbox{\includegraphics[scale=.5]{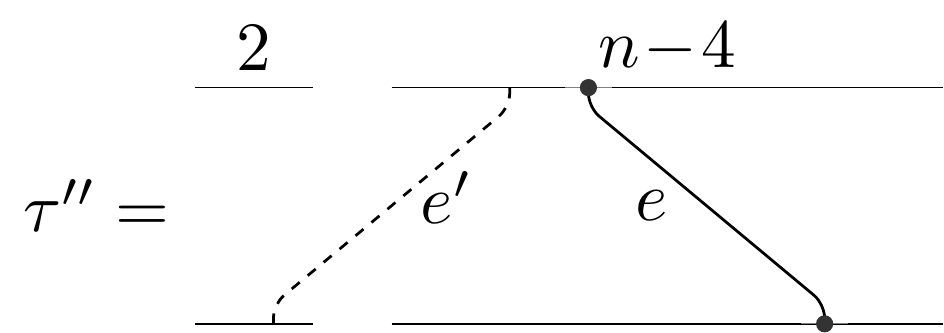}}\arrow{dl}\\[4pt]
& \fbox{\includegraphics[scale=.5]{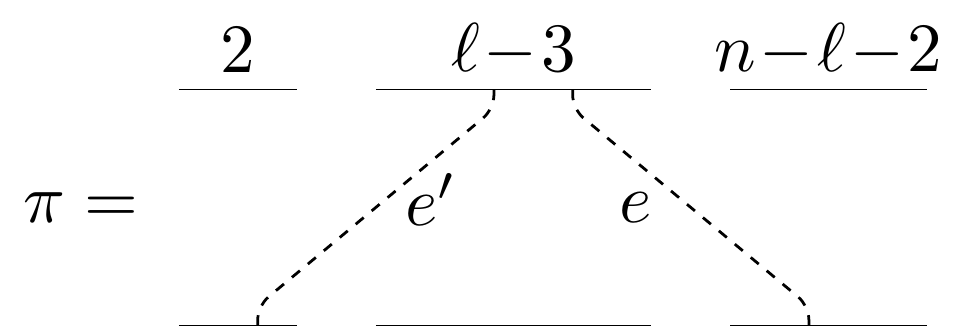}}&
\end{tikzcd}
\end{center}
\caption[caption]{\label{fig_2_nminus2} In the figure the gray edges are represented by dashed edges.\\
Third line: The permutation $\pi$ with cycle invariant $(\{2,\ell-3,n-\ell-2\})$, the gray edge $e'$ is within the cycles of length 2 and $\ell-3$ and the gray edge $e$ is within the cycles of lengths $\ell-3$ and $n-\ell-2$.\\
Second line, left: The permutation $\tau'$ with invariant $(\{\ell,n-\ell-2\})$, the gray edge $e$ is within the cycles of length $\ell$ and $n-\ell-2$.\\
 Second line, right: The permutation $\tau''$ with invariant $(\{2,n-4\})$, the gray edge $e'$ is within the cycles of length $2$ and $n-4$.\\
First line: The permutation $\s''$. Clearly as wanted the edge $e'$}
\end{figure}

\begin{proof}[Proof of proposition \ref{pro_structure_1}]
Let $\s$ be a permutation with invariant $(\lam=\{n-1\},s)$ we will show that there exists $\s''$ connected to $\s$ and a coloring $c_{e'}$ where one edge is grayed such that the reduction $\tau''$ has invariant $(\{2,n-4\},0)$. (For the sign this is obvious since by theorem \ref{thm.Main_theorem_2} if there are even cycles the sign must be 0). Then by proposition \ref{pro_one_edge} the gray edge must be within the two cycles of length 2 and $n-4$ in $\tau''$ and the proposition is proven.

By lemma \ref{lem_id_2}, let $c$ be a coloring of $\s$ with one gray edge $e$ such that the reduction $\tau$ is not in the exceptional class. Since there is only one cycle, the two top adjacent arcs of $e$ are part of the same cycle, thus by proposition \ref{pro_top_adj_arc}, $\tau$ has invariant $\lam'=\{\ell,n-\ell-2\}$. If $\ell=2$ or $n-\ell-2=2$ we are done. otherwise let us suppose that $\ell\geq n-\ell-2$. 

By proposition \ref{lem_choice_left_adj}, there exists $\tau'$ with the same invariant, $\pi$ with cycle invariant $\lam''=\{2, \ell-3, n-\ell-2\}$ and $\Pi$ a consistent labelling of $\pi$ such that $\tau'=\pi|(1,t_{1,\ell-3,1},b_{0,2,1})$, let $e'$ be this edge of $\tau'$. 

By induction, the classification theorem is proven at size $n-1$ thus there exists $S$ such that $\tau'=S(\tau)$. Thus there exists a boosted sequence $S'$ such that $(\s',c')=S'(\s,c)$ where $\tau'$ is the reduction of  $(\s',c')$.

 Define $\Pi'$ the consistent labelling of $\tau'$ such that $\s'= \tau|(1,t_{1,\ell,1},b_{0,n-\ell-2,1})$ (the edge $e$ is within those two cycles due to proposition \ref{pro_one_edge} and the cycle invariant of $\s'$ and $\tau'$).

Now, by lemma \ref{lem_where_arc}, the two top arcs of the cycle of length 2 of $\pi$ correspond to three top arcs of the the cycle of length $\ell$ of $\tau'$, since $\ell>3$ there exists a top arc of the cycle (say $t_{2x+1,\ell,1}$) which does not correspond to an arc of the cycle of length 2 but rather to an arc of the cycle of length $\ell-3$ of $\pi$. Then the gray edge $e$ of $(\s'',c_e)=S^{2x}_e(\s',c')$ is within the arcs labelled $t_{2x+1,\ell,1}, b_{2x \mod ,n-\ell-2,1}$ of $(\tau',\Pi')$ by proposition \ref{pro_Se_action}. Thus in $\pi$ the edge $e$ is within the cycle of length $n-\ell-2$ and the cycle $\ell-3$ by choice of the arc $t_{2x+1,\ell,1}$.

Finally let $c_{e'}$ be the coloring of $\s''$ where the edge $e'$ is grayed and $\tau''$ be the reduction of $(\s'',c_{e'})$ then $\s''$, $\tau''$ and $e'$ satisfy the conditions of the proposition since $\tau''$ has cycle invariant $(\{2,n-4\})$. Indeed, $\tau'$' is obtained from $\pi$ by inserting the edge $e$ within the cycle of length $n-\ell-2$ and the cycle $\ell-3$ thus by proposition \ref{pro_one_edge} since $\pi$ has cycle invariant $\{2,\ell-3,n-\ell-2\}$, $\tau''$ must have cycle invariant $(\{2,n-4\})$.

See figure \ref{fig_2_nminus2_sequence} and \ref{fig_2_nminus2}.
\end{proof}

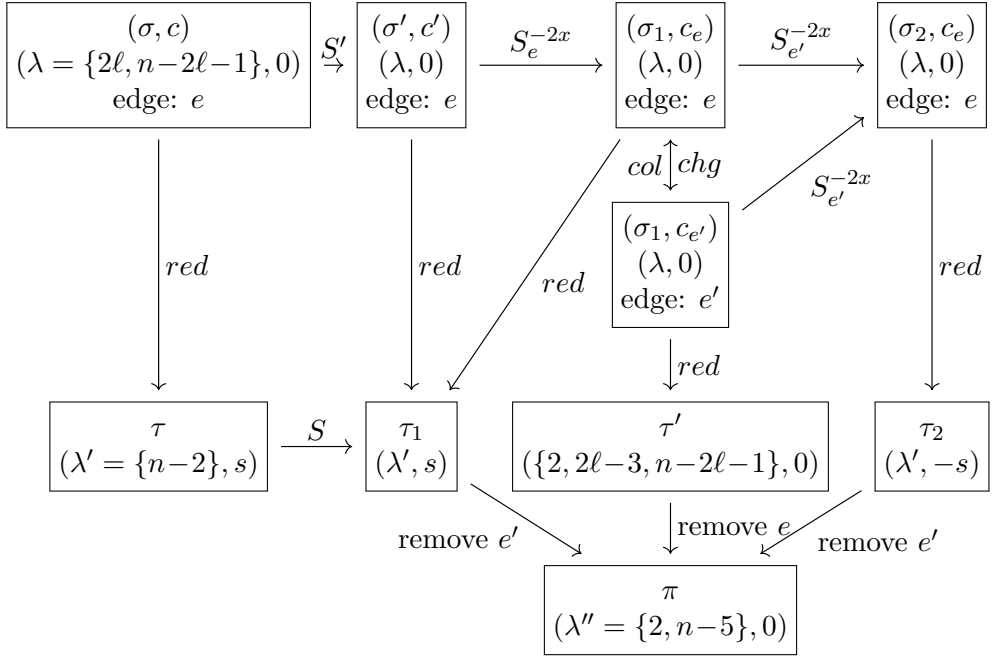
\begin{figure}[tb!]
\begin{center}
\begin{tikzcd}[column sep=.7em,every label/.append style = {font = \normalsize}]
\defineframenode{$(\s,c)$\\ $(\lam=\{2\ell,n\!-\!2\ell\!-\!1\},0)$\\edge: $e$} \arrow[dd,"red"] \arrow[r,"S'"] 
&\defineframenode{$(\s',c')$\\ $(\lam,0)$\\edge: $e$}   \arrow[dd,"red"] \arrow[r,"S_e^{-2x}"] &\defineframenode{$(\s_1,c_{e})$\\ $(\lam,0)$\\edge: $e$}   \arrow[ldd,"red"]  \arrow[<->]{d}{chg}[swap]{col}   \arrow[r,"S_{e'}^{-2x}"]
&\defineframenode{$(\s_2,c_{e})$\\ $(\lam,0)$\\edge: $e$}   \arrow[dd,"red"]\\

&&\defineframenode{$(\s_1,c_{e'})$\\ $(\lam,0)$\\edge: $e'$}   \arrow[d,"red"]  \arrow[ur,swap,"S_{e'}^{-2x}"] &\\
\defineframenode{$\tau$\\ $(\lam'=\{n\!-\!2\},s)$}  \arrow[r,"S"] 
&\defineframenode{$\tau_1$\\ $(\lam',s)$}  \arrow[dr,swap,"\text{remove $e'$}"] 
&\defineframenode{$\tau'$\\ $(\{2,2\ell\!-\!3,n\!-\!2\ell\!-\!1\},0)$}  \arrow[d,"\text{remove $e$}"] 
&\defineframenode{$\tau_2$\\ $(\lam',-s)$}\arrow[dl,"\text{remove $e'$}"] \\
&&\defineframenode{$\pi$\\ $(\lam''=\{2,n\!-\!5\},0)$} &
\end{tikzcd}
\end{center}
\caption[caption]{\label{fig_first_case_sequence} The configurations involved in the proof of proposition \ref{pro_structure_2}, first case.}
\end{figure}

\begin{proof}[Proof of proposition \ref{pro_structure_2}]

Let $\s$ be a permutation of size $n$ and cycle invariant $\lam\neq \{n-1\}$. 
We distinguish two cases: either $\lam=\{2\ell,n-2\ell-1\}$ for some $\ell$ or not. 

\paragraph*{First case: $\lam=\{2\ell,n-2\ell-1\}$.} We suppose $2\ell\geq n-2-\ell$.
Clearly by lemma \ref{lem_top_adja_arc} there exists an edge $e$ of $\s$ such that its two top adjacent arcs are part of the cycles of length $2\ell$ and $n-2\ell-1$. Then the reduction $\tau$ has a unique cycle of length $n-2$ thus $\tau$ and $e$ verify the condition of proposition \ref{pro_structure_2}... except for the sign. Indeed the sign $s$ of $\tau$ is non-zero (by the classification theorem \ref{thm.Main_theorem_2}) and we have no control whether it is positive or negative (in the second case we will use lemma \ref{lem_sign_control} to control the sign of $\tau$ from the sign of $\s$ but here $\s$ has sign 0 and $\tau$ has sign non-zero so the lemma cannot apply). 

Thus all the difficulty of this case will be to get a given sign independant of $\s$. We achieve such goal by constructing two permutations $\s_1$ and $\s_2$ connected to $\s$ which differ only in the position of one edge $e'$, then the reductions where the edge e is grayed $\tau_1$ and $\tau_2$ have respectively invariant $(\{n-2\},s)$ and $(\{n-2\},-s)$ thanks to proposition \ref{pro_opposite_sign}. Then we can choose that of the two that has positive sign. 

See figure \ref{fig_first_case_sequence} for the structure of the proof and figure \ref{fig_perm_constructed} for the form of the constructed permutations.\\
\begin{figure}[bt!]
\begin{center}\begin{tabular}{ll}
\framebox{$\s_1 : \ $ \raisebox{-20pt}{\includegraphics[scale=.5]{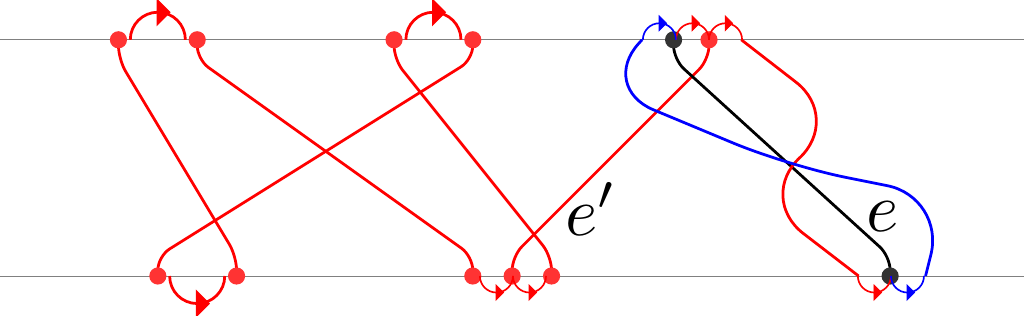}}}
& \framebox{$\s_2 : \ $ \raisebox{-20pt}{\includegraphics[scale=.5]{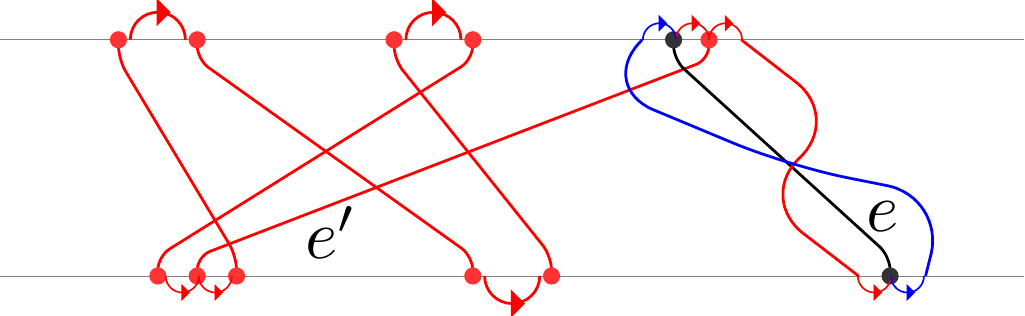}}} \\[1mm]
\framebox{$\tau_1 : \ $ \raisebox{-20pt}{\includegraphics[scale=.5]{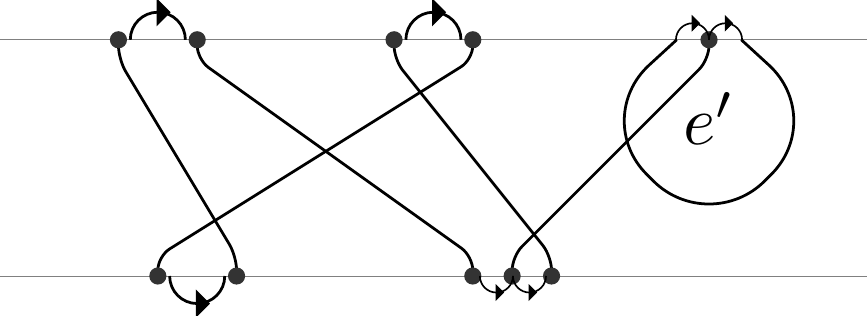}}}
& \framebox{$\tau_2 : \ $ \raisebox{-20pt}{\includegraphics[scale=.5]{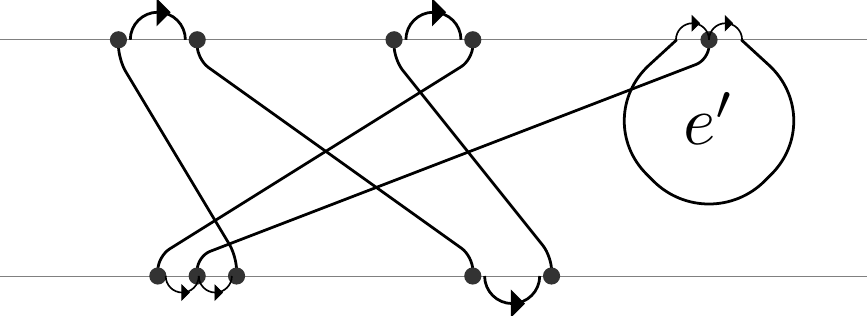}}} \\[2pt]
\multicolumn{2}{c}{\framebox{$\pi : \ $ \raisebox{-20pt}{\includegraphics[scale=.5]{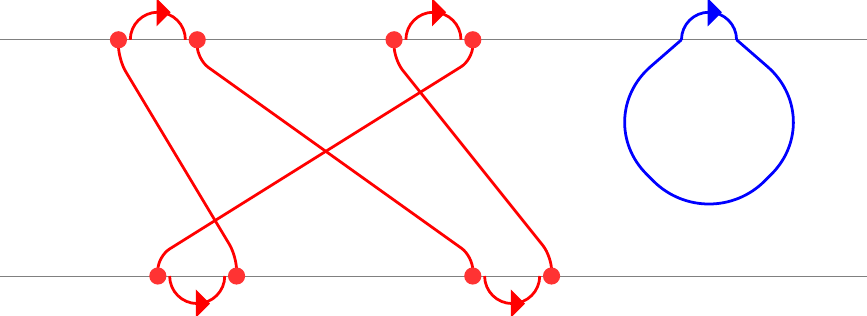}}}}
\end{tabular}\end{center}
\caption[caption]{\label{fig_perm_constructed} The permutations constructed in the proof of proposition \ref{pro_structure_2}, first case. $\tau_1$ and $\tau_2$ have invariant $(\{n\!-\!2\},s)$ and $(\{n\!-\!2\},-\!s)$ respectively by proposition \ref{pro_opposite_sign}}
\end{figure}

Let us begin. By lemma \ref{lem_top_adja_arc} there exists an edge $e$ of $\s$ such that its two top adjacent arcs are part of the cycles of length $2\ell$ and $n-2\ell-1$. 

Let $c$ be the coloring where $e$ is gray and $\tau$ the reduction of $(\s,c)$. By proposition \ref{pro_top_adj_arc}, $\tau$ has size $n'-1$ and invariant $(\lam'= \{n-2\},s)$ for some $s\neq 0$ ($s\neq 0$ by the classification theorem \ref{thm.Main_theorem_2}). Moreover we can assume that $\tau$ is not in the exceptional class by lemma \ref{lem_id_1}.

Now by proposition \ref{pro_structure_1} there exists $S$ and $\tau_1=S(\tau)$ , $\pi$ with invariant $(\lam''=\{2,n-5\},0)$ and $\Pi$ a consistent labelling of $\pi$ such that $\tau_1=\pi|(1,t_{1,n-5,1},b_{0,2,1})$. Let $e'$ be this edge of $\tau_1$.

By induction, there exists $S'$ boosted sequence of $S$ such that $(\s',c')=S'(\s,c)$ and $(\s',c')$ has reduction $\tau_1$.

Now, by lemma \ref{lem_where_arc}, the two top/bottom arcs of the cycle of length 2 of $\pi$ correspond to three top/bottom arcs of the unique cycle of length $n-2$ of $\tau_1$. Moreover those arcs are consecutive in the cyclic ordering of the cycle, indeed we have: 
\begin{center}\begin{tabular}{cc}
\includegraphics[scale=.75]{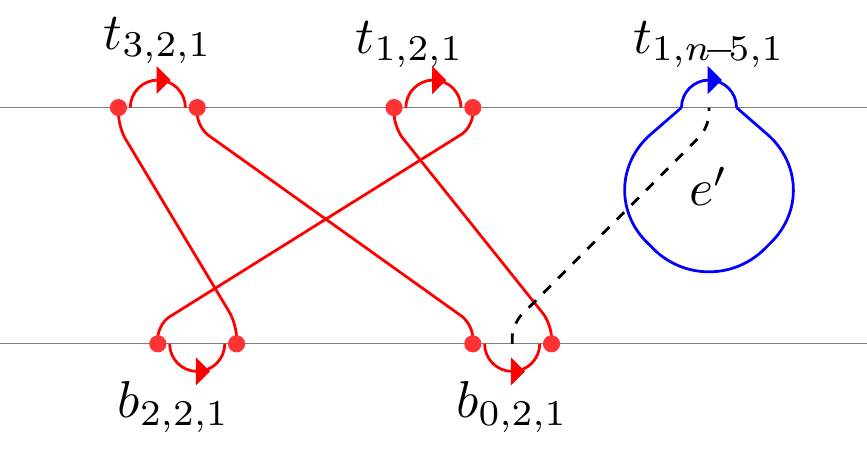} &
\includegraphics[scale=.75]{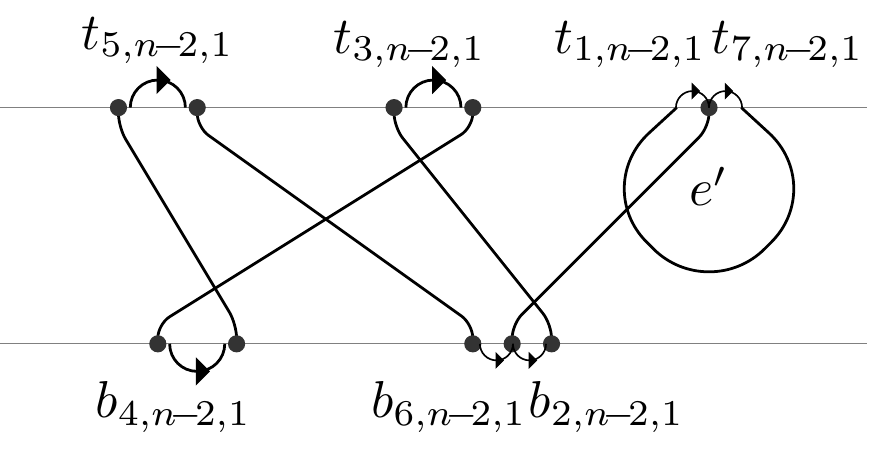}\\
$\pi,\Pi$ & $\tau_1,\Pi'$
\end{tabular}\end{center}
Where $\Pi'$ is a consistent labelling of $\tau_1$.

Due to the cycle invariant of $\s'$ and $\tau_1$ we have (by proposition \ref{pro_one_edge_two})
$\s'=\tau_1|(1,t_{2x+1,n-2,1},b_{2y,n-2,1})$ with $(\frac{q_{2(n-2)}(2x+1,2y)+1}{2},\frac{q_{2n-2}(2x+1,2y)+1}{2})=(2\ell,n-2\ell-1)$ or $(n-2\ell-1,2\ell)$ (let us say $(2\ell,n-2\ell-1)$). 
  
Let $(\s_1,c_{e})=S_e^{-2x}(\s',c')$ so that $\s_1=\tau_1|(1,t_{1,n-2,1},b_{2y' ,n-2,1})$ with $2y'=2y-2x \!\!\! \mod 2(n-2)= 1+q_{2(n-2)}(2x+1,2y)$ (by proposition \ref{pro_Se_action}). Then by lemma \ref{lem_where_arc2} the arcs labelled $b_{2,n-2,1},t_{3,n-2,1}$, $b_{4,n-2,1},t_{5,n-2,1},b_{6,n-2,1}$ of $(\tau_1,\Pi')$ correspond to two top and three bottom arcs of the cycle of length $2\ell$ in $\s_1$ since $2,3,4,5,6 \in ]1,2y'[$ (indeed $2\ell\geq n-2\ell-1$ so $2y'>6$ for $n$ large enough). More precisely, we have
\begin{center}\begin{tabular}{cc}\includegraphics[scale=.75]{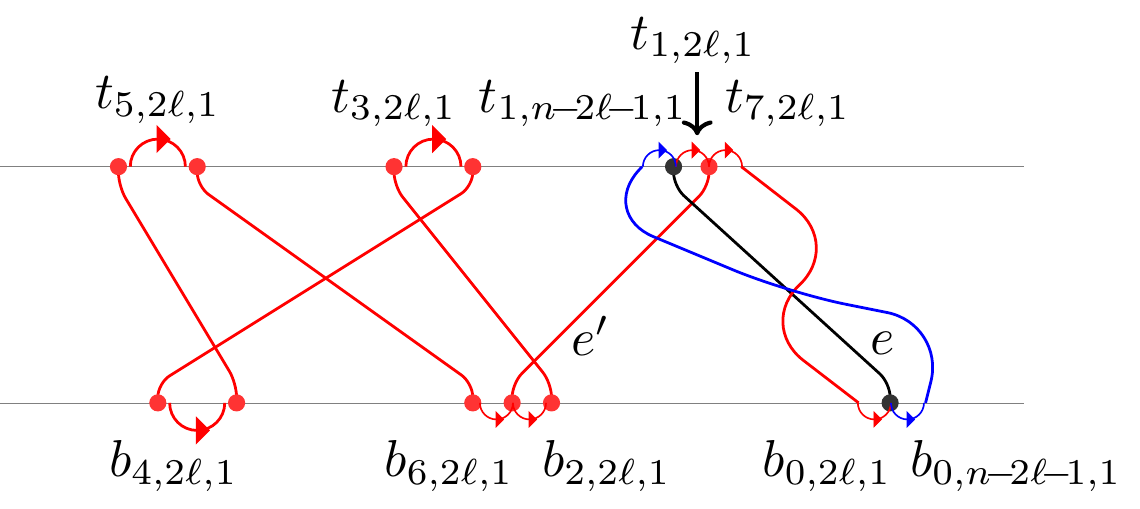} \\
$\s_1,\Pi''$
\end{tabular}\end{center}
where $\Pi''$ is a consistent labelling of $\s_1$, and the arcs $b_{2,n-2,1},t_{3,n-2,1}$, $b_{4,n-2,1},t_{5,n-2,1},b_{6,n-2,1}$ of $(\tau_1,\Pi')$ correspond to $b_{2,2\ell,1},t_{3,2\ell,1}$, $b_{4,2\ell,1},t_{5,2\ell,1},b_{6,2\ell,1}$ of $(\s_1,\Pi'')$. Moreover the left-adjacent top arc of $e'$ is $t_{1,2\ell,1}$ and the right-adjacent is $t_{7,\ell,1}$ 

Let us consider $c_{e'}$ the coloring of $\s_1$ where $e'$ and $\tau_{e'}$ the reduction, then by proposition \ref{pro_top_adj_arc} $\tau_{e'}$ has invariant $\{2,2\ell-3,n-2\ell-1\}$ and the gray edge $e'$ lies within a bottom arc of the cycle of length 2 and a top arc of the cycle of length $2\ell-3$. we have:
\begin{center}\begin{tabular}{cc}\includegraphics[scale=.75]{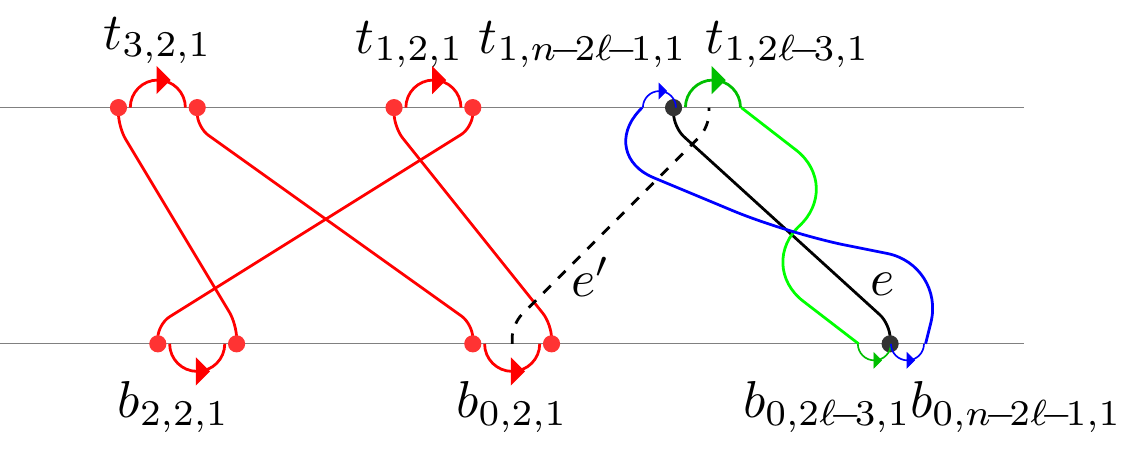} \\
$\tau_{e'},\Pi_{e'}$
\end{tabular}\end{center}
where $\Pi_{e'}$ is a consistent labelling of $\tau_{e'}$. Moreover the arc labelled $b_{0,2,1}$ and $b_{2,2,1}$ of $\tau_{e'},\Pi_{e'}$ correspond to the arcs with the same labels in $\pi,\Pi$  

$2\ell-3$ is odd and 2 is even, thus, we have (by proposition \ref{pro_Se_action}) $\s_2=S^{2(2\ell-3)}_{e'}(\s_1)=\tau_{e'}|(1,t_{1,2\ell-3,1},b_{2,2,1})$ since $\s_1=\tau_{e'}|(1,t_{1,2\ell-3,1},b_{0,2,1})$. 

Let us now define the coloring $c'_e$ of $\s_2$ where $e$ is grayed and let $\tau_2$ be the reduction. 
Then we have $\tau_2=\pi|(1,t_{1,n-5,1},b_{2,2,1})$ and $\tau'=\pi|(1,t_{1,n-5,1},b_{0,2,1})$ thus both have opposite sign by proposition \ref{pro_opposite_sign} $-s$ and $s$ respectively and same cycle invariant ($\lam'=\{n-2\}$).

\paragraph*{Second case.} Let $\ell$ be the largest cycle of $\s$ and $\ell_1$ be any other cycle. We will show that there exists $\s'$ connected to $\s$ and an edge $e$ of $\s'$ such that the two top adjacent arcs of $e$ are part of the cycles of length $\ell$ and of $\ell_1$.

Define $\tau_1$ the reduction of $\s'$ where the edge $e$ is grayed, then by proposition \ref{pro_top_adj_arc} the gray edge is within a cycle of length $k=\ell+\ell_1-1$ in $\pi$ and there are no other cycle of length $k$ in $\pi$ since it has cycle invariant $\lam'=\lam \setminus \{\ell,\ell'\}\bigcup \{k=\ell+\ell_1-1\}$ and thus $k$ is the unique largest cycle of $\pi$. Remains the sign invariant, we will show that by choosing $\ell_1$ carefully and applying lemma  \ref{lem_sign_control}  the sign invariant of $\pi$ is still $s$ and this conclude proposition \ref{pro_structure_2}.\\

By lemma \ref{lem_top_adja_arc} there exists an edge $e$ of $\s$ such that its two top adjacent arcs are part of the cycles of length $\ell$ and $\ell'$ for some $\ell'$.

If $\ell'=\ell_1$ we are done. Otherwise let $c$ be the coloring where $e$ is gray and $\tau$ the reduction of $(\s,c)$. By proposition \ref{pro_top_adj_arc}, $\tau$ has size $n-1$ and cycle invariant $\lam'=\lam \setminus \{\ell,\ell'\}\bigcup \{m=\ell+\ell'-1\}$ thus $m$ is the (unique) largest cycle of $\tau$. Moreover since $\ell\geq\ell'\geq 2$, $m\geq 3$ thus by proposition \ref{pro_T_end_perm} there exists $\tau'$ with the same invariants as $\tau$ and a consistent labelling $\Pi$ such that: 
\[(\tau',\Pi)= \raisebox{-35pt}{\includegraphics[scale=.5]{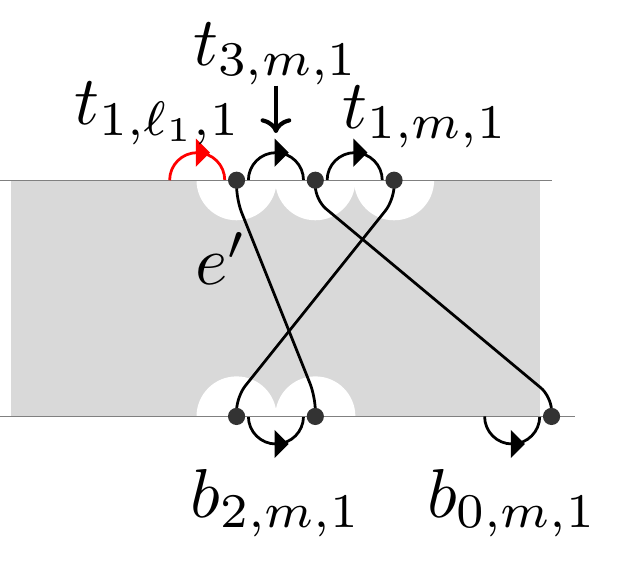}}\]
where $\ell_1$ is any cycle of both $\lam$ and $\lam'$.

Then by induction, since the classification theorem is proven for $n'-1, \lam'$ there exists a sequence $S$ such that $S(\tau)=\tau'$. Thus there exists a boosted sequence $S'$ such that $(\s',c')=S'(\s,c)$ where $\tau'$ is the reduction of  $(\s',c')$. Moreover since $(\s',c')$ has invariant $\lam$ and $\tau'$ has invariant  $\lam'=\lam \setminus \{\ell,\ell'\}\bigcup \{m=\ell+\ell'-1\}$ by proposition \ref{pro_one_edge_two} the gray edge must be inserted within the cycle of length $m$ in the following way: \[\s'=\tau'|(1,t_{2x+1,m,1},b_{2y,m,1}) \text{ with } (\frac{q_{2m}(2x+1,2y)+1)}{2},\frac{q_{2m}(2y,2x+1)+1}{2})=(\ell,\ell') \text{ or }(\ell',\ell).\]

Let us assume wlog that $(\frac{q_{2m}(2x+1,2y)+1)}{2},\frac{q_{2m}(2y,2x+1)+1}{2})=(\ell,\ell')$ then we choose $q$ such that $3 \in ]2x+1+q,2y+q[$ (the interval is defined modulo $2m$) then by lemma \ref{lem_where_arc} the arc of $S_e^q(\s')$ corresponding to $t_{3,m,1}$ is part of the cycle of length $\ell$ while the arc corresponding to $t_{1,\ell_1,1}$ is untouched and still part of $\ell_1$. Then the edge $e'$ of $\s'$ verifies that its two adjacents top arcs are parts of the cycles of length $\ell$ and $\ell_1$ respectively.

Let us now consider the sign. If $s\neq0$ there are no even cycle in $\lam$, let $\ell_1$ be any cycle,  let us say the second largest cycle then $s'=s$ by lemma \ref{lem_sign_control} since $\lam'$ does not contain even cycle either. If $s=0$ then $\lam$ has an even number of even cycle (cf theorem \ref{thm.Main_theorem_2} and since we are in the second case $\lam \neq \{2\ell, n-2\ell-1\}$ thus there are either at least an odd cycle or at least four even cycles. 

If $\ell$ is odd let $\ell_1$ be the largest even cycle then $\lam'$ still contains even cycle and thus $s'=s=0$, otherwise if $\ell$ is even and there are odd cycles then let $\ell_1$ be the largest odd cycle then $\lam'$ still contains even cycle and thus $s'=s=0$
Finally if there are no odd cycles let $\ell_1$ be the second largest even cycle then $\lam'$ still contains even cycle since in this case $\lam$ had at least 4 of them.
\end{proof}

Let us now carry out the proof of the classification theorem.

\begin{proof}[Proof of the classification theorem \ref{thm_main}]
 First case of the induction: Let $\s_1,\s_1'$ be two permutations of size $n$ with the same invariants (not in the exceptional class) $(\lam,s)$ and $|\lam|>1$. Let us apply proposition \ref{pro_structure_2} to $\s_1$ and $\s_1'$. Let $(\s,c)$ and $(\s,c')'$ connected to $\s_1$ and $\s_1'$ respectively be the obtained permutations and let $\tau$ and $\tau'$ be their reductions. Then we have the following : $\tau$ and $\tau'$ have cycle invariant $\lam'=\lam\setminus \{\ell,\ell'\} \cup\{k=\ell+\ell'-1\}$ where $\ell$ and $\ell'$ are the two cycles of $\lam$ and $\lam'$ has only one cycle of length $k$. Moreover $\tau$ and $\tau'$ have same sign invariant. 

Thus by the classification theorem at size $n-1$, there exists $S$ such that $S(\tau)=\tau'$. Let $S'$ be the boosted sequence we have $S'(\s,c)=(\s'',c'')$ and the reduction of $(\s'',c'')$ is $\tau'$. Let $\Pi$ be a consistent labelling of $\tau'$. By proposition \ref{pro_one_edge_two} since $\tau'$ has invariant $\lam\setminus \{\ell,\ell'\} \cup\{\ell+\ell'-1\}$ and both $(\s',c')$ and $(\s'',c'')$ have invariant $\lam$ we must have $\s'=\tau'|(1,t_{2x+1,k,1},b_{2y,k,1}$ and $\s''=\tau'|(1,t_{2x'+1,k,1},b_{2y',k,1})$ with 
\[\left(\frac{q_{2m}(2x+1,2y)+1)}{2},\frac{q_{2m}(2y,2x+1)+1}{2}\right)=(\ell,\ell') \text{ or }(\ell',\ell)\]
and 
\[\left(\frac{q_{2m}(2x'+1,2y')+1)}{2},\frac{q_{2m}(2y',2x'+1)+1}{2}\right)=(\ell,\ell') \text{ or }(\ell',\ell)\]
Thus we are in the conditions of propostion \ref{pro_connect_1} and there exists $i$ such that $S_e^i(\s'')=\s'$\\

 Second case of the induction: Let $\s_1,\s_1'$ be two permutations of even size $n+2$ with the same invariants (not in the exceptional class) $(\lam=\{n+1\},s)$. 
Let us apply proposition \ref{pro_structure_1} to $\s_1$ and $\s_1'$. Let $(\s,c)$ and $(\s,c')'$ connected to $\s_1$ and $\s_1'$ respectively be the obtained permutations and let $\tau$ and $\tau'$ be their reductions. Then we have the following : $\tau$ and $\tau'$ have cycle invariant $\lam'=\{2,n-2\}$ and same sign invariant 0. Moreover since $n+2$ is odd $n-2$ is even.

Thus by the classification theorem at size $n-1$, there exists $S$ such that $S(\tau)=\tau'$. Let $S'$ be the boosted sequence we have $S'(\s,c)=(\s'',c'')$ and the reduction of $(\s'',c'')$ is $\tau'$. Let $\Pi$ be a consistent labelling of $\tau'$. By proposition \ref{pro_one_edge} since $\tau'$ has invariant $\{2,n-2\}$ and both $(\s',c')$ and $(\s'',c'')$ have invariant $\lam$ we must have $\s'=\tau'|(1,t_{2x+1,2,1},b_{2y,n-2,1}$ and $\s''=\tau'|(1,t_{2x'+1,2,1},b_{2y',n-2,1})$ (or $\s'=\tau'|(1,t_{2x+1,n-2,1},b_{2y,2,1}$ or $\s''=\tau'|(1,t_{2x'+1,n-2,1},b_{2y',2,1})$)

Thus we are in the conditions of propostion \ref{pro_connect_2} and there exists $i$ such that $S_e^i(\s'')=\s'$.
\end{proof}

\section*{Second part}
\label{sec:intro}
\addcontentsline{toc}{section}{\nameref{sec:intro}}

In this part, we come back to the standard Rauzy dynamics that we studied in \cite{DS17} and \cite{D18}. In term of notational change, the top principal cycle is renamed the rank and is an invariant of the Rauzy dynamics. Thus if we had a permutation $\s$ with invariant $(\lam,s)$ in the extended Rauzy dynamics we obtain, in the standard Rauzy dynamics, a permutation with invariant $(\lam\setminus \{r\},r,s)$  where $r$ is the length of the top principal cycle.

\section{Lower bound on the diameter}

In this section we prove: 

\begin{theorem}\label{thm_lower_bound}
Let $C$ be a class at size $n$ then the diameter of $C$ for the alternation distance is at least $n/16$.
\end{theorem}

In order to do so we define the following new dynamics:

\begin{definition}[Pivotless Rauzy dynamics]
Let $(L^{i,j})_{1\leq i\leq n, 1\leq j<n-i}$ and $(R^{i,j})_{1\leq i\leq n, 1\leq j<i}$ be the two following sets of operators on $\kS_n$:
\begin{align*}
L^{i,j}\;\left(
\;
\raisebox{-25pt}{\includegraphics[scale=.499]{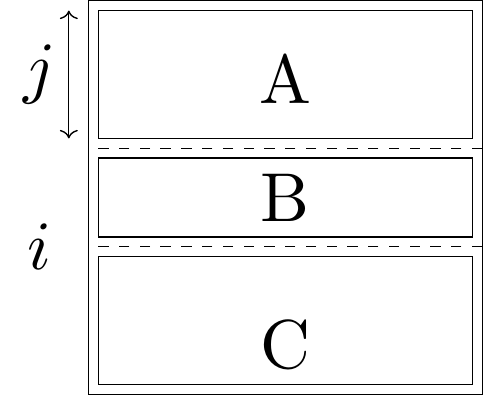}}
\;
\right)
&=
\raisebox{-25pt}{\includegraphics[scale=.499]{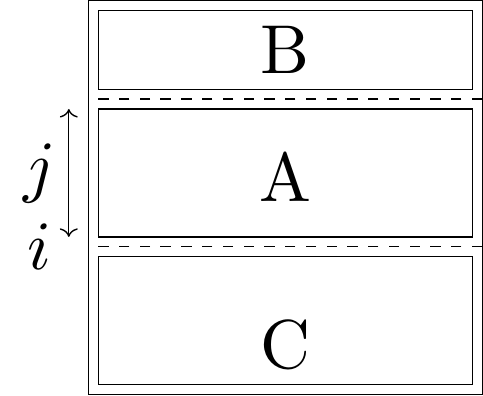}}
&
R^{i,j}\;\left(
\;
\raisebox{-25pt}{\includegraphics[scale=.499]{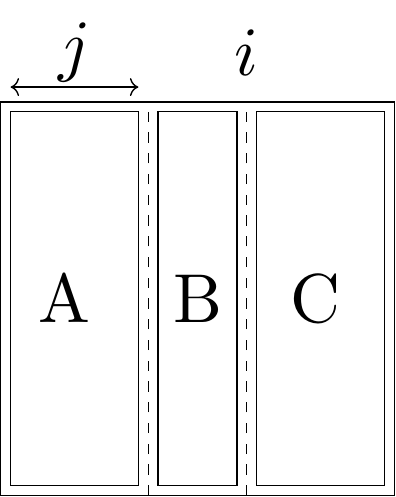}}
\;
\right)
&=
\raisebox{-25pt}{\includegraphics[scale=.499]{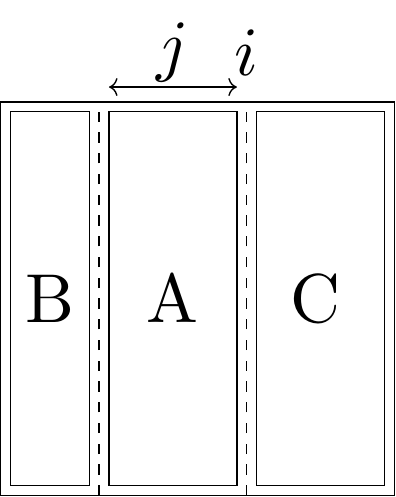}}
\rule{0pt}{38pt}
\end{align*}
We call the dynamics on $\kS_n$ with operators $(L^{i,j})_{1\leq i\leq n, 1\leq j<n-i}$ and $(R^{i,j})_{1\leq i\leq n, 1\leq j<i}$ the \emph{pivotless Rauzy dynamics}.
\end{definition}

Note that if $\s(1)=i$ then for all $j$, $L^{i,j}(\s)=L^j(\s)$ and if $\s(i)=n$ then for all $j$, $R^{i,j}(\s)=R^j(\s)$, as illustrated below: 
 \begin{align*}
L^{i,j}\!\left(
\raisebox{-25pt}{\includegraphics[scale=.475]{figure/fig_matrix_dyna_3_Lij.pdf}}
\,
\right)
&\!=\!L^{j}\!\left(
\raisebox{-25pt}{\includegraphics[scale=.475]{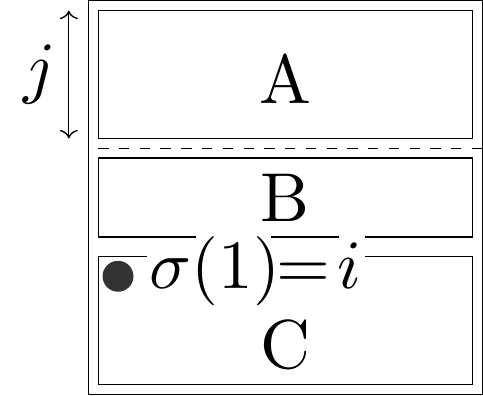}}
\,
\right)
&
R^{i,j}\!\left(
\raisebox{-25pt}{\includegraphics[scale=.475]{figure/fig_matrix_dyna_1_Rij.pdf}}
\right)
&\!=\!R^{j}\!\left(
\raisebox{-25pt}{\includegraphics[scale=.475]{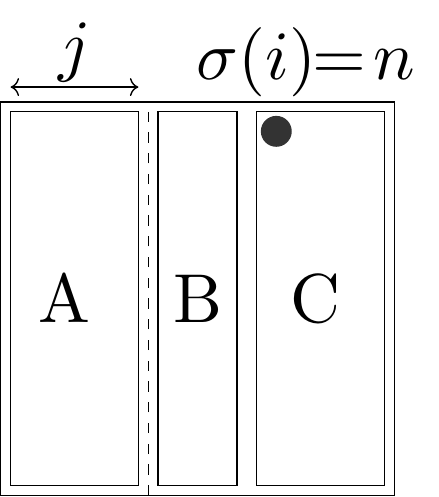}}
\! \! \!
\right)
\end{align*}
Thus the \emph{pivotless Rauzy dynamics} mimics the action of the Rauzy dynamics (with operators $(L^j)_j,(R^j)_j$) for any given pivot. 

Given a permutation $\s$ we can partition its set of edges in two colors: black and red. Given a partition $c$, the permutation $\tau$ corresponding to the restriction of $\s$ to the red edges is called the \emph{pivotless reduced permutation} for $\s,c$. 

For a pair $(\s,c)$, the Rauzy dynamics (with operators $(L^j)_j,(R^j)_j$) reduces to the pivotless Rauzy dynamic on $\tau$ (its pivotless reduced permutation) as follows: for every operators $R^j$ or $L^j$ we 
define
\[P(R^j)=\begin{cases}
Id & \text{If } \{(1,\s(1)),\ldots,(j,\s(j))\} \text{ has no red edge.}\\
R^{i',j'} & \text{If } \{(1,\s(1)),\ldots,(\s^{-1}(n),n)\} \text{ has $i'$ red edges}\\
& \text{and }  \{(1,\s(1)),\ldots,(j,\s(j))\} \text{ has $j'$ red edges.}
\end{cases}\]
and 
\[P(L^j)=\begin{cases}
Id & \text{If } \{(\s^{-1}(n-j+1),n-j+1)),\ldots,(\s^{-1}(n),n)\} \text{ has no red edge.}\\
R^{i',j'} & \text{If } \{(\s^{-1}(1),1),(\s^{-1}(2),2),\ldots,(1,\s(1))\} \text{ has $i'$ red edges}\\
& \text{and }   \{(\s^{-1}(n-j+1),n-j+1)),\ldots,(\s^{-1}(n),n)\}  \text{ has $j'$ red edges.}
\end{cases}\]
Likewise for a sequence $S=H_k^{j_k}\ldots H_1^{j_1}$ $(H_i\in \{L,R\}, H_i \neq H_{i+1})$, acting on $(\s,c)$, the \emph{pivotless sequence} of S, acting on $\tau$, is $P(S)=P(H_k^{j_k})\ldots P(H_1^{j_1})$.

In essence, The pivotless dynamic is the appropriate notion such that the following diagram makes sense: for a permutation $(\s,c)$ and a sequence $S$ calling $red$ the operator that extracts the pivotless reduced permutation $\tau$ from $(\s,c)$, we have  
\begin{center}
$\begin{tikzcd}
(\s,c)\arrow{d}{red} \arrow{r}{S} &  (\s',c')
 \arrow{d}{red}\\
  \tau \arrow{r}{P(S)} &\tau'\\
\end{tikzcd}$
\end{center}

\begin{figure}[t]
\begin{center}
$\begin{tikzcd}
 \raisebox{-15pt}{\includegraphics[scale=.5]{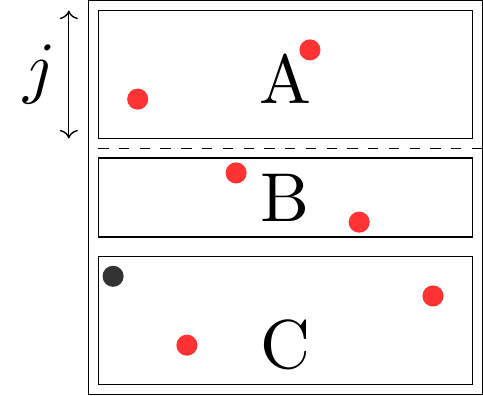}}\arrow{d}{red} \arrow{r}{L^j} &  
\raisebox{-15pt}{\includegraphics[scale=.5]{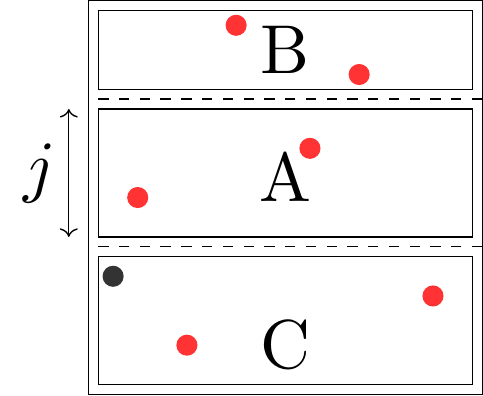}}\arrow{d}{red}\\
  \raisebox{-32pt}{\includegraphics[scale=.5]{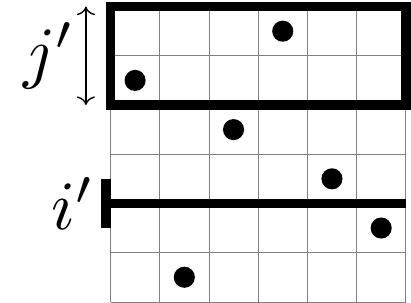}}\arrow{r}{P(L^j)=L^{i',j'}}[swap]{L^{2,2}} &  \raisebox{-32pt}{\includegraphics[scale=.5]{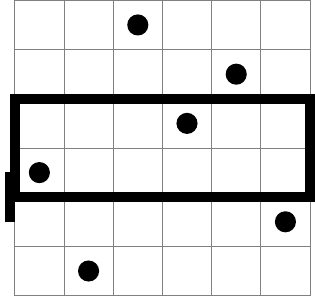}}
\\
\end{tikzcd}$
\end{center}
~\\[-15mm]
\caption{\label{fig_P_reduction} An example of the action of the pivotless dynamics. In the top row we have $(\s,c)$ and $(\s',c')=L^j(\s,c)$. In the bottom row are represented the respective pivotless reduced permutation $\tau$ and $\tau'=L^{i',j'}(\tau)$. Note that the pivotless sequence $P(L^j)=L^{2,2}$ since there are 2 red edges in the block $C=\{(\s^{-1}(1),1),(\s^{-1}(2),2),\ldots,(1,\s(1))\}$ and also two red edges in the block $A= \{(\s^{-1}(n-j+1),n-j+1)),\ldots,(\s^{-1}(n),n)\}$}
\end{figure}

See also figure \ref{fig_P_reduction} for an example.

Thus the pivotless dynamics can be considered a 'dual' of the boosted dynamics. For the boosted dynamics we had a Rauzy sequence $S$ on the reduced permutation $\tau$ and the boosted dynamics gave us a Rauzy sequence $B(S)$ on $\s$. For the pivotless dynamics, this is the opposite, we have a Rauzy sequence $S$ on the permutation $\s$ and the pivotless dynamics gives us a pivotless Rauzy sequence $P(S)$ on $\tau$. $P(S)$ cannot be a Rauzy sequence in this case, since, in moving from $\s$ to $\tau$, we discarded all the black edges of $\s$ some being the acting pivots in the sequence $S$.

Working with the pivotless dynamics gives lower bound in the alternating distance: if we show that any sequence $S'$ from $\tau$ to $\tau'$ has size at least $k$, then the alternating distance between $(\s,c)$ and $(\s',c')$ is also at least $k$ since any sequence $S$ such that  $(\s',c')= S(\s,c)$ must gives a pivotless sequence $P(S)$ such that $\tau'=P(S)(\tau)$ and $|P(S)|\leq|S|$ by definition of $P$ (where $|S|$ is in the alternation sense). 

We state this fact in the following lemma 
\begin{lemma}\label{lem_pivotless_alternation_distance}
Let $(\s,c)$ and $(\s',c')$ be two permutations with a black/red edge coloring, let $\tau$ and $\tau'$ be their respective pivotless reduced permutations and let $S$ be a sequence such that $S(\s,c)=(\s',c')$.

 If the standard distance (in the pivotless dynamics) between $\tau$ and $\tau'$ is at least $k$ then the alternation distance between $\s$ and $\s'$ is also at least $k$. 

\[ d^{pivotless}_{G}(\tau,\tau')\geq k \implies d^{A}(\s,\s')\geq k\]
\end{lemma}

To complete our discussion on pivotless dynamics, let us prove a lower bound in the stardard distance for the pivotless dynamics:

\begin{lemma}
Let $\tau=id_n$ and $\tau'=\omega \, id_n$ (i.e $\tau'=(n,n-1,\ldots,1)$) then  $d^{pivotless}_{G}(\tau,\tau')\geq n/2.$ 
\end{lemma}

\begin{pf}
We show that every operator can break at most 2 descents (by descent we mean $\s(i)=\s(i+1)+1$) thus since $\tau$ has $n-1$ descents and $\tau'$ has none, we must have  $d^{pivotless}_{G}(\tau,\tau')\geq (n-1)/2$.

Clearly since the blocks $A,B,C$ are just translated if a descent is broken it has to happen at the boundary of A and B and B and C, thus at most two descents are affected at each application of an operator.\qed
\end{pf}

A more involved argument could probably prove $d^{pivotless}_{G}(\tau,\tau')= n-1.$ 

In the following, we will show that the diameter of any class $C_n$ in the alternation distance is $\Omega(n)$ by constructing two permutations $\s$, and $\s'$ and a black/red edge coloring $c$ of $\s$ such that:\begin{itemize}
\item There exists a sequence $S$ such that  $S(\s)=\s'$.
\item For any sequence S such that  $S(\s)=\s'$, let $c_S$ be the coloring such that $(\s',c_S)=S(\s,c)$. Then $\tau$ and $\tau'$ the respective pivotless reduced permutations of $(\s,c)$ and $(\s',c_S)$ have $\tau=id_k$ and $\tau'=\omega \, id_k$ for some $k =\Omega(n)$.
\end{itemize} 

Then by lemma \ref{lem_pivotless_alternation_distance}, we will have proven theorem \ref{thm_lower_bound}.

Let us begin by constructing a permutation with a long increasing subsequence (thus coloring in red this increasing subsequence we have our $(\s,c)$ and our $\tau$).

\begin{proposition}\label{pro_increase}
For any class $C$ with invariant $(\lam,r,s)$, there exists a permutation $\s\in C$ with an increasing subsequence of size at least $n/2-d$ for some small $d$.
\end{proposition}

\begin{pf}
We will follow the construction of theorems 61 and 66 from \cite{D18} with exactly one minor change. 
For compactness sake we will only describe how to modify the construction of theorem 61 (the modification occurs at the same place for theorem 66).

Let us recall the construction. 
We start with a base permutation with invariant $(\lam'\subseteq \lam,r)$ then we add cycles to $\lam'$ with propositions \ref{pro_ins_cross_perm} and finally we adjust the sign $s(\s)$ by moving at most three edges.

For this new construct, we also start with base permutation $s^0$ with invariant $(\lam'\subseteq \lam,r)$. Clearly the base permutations all have an increasing subsequence of size $k/2$ (where $k$ is the size of the permutation). We reproduced figure 28 of \cite{D18} which illustrates the base permutations in figure \ref{fig_perm_x_family}.

The next step is to add cycles by inserting, one by one, some $C_i$ into the current permutation $\s^j$ in the order described in the theorem 61 of \cite{D18}. The only change in the construction happens here: instead of inserting the $C_i$ on the last edge of current permutation $\s^j$ we insert it in the following way: 

Let $\s^0$ be the base permutation and suppose that we insert (in this order) $C_{i_1},\ldots,C_{i_\ell}$ then \[\forall j\geq 1 \ \s^j=\s^{j-1}_{j+1}(C_{i_j}).\]
Refer to figure \ref{fig_const_LIS} for an example of this insertion scheme.

\begin{figure}[tb!]
\begin{center}
\begin{tabular}{lll}
\makebox[0pt][l]{\raisebox{0mm}{\includegraphics[scale=.46]{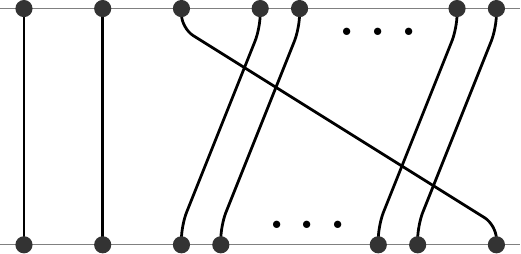}}}
\hspace{28pt}\raisebox{31pt}{$\rotatebox{90}{$\left. \rule{0pt}{21.5pt}\right\}\rotatebox{-90}{$\!\!\!\! 2k$} $}$}
&
\makebox[0pt][l]{\raisebox{0mm}{\includegraphics[scale=.46]{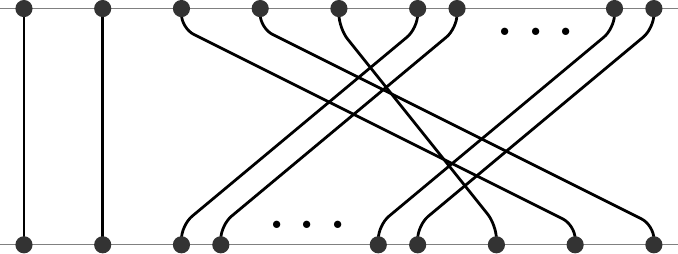}}}
\hspace{49pt}\raisebox{31pt}{$\rotatebox{90}{$\left. \rule{0pt}{21.5pt}\right\}\rotatebox{-90}{$\! \! \! \! 2k$} $}$}&
\makebox[0pt][l]{\raisebox{0mm}{\includegraphics[scale=.46]{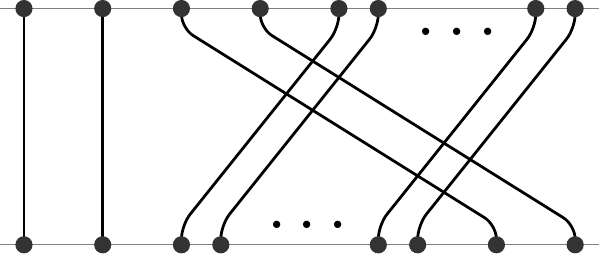}}}
\hspace{38.5pt}\raisebox{31pt}{$\rotatebox{90}{$\left. \rule{0pt}{21.5pt}\right\}\rotatebox{-90}{$\! \! \! \! 2k$} $}$}\\
$\lambda=\{2k\!+\!1\},r=1$ & $\lambda=\{2k\!+\!1\},r=3$ & $\lambda=\emptyset,r=2k\!+\!3$\\
&&\\

\makebox[0pt][l]{\raisebox{0mm}{\includegraphics[scale=.46]{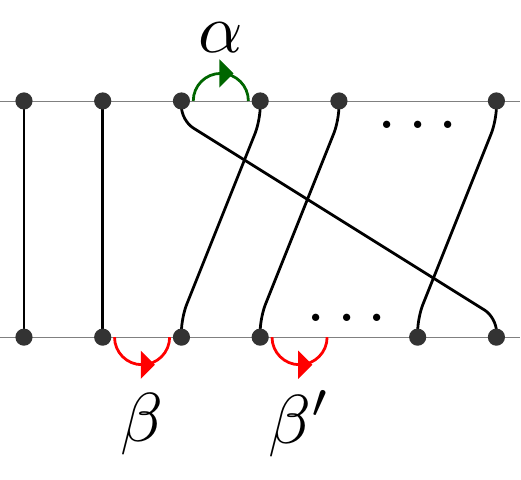}}}
\hspace{29pt}\raisebox{49pt}{$\rotatebox{90}{$\left. \rule{0pt}{21.5pt}\right\}\rotatebox{-90}{$\!\!\!\!\!\!\!\!\! 4k\!+\!3$} $}$}
&
\makebox[0pt][l]{\raisebox{0mm}{\includegraphics[scale=.46]{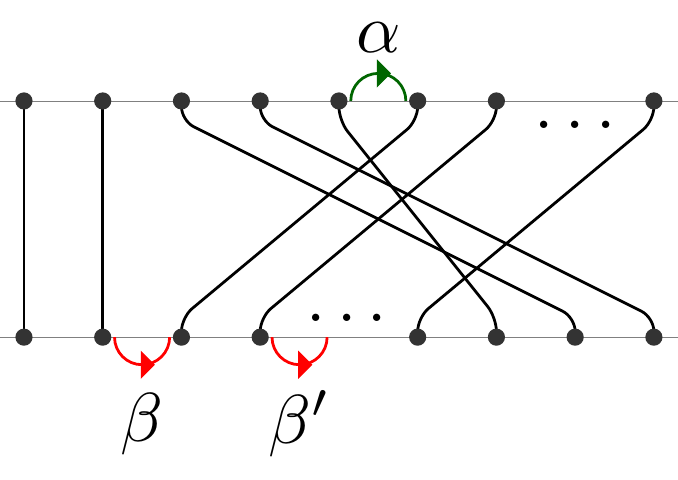}}}
\hspace{50pt}\raisebox{49pt}{$\rotatebox{90}{$\left. \rule{0pt}{21.5pt}\right\}\rotatebox{-90}{$\!\!\!\!\!\! \! \! \! 4k\!+\!3$} $}$}&
\makebox[0pt][l]{\raisebox{0mm}{\includegraphics[scale=.46]{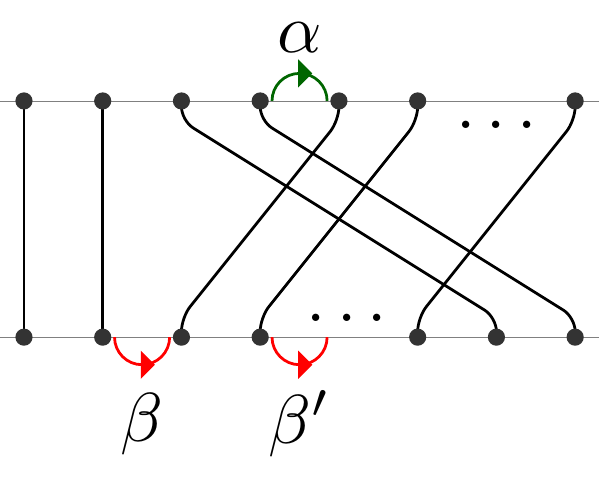}}}
\hspace{39.5pt}\raisebox{49pt}{$\rotatebox{90}{$\left. \rule{0pt}{21.5pt}\right\}\rotatebox{-90}{$\!\!\!\!\!\! \! \! \! 4k\!+\!3$} $}$}\\
$\lambda=\{{\color{red}2k\!+\!2},{\color{green!50!black}2k\!+\!2}\},r=1$ & $\lambda=\{{\color{red}2k\!+\!2},{\color{green!50!black}2k\!+\!2}\},r=3$ & $\lambda=\{{\color{red}2k\!+\!2}\},r={\color{green!50!black}2k\!+\!4}$
\end{tabular}
\end{center}
\caption{\label{fig_perm_x_family} Two families of base permutations with their respective cycle invariant.}
\end{figure}

\begin{figure}[tb!]
\begin{center}$
\begin{array}{lclcl}
\raisebox{0mm}{\includegraphics[scale=.4]{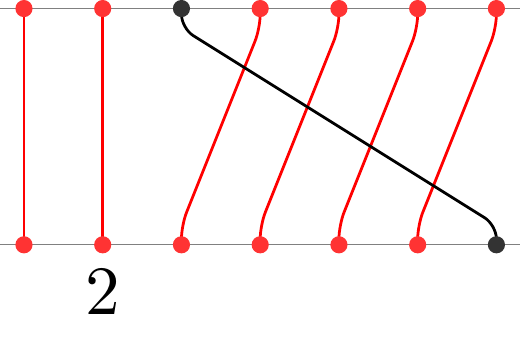}}
&&
\raisebox{0mm}{\includegraphics[scale=.4]{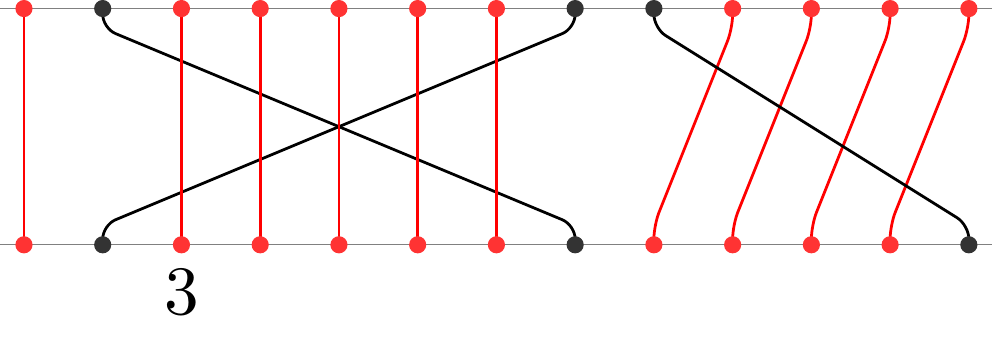}}
&&
\raisebox{0mm}{\includegraphics[scale=.4]{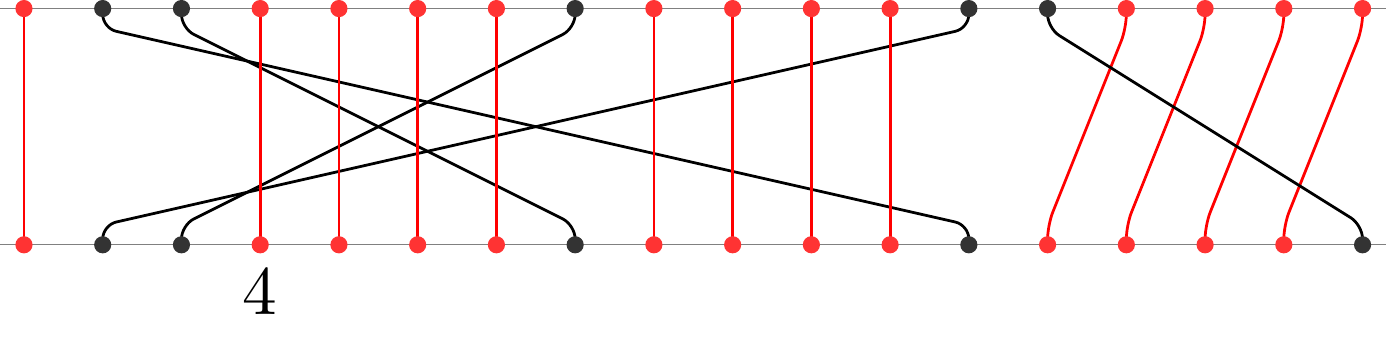}}
\\
\s^{0}: & & \s^1=\s^0_2(C_5):&&\s=\s^2=\s^1_3(C_4):\\
(\lam',r)=(\{5\},1)&&(\lam^1,r)=(\{3,3,5\},1)&&(\lam,r)=(\{3,3,5,5\},1)\\
Size: 7& & Size: 13&&Size: 18\\
LIS: 6 &&LIS: 10&&LIS: 13\\
\end{array}$
\end{center}
\caption{\label{fig_const_LIS} An example of the procedure to construct a permutations $\s$ with invariant $(\{3,3,5,5\},1)$ and a LIS (in red) of size at least $n/2$.}
\end{figure}

In this insertion procedure, each $C_i$ added contributes to the longest increasing subsequence (LIS) by $i-1$ (one indice of the LIS in $\s^{J-1}$ is replaced by $i$ indices in $\s^{j}$) and increase the size of the current permutation by $i+1$:
\begin{center}$
\begin{array}{lcl}
\makebox[0pt][l]{\raisebox{0mm}{\includegraphics[scale=.5]{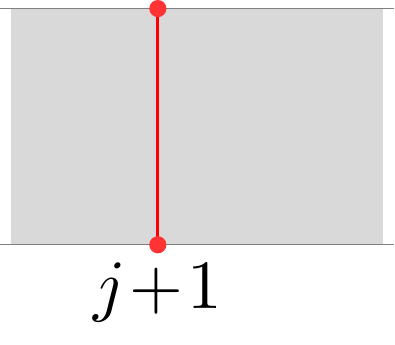}}}
&$\qquad\qquad\qquad$&
\makebox[9pt][l]{\raisebox{4.2mm}{\includegraphics[scale=.5]{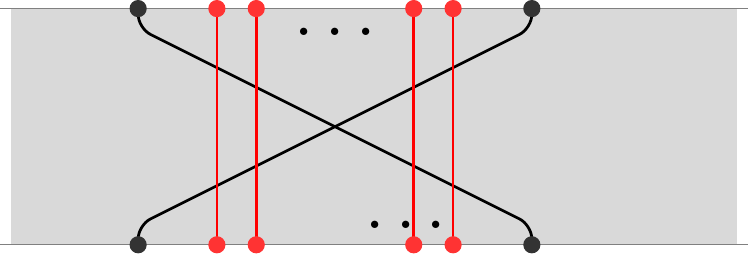}}}
\hspace{18pt}\raisebox{34pt+4.2mm}{$\rotatebox{90}{$\left. \rule{0pt}{22.5pt}\right\}\rotatebox{-90}{$\!\!i$} $}$}\\
\s^{j-1}: & & \s^j:\\
Size: k&&Size: k+i+1\\
LIS: k' &&LIS: k'+i-1\\
\end{array}$
\end{center}
Thus, by induction, the LIS of $\s^{\ell}$ must be at least $|\s^{\ell}|/2$ as long as the $(C_{i_j})_j$ inserted have $i_j>2, \forall j$. This is the case since by construction (refer to theorem \ref{thm_||Xshaped}) at most one $C_{i_j}$ has $i_j=2$.

Finally, we adjust the sign of $\s^\ell$, which is accomplished by moving at most three edges of $\s^\ell$. Thus up to a small constant $c<8$, we have constructed a permutation $\s$ with invariant $(\lam,r,s)$ and with a LIS of size at least $n/2-d$.

\qed
\end{pf}

We have constructed the permutation $(\s,c)$ and we now need a permutation $\s'$ such that for any sequence $S,$ the coloring $c'$ in $S(\s,c)=(\s',c')$ is a decreasing subsequence. This is achieved with the proposition below: 

\begin{proposition}\label{pro_decrease}
For any class $C$ with invariant $(\lam,r,s)$, there exists a permutation $\s\in C$ such that the set of edges can be partitionned into 5 subsets: a subset has size at most 6 and the other subsets are all decreasing subsequences. 
\end{proposition}

Indeed, by the pigeonhole principle, the red increasing subsequence of $(\s,c)$ of size $n/2-d$ will distribute into the five subsets of $\s'$ and thus $(\s',c')$ will have a red decreasing subsequence of size at least $(n/2-d-6)/4$ which proves theorem \ref{thm_lower_bound}.

To construct the permutations of this proposition \ref{pro_drecrease}, we will use the $T,q_1,q_2$ operators of article \cite{DS17}. For that purpose, let us recall some facts: 

\begin{lemma}[Representant for the $T,q_1,q_2$ operators]
For any class $C$ with invariant $(\lam,r,s)$, there exists a sequence of operators \[S=T^{i_k}q_{j_k}\ldots T^{i_1}q_{j_1}T^{i_0} \text{ with } \forall 1\leq \ell \leq k-1, \ i_\ell\geq 1,\ i_0,i_k\geq 0 \text{ and  } \forall \ell, \ j_\ell \in \{1,2\}\]
 such that $S(id_m)\in C$ for some $5\leq m\leq 7.$
\end{lemma}

A second remark is that the $T$ operator needn't be applied to the first edge $(1,\s(1))$: it can be applied to any edge without changing the image of the class.

\begin{definition}
We define the $T_i$ operator $T_i:$ $\kS_n \to \kS_{n+2}$ as:
\begin{center}
\begin{tabular}{ccc}
$\!\!\! \!\!\!\!\!\!\! \s:$ && $T_i(\s):$\\
\raisebox{0mm}{\includegraphics[scale=1]{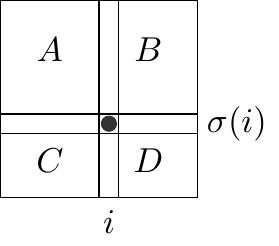}}
& & \raisebox{1mm}{\includegraphics[scale=1]{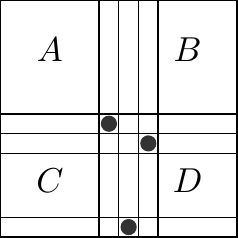}}
\end{tabular}
\end{center}
By definition $T_1=T$.
\end{definition}  

This is the same $T_i$ as the ones defined in section \ref{sdection_tech}.

\begin{lemma}
Let $\s$ be an irreducible permutation then for every $i$, $T(\s) \sim T_i(\s)$.
\end{lemma}

As we mentionned below proposition \ref{pro_T_action}, this lemma is a consequence of the corollary 3.14 of \cite{DS17}

We construct our permutation $(\s',c')$ (where $c'$ is the coloring into the 5 subsets described in proposition \ref{pr_decrease}) inductively by applying operators $T,q_1,q_2$, thus we first define colored versions of them.

\begin{definition}
Let $Col_n$ be the set of 5 coloring of $[1,\ldots,n]$ : $c\in Col_n$ if and only if $c: [1,\ldots,n] \to \{red,blue,green,brown,black\}$.
Let $(\s,c)\in (\kS_n,Col_n)$ and let $(i,\s(i))$ be the rightmost $red$ edge of $(\s,c)$, then we define

\begin{itemize} 
\item$T_{r}: \ (\kS_n,Col_n) \to (\kS_{n+2},Col_{n+2})$ to be the following operator:  
\[ T_r(\s)=T_i(\s) \quad 
T_r(c)(j)=\begin{cases}
c(j) & \text{if } j\leq i\\
blue & \text{if } j= i+1\\
red & \text{if } j= i+2\\
c(j-2) & \text{if } j> i+2\\
\end{cases}\]
Thus $T_r$ applies $T$ on the rightmost $red$ edge of $\s$ and color the two newly added edges is color $blue$ and $red$ (the $blue$ edge is the new edge $(i+1,1)$ and the $red$ edge is the new edge $(i+2,\s(i)+1)$): 
\begin{center}
\begin{tabular}{ccc}
$\!\!\! \!\!\!\!\!\!\! \s:$ && $T_r(\s):$\\
\raisebox{0mm}{\includegraphics[scale=1]{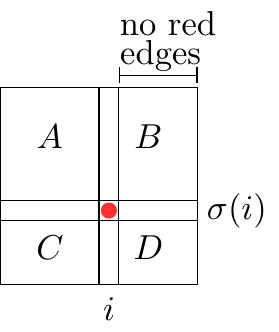}}
& & \raisebox{1mm}{\includegraphics[scale=1]{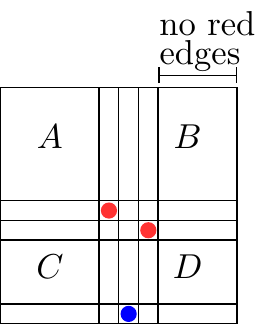}}
\end{tabular}
\end{center}

\item $q^{col}_k: \ (\kS_n,Col_n) \to (\kS_{n+1},Col_{n+1})$ where $k\in \{1,2\}$ and $col\in \{red,blue,green,brown,black\}$ to be the following operators:  
\[ q^{col}_k(\s)=q_k(\s) \quad 
q^{col}_k(c)(j)=\begin{cases}
c(j) & \text{if } j < q_k(\s)^{-1}(1)\\
col & \text{if } j= q_k(\s)^{-1}(1)\\
c(j-1) & \text{if } j > q_k(\s)^{-1}(1)\\
\end{cases}\]
Thus $q^{col}_k$ applies $q_k$ to $\s$ and color the newly added edge in color $col$ 
\end{itemize}
\end{definition}

The permutation $\s',c'$ of proposition \ref{pro_decrease} is then 

\begin{lemma}
Let C be the class with invariant $(\lam,r,s)$ and let 
\[S=T^{i_k}q_{j_k}\ldots T^{i_1}q_{j_1}T^{i_0} \text{ with } \forall 1\leq \ell \leq k-1, \ i_\ell\geq 1,\ i_0,i_k\geq 0 \text{ and  } \forall \ell, \ j_\ell \in \{1,2\}\]
 be a sequence such that $S(\s_0)\in C$ with $\s_0=id_m$ for some $5\leq m\leq 7.$ Let $c_0$ be the coloring of $\s_0$ defined by 
\[c_0(i)=\begin{cases} rd & \text{if } i=1\\ bk & \text{otherwise.}\end{cases}\]
Let $(\ell_p)_{1\leq p\leq r}$ be the sequence of indices of the $q_1$ in $S$ i.e \[\forall 1\leq p\leq r, j_{\ell_p}=1 \quad \forall p' \notin \{(\ell_p)_p\}\ j_{p'}=2\] 
Finally define $S'$ a colored version of sequence $S$ to be $S'=S'_2S'_1$ with 
\[S'_1=T_r^{i_{\ell_1}\!-\!1}q^{red}_2 \ldots T_r^{i_1}q^{red}_2T_r^{i_0}\]
and 
\[S'_2=T_r^{i_k}\ldots T_r^{i_{\ell_2}\!+\!1}q^{brown}_2T_r^{i_{\ell_2}} q_{1}^{green}T_r^{i_{\ell_2}\!-\!1}q^{brown}_2   \ldots T_r^{i_{\ell_1}\!+\!1}q^{brown}_2T_r^{i_{\ell_1}} q_{1}^{green}\]
Thus in $S'$, the colored version of $q_1$ is always $q^{green}_1$. For $q_2$, in $S'_1$ the colored version is $q^{red}_2$ and in $S'_2$ it is exclusively $q^{brown}_2$. 

Then the permutation $(\s',c')=S'(\s_0,c_0)$ verifies the condition of proposition \ref{pro_decrease}: the color $bk$ correspond to the subset of size $<6$ and the other colors are all decreasing subsequences.
\end{lemma}

\begin{pf}
we will represent the permutation $S'(id_s)$ for $s=5,6,7$. In our figures below $s=6$, but the case $s=5$ and $s=7$ are just as simple: For $s=5$ it suffices to remove the lowest black point and for $s=7$ we replace it by an identity of size 2.

The figures are a straigthforward application of the colored operators, the only justification required is the placement of the red, green and brown dots when applying $q_2^{red}$,$q_1^{green}$ and $q^{brown}_2$. Recall that $q_1$ and $q_2$ change the rank to 1 (respectively 2). In the case of $q_2^{red}$ the permutations have the form: \[\tau=*,n-3,n-2,n-1,n \text{ thus } q_2^{red}(\tau)=(*+1),{\color{red}1},n-2,n-1,n,n+1.\]
\begin{center}
\begin{tabular}{cc}
\raisebox{4pt}{\includegraphics[scale=.5]{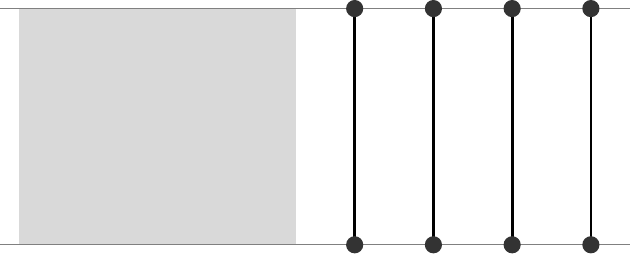}} &\includegraphics[scale=.5]{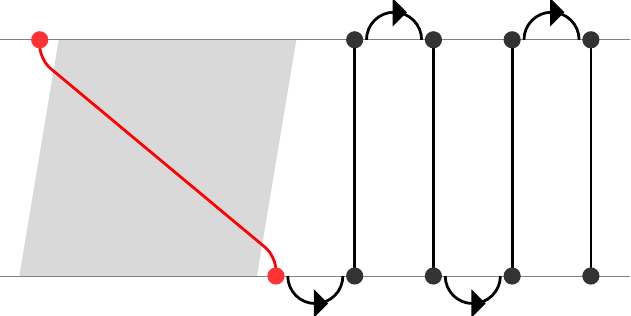}
\end{tabular} 
\end{center}
\begin{center}
\begin{tabular}[t]{cccc}
$T^{i_0}(id_s)$ & $q_2^{red}T^{i_0}(id_s)$&$q^{red}_2T^{i_1}q_2^{red}T^{i_0}(id_s)$&$S_1'(id_s)$\\
\adjustimage{valign=T}{figure/big_figure1.pdf}  &\adjustimage{valign=T}{figure/big_figure2.pdf}  &\adjustimage{valign=T}{figure/big_figure3.pdf}  &\adjustimage{valign=T}{figure/big_figure4.pdf} 
\end{tabular}
\end{center}

For the case $q_1^{green}$ the permutations have the form: \[\tau=*,n-1,n \text{ thus } q_1^{green}(\tau)=(*+1),{\color{green}1},n,n+1.\]
\begin{center}
\begin{tabular}{cc}
\raisebox{4pt}{\includegraphics[scale=.5]{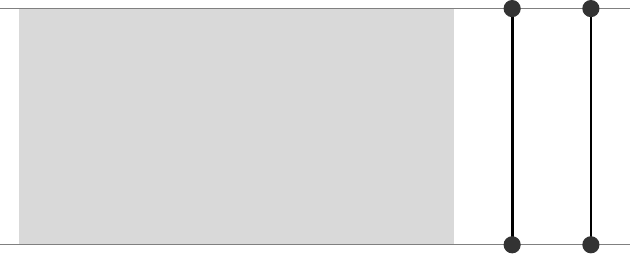}} &\includegraphics[scale=.5]{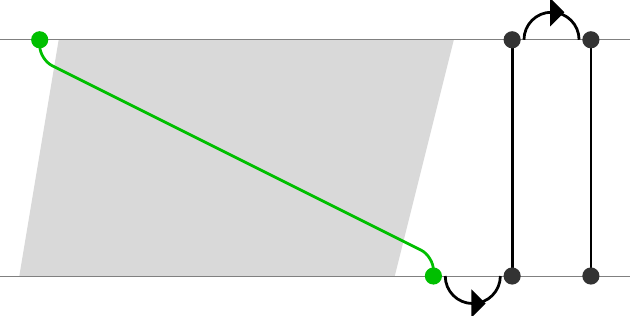}
\end{tabular} 
\end{center}

\begin{center}
\begin{tabular}[t]{cc}
$q_1^{green}S'_1(id_s)$ & $T^{i_{\ell_1}}q_1^{green}S'_1(id_s)$\\
\adjustimage{valign=T}{figure/big_figure5.pdf}  &\adjustimage{valign=T}{figure/big_figure6.pdf}  
\end{tabular}
\end{center}

The last case is somewhat more difficult, the permutations have the form: \[\tau=*,{\color{blue}k-1},*,{\color{green}k},n-1,n \text{ thus } q_2^{brown}(\tau)=(*+1),{\color{brown}1},{\color{blue}k},(*+1),{\color{green}k+1},n,n+1.\]
\begin{center}
\begin{tabular}{cc}
\raisebox{4pt}{\includegraphics[scale=.5]{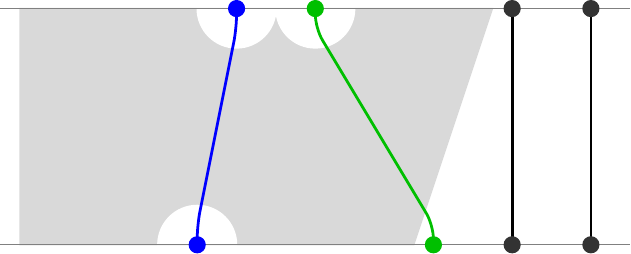}} &\includegraphics[scale=.5]{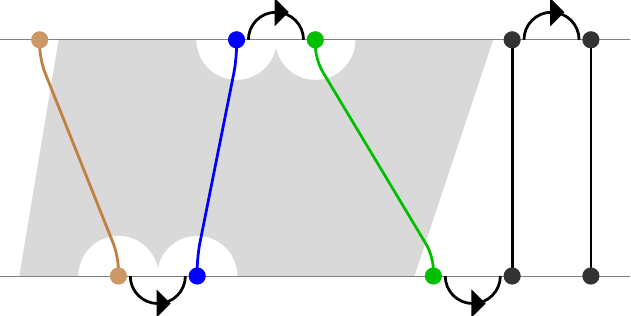}
\end{tabular} 
\end{center}

\begin{center}
\begin{tabular}[t]{cc}
$q_2^{brown}T^{i_{\ell_1}}q_1^{green}S'_1(id_s)$ & $q_2^{brown}T^{i_{\ell_1}+1}q_2^{brown}T^{i_{\ell_1}}q_1^{green}S'_1(id_s)$\\
\adjustimage{valign=T}{figure/big_figure7.pdf}  &\adjustimage{valign=T}{figure/big_figure8.pdf}  
\end{tabular}
\end{center}

\begin{center}
\begin{tabular}[t]{l}
$T_r^{i_{\ell_2}\!-\!1}q^{brown}_2   \ldots T_r^{i_{\ell_1}\!+\!1}q^{brown}_2T_r^{i_{\ell_1}}q_1^{green}S'_1(id_s)$\\
\adjustimage{valign=T}{figure/big_figure9.pdf}  \\[3mm]
\rule{0mm}{5mm}$q_1^{green}T_r^{i_{\ell_2}\!-\!1}q^{brown}_2   \ldots T_r^{i_{\ell_1}\!+\!1}q^{brown}_2T_r^{i_{\ell_1}}q_1^{green}S'_1(id_s)$\\
\adjustimage{valign=T}{figure/big_figure10.pdf}  
\end{tabular}
\end{center}
The two figures above represent the almost the same permutation (we simply applied $q_1^{green}$ to the second) but we compacted the second for space-saving purpose.

\begin{center}
\begin{tabular}[t]{c}
 \multicolumn{1}{c}{$q_2^{brown}T^{i_{\ell_1}}q_1^{green}S'_1(id_s)$}\\
\adjustimage{valign=T}{figure/big_figure11.pdf}
\end{tabular}
\end{center}
In order to represent the full sequence $S'(id_s)$, each sequence $T_r^{i_{\ell_{p+1}}\!-\!1}q^{brown}_2   \ldots T_r^{i_{\ell_p}\!+\!1}q^{brown}_2T_r^{i_{\ell_p}}$ for $1\leq p\leq r$ correspond to a pair of blocks $A_p,B_p$. Then we have: 

\begin{center}
\begin{tabular}[t]{c}
$S'(id_s)$\\
\adjustimage{valign=T}{figure/big_figure12.pdf}  
\end{tabular}
\end{center}
Clearly the blue, green, red and brown subsequences are all decreasing and the black subsequence has size at most 6.
\qed
\end{pf}

Note that the permutations $(\s,c)$ constructed are standard and the permutations $(\s',c')$ are at alternation distance 1 of a standard (since $\s'(n)=n$) thus the lower bound also works for the set of standard permutations of a class.

\section{Algorithm and upper bound on the diameter}

Let us first recall the notion of alternation distance:
\begin{definition}[distance and alternation distance]
Let $\s$, $\s'$ configurations in
the same class. The \emph{graph distance} $d_G(\s,\s')$ is the
ordinary graph distance in the associated Cayley Graph,
i.e.\ $d_G(\s,\s')$ is the minimum $\ell \in \bN$ such that there
exists a word $w \in \{L,L^{-1},R,R^{-1}\}^*$,
of length $\ell$, such that $\s' = w \s$.  We also define the
\emph{alternation distance} $d^{A}(\s,\tau)$ as the analogous quantity, for words
in the infinite alphabet $\{(L)^j,(R)^j\}_{j \geq 1}$.
\end{definition}

In this section we will prove that:

\begin{theorem}\label{thm_diameter_up}
Let $C$ be a Rauzy class of size $n$, the diameter of $C$ for the alternation distance is at most $27n$.
\end{theorem}

and

\begin{theorem}\label{thm_algo}
Let $\s$ and $\s'$ be two permutations of size $n$ we can decide if they are in the same class and (if positive) find a word $w \in \{L,L^{-1},R,R^{-1}\}^*$ such that $\s'=w\s$ in time $O(n^2)$. 
\end{theorem}

More precisely we will construct the word $w$ of theorem \ref{thm_algo} and show that it is of size at most $26n$ for the alternation distance thus deducing theorem \ref{thm_diameter_up}

Let us recall the notion of zig-zig path:

\begin{figure}[b!!]
\begin{center}
\includegraphics{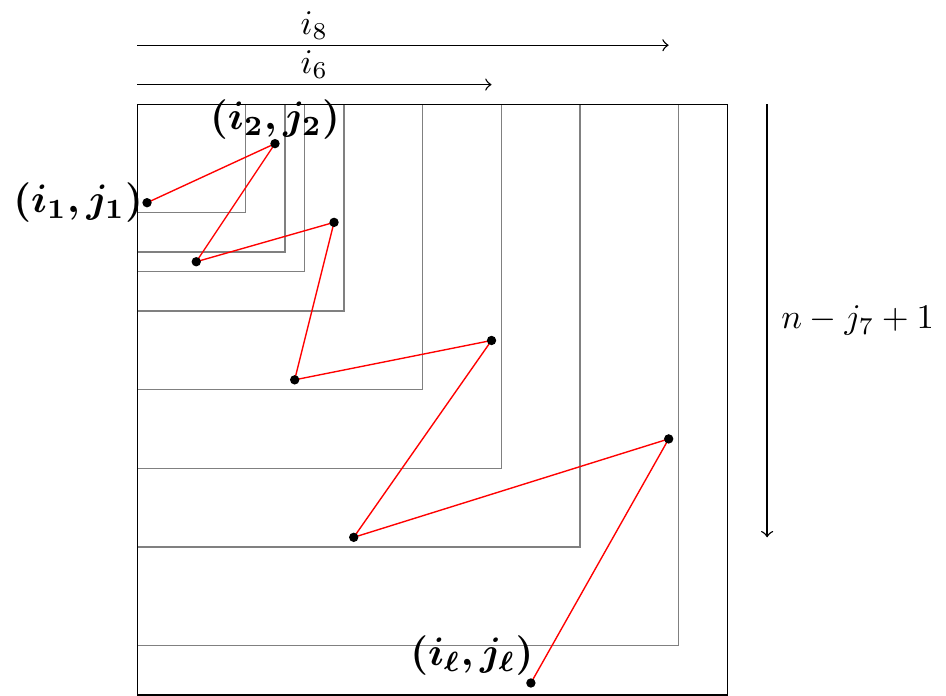}
\caption{\label{figure.zigzag}%
  Structure of a zig-zag path, in matrix representation.  It is
  constructed as follows: for each pair $(i_a,j_a)$, put a bullet at
  the corresponding position, and draw the top-left square of side
  $\max(i_a,j_a)$. Draw the segments between bullets $(i_a,j_a)$ and
  $(i_{a+1},j_{a+1})$ (here in red). The resulting path must connect
  the top-left boundary of the matrix to the bottom-right boundary.
  If the bullets are entries of $\s$, then a top-left $k \times k$
  square cannot be a block of the matrix-representation of $\s$,
  because, in light of the inequalities (\ref{eq.6537564}), one of the
  two neighbouring rectangular blocks (the $k \times (n-k)$ block to
  the right, or the $(n-k) \times k$ block below) must contain a
  bullet of the path, and thus be non-empty.}
\end{center}
\end{figure}

\begin{definition}
A set of edges
$\big((i_1,j_1),(i_2,j_2),\ldots,(i_{\ell},j_{\ell})\big)$ is a
\emph{$L$ zig-zag path} if $i_1=1$, the indices satisfy the pattern of
inequalities
\begin{align}
j_{2b} &> j_{2b-1}
\ef,
&
i_{2b} &> n-j_{2b-1}+1
\ef,
&
i_{2b} &> i_{2b+1}
\ef,
&
n-j_{2b+1}+1 >i_{2b}
\ef,
\label{eq.6537564}
\end{align}
and either $i_\ell=n$ or $j_\ell=1$.  The analogous structure
starting with $j_1=n$ is called
\emph{$R$ zig-zag path}.

Given permutation $\s$, we define $Z(\s)$ to be the size of the smallest Zig-Zag path (i.e. $\ell$ is minimized). Clearly $Z(\s)\leq n$.
\end{definition}

We have shown (see \cite{DS17} section 3.2) that if $\s$ has a Zig-zig path of length $\ell$ then there exists a word $w\in\{(L)^j,(R)^j\}_{j \geq 1}$ of size (for the alternation metric) $|w|\leq \ell$ such that $\s'=w\s$ is standard. Moreover we also proved that we could compute a greedy zig-zig path of a permutation in linear time such that its length was at most $Z(\s) +1$.

Thus we can always start the algorithm from a standard permutation. Note that computing $\s'=w\s$ by applying the operators $L^i$, $R^i$ of $w$ takes $O(|w|n)$ where $|\cdot|$ is the alternation length (which is fine since it is at most quadratic). We can actually standardize a permutation in linear time but the algorithm is more complicated.

Let us first start with a few preliminaries.

Define $T_{min}(\s)= T_i(\s)$ with $i$ such that $\s(i)=1$.

\begin{lemma}\label{lem_T_and_boosted}
Let $\tau$ be a permutation and $\s=T_i(\tau)$ for some $i$, let $c$ be the $(n-2,2)$ coloring of $\s$ with the edges $(i+1,1)$ and $(i+2,\s(i+2))$ grayed so that $\tau$ is the reduction of $(\s,c)$ and let $S$ be a sequence such that $\tau'=S(\tau)$ then $(\s',c')=B(S)(\s)$ verifies $\s'=T_{i'}(\tau')$ for some $i'$.

Moreover if $\s=T_{min}(\tau)$ then $(\s',c')=B(S)(\s)$ verifies $\s'=T_{min}(\tau')$.
\end{lemma}

This is a consequence of lemma 3.12 of \cite{DS17} since the $T$ structure is a particular square constructor. The moreover part is not in the lemma 3.12 but can be deduced from the proof.

\begin{lemma}\label{lem_T_seq}
Let $\tau$ be a standard permutation and define $\s=T_i(\tau)$ and $\s'=T_j(\tau)$, then there exists a sequence $S$ of alternation length at most 5 such that $\s'=S\s$
\end{lemma}
This is a consequence of lemma 3.13 of \cite{DS17} since the $T$ structure is a particular square constructor.

\begin{proposition}\label{pro_S_O}
Let $\tau$ and $\tau'$ two connected permutations and $S$ such that $\tau'=S(tau)$. Let $S_O$ be any sequence of operators $q_1,q_2$ and $T_{min}$, i.e.:
\[S_O=T_{min}^{i_k}q_{j_k}\ldots T_{min}^{i_1}q_{j_1}T_{min}^{i_0} \text{ with } \forall 1\leq \ell \leq k-1, \ i_\ell\geq 1,\ i_0,i_k\geq 0 \text{ and  } \forall \ell, \ j_\ell \in \{1,2\}.\]
and define $\s=S_O(\tau)$, as well as a colorings $c$ such that $\tau$ is the reduction of $(\s,c)$ then $(\s',c')=B(S)(\s,c)$ where $B(S)$ is the boosted sequence of $S$ verifies that $\s'=S_O(\tau')$
\end{proposition}

\begin{pf}
By induction on the number of operators in $S_O$. The inductive step is trivial once the initial step is proven so we only need to check for $T_{min}$, $q_1$ and $q_2$.
For $T_{min}$ this is a consequence of lemma \ref{lem_T_and_boosted}, and for $q_1$ and $q_2$ this is a direct consequence of the proof of lemma 5.10  in \cite{DS17}.
\end{pf}

The following lemma correspond to proposition 40 of \cite{D18}.
\begin{lemma}\label{lem_correct_cycle}
Let $\s$ be a standard permutation with cycle invariant $(\lam,r)$ and type $X(r,i)$. Let $c$ be the coloring of $\s$ with the first edge $(1,1)$ grayed and let $\tau$ be the reduction of $(\s,c)$. Then $\tau$ have the cycle invariant $(\lam'=\lam \setminus \{i\},r'=r+i-1)$.
\end{lemma}

Let us define \[d(\s)=\s(1)-1,\ldots,\widehat{\s(\s^{-1}(1))},\ldots,\s(n)-1.\]
Thus in the above lemma we have $\tau=d(\s).$

\begin{lemma}\label{lem_only_one_red}
Let $\s$ be a standard permutation of rank 2, then there is at most one permutation $\s'$ of type $X$ in the standard familly such that $d(\s')$ is reducible.
\end{lemma}

\begin{pf}
If there are no permutation of type $X$ such that their image by d is reducible in the standard family then we are done. If there is let $\s'$ be one, we shall prove that it is unique. 

$\s'$ must have the form described below. Suppose $d(L^i(\s'))$ is reducible with $i\leq k_1$ then $L^i(\s')$ must have the form described below. Likewise if $d(L^j(\s'))$ is reducible with $j> k_1$,  $L^j(\s')$ must have the form described below.
\begin{center}
\begin{tabular}{ccc}
\includegraphics[scale=1]{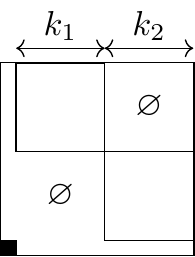} & \includegraphics[scale=1]{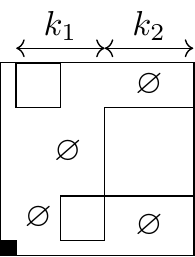} & \includegraphics[scale=1]{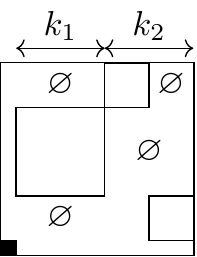}\\
$\s'$ & $L^i(\s')$ & $L^j(\s')$
\end{tabular}
\end{center}

Thus $\s'$ must be of this type: \begin{tabular}{c}
\includegraphics[scale=1]{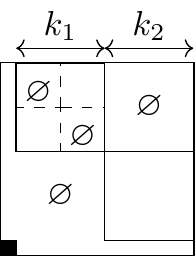}\end{tabular} or  \begin{tabular}{c}
 \includegraphics[scale=1]{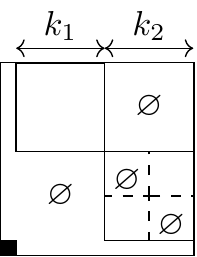}
\end{tabular}\rule{2.03cm}{0pt}(A)

However since $\s'$ has rank 2, it has the form
\begin{center}\begin{tabular}{c}
\includegraphics[scale=.5]{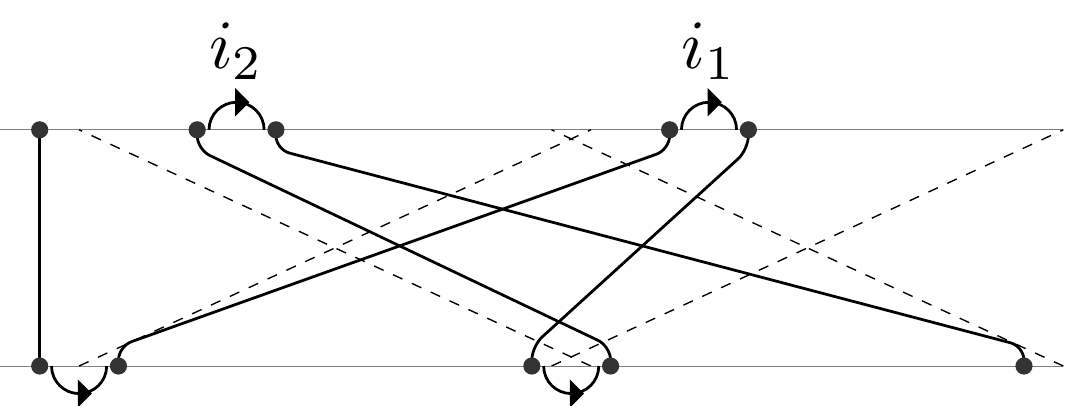}\end{tabular} or equivalently \begin{tabular}{c}
 \includegraphics[scale=1]{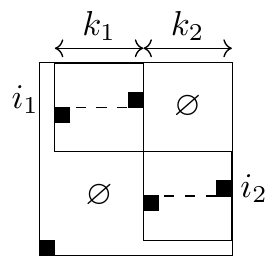}
\end{tabular}\rule{2mm}{0pt}(B) \end{center}

Clearly considering the two requirements (A) and (B) the two blocks can only be separated at $i_1$ and $i_2$ respectively, but $L^{i_1}(\s')$ and $L^{i_2}(\s')$ are of type $H$:
\begin{center} \begin{tabular}{cc}
\includegraphics[scale=.5]{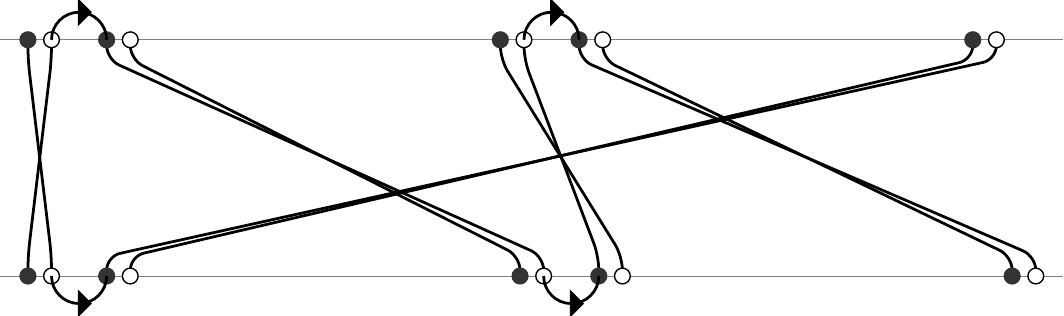} &
 \includegraphics[scale=.5]{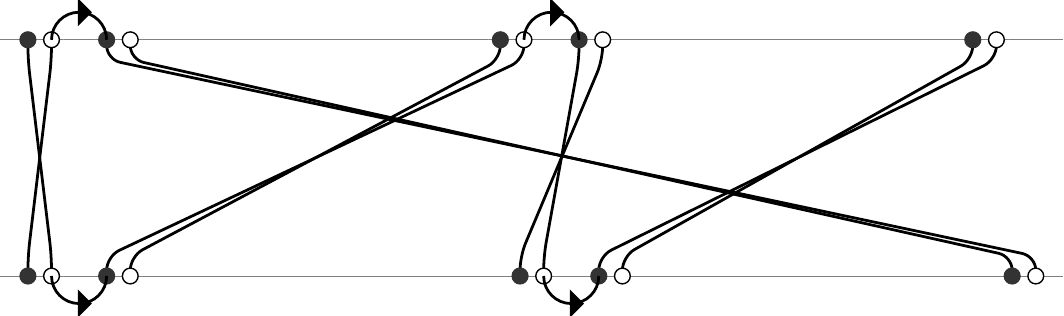}\\
$L^{i_1}(\s')$ type $H(1,2)$ & $L^{i_2}(\s')$ type $H(2,1)$ 
\end{tabular}\end{center}
\qed
\end{pf}

\begin{lemma}\label{lem_Z_rank_2}
Let $\s$ be a standard permutation of rank 2, then for any (but at most one) permutation $\s'$ of type $X$ in the standard familly $\tau'=d(s')$ verifies $Z(\tau')\leq 5$.
\end{lemma}

\begin{pf}
The 'but at most one' concerns the only $\tau'$ that could be reducible. For any other $\tau'$, by the previous proposition it has either of the following forms:
\begin{center}
\begin{tabular}{c}
\includegraphics[scale=1]{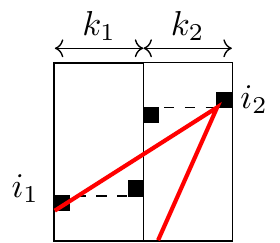}\end{tabular} or \begin{tabular}{c} \includegraphics[scale=1]{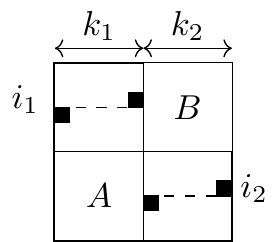} \end{tabular}
 \end{center}

In the first case $Z(\tau)\leq 2$, in the seconde case both $A$ and $B$ are non-empty (at least one must be non-empty if not we have a reducible permutation and if one is then so must be the other since $\tau'$ is a permutation and thus has exactly one point per line and column). Also note that in the figure $i_2$ could be big enough that the two black points are in $B$ but this just makes things easier.

We must have : 
\begin{center}
\begin{tabular}{c}
\includegraphics[scale=1]{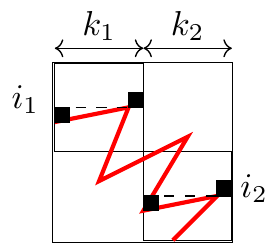}\end{tabular} 
 \end{center}

Thus $Z(\tau')\leq 6$
\qed
\end{pf}

We can now start our proof. Unfortunatly the analysis of time of this algorithm is somewhat difficult as presented in the recursive form of the proof. Thus we need a second proposition where the structure is iterative.

\begin{theorem}\label{thm_upper_bound}
Let $\s$ and $\s'$ be two standard permutations of size n, if there are in the same class then we can find in $O(n^2)$ a word $w \in \{L,L^{-1},R,R^{-1}\}^*$ such that $\s'=w\s$ and $w$ has alternation length at most $26n$. Otherwise we certify that they are not in the same class.
\end{theorem}

\begin{pf}
First of all we solve the case of exceptional classes. By the structure theorems of appendices C.1 and C.2 we can decide in linear time if two permutations are in exceptional classes and the diameter of those classes is at most $n$. 
Thus we can suppose that neither $\s$ and $\s'$ are in an exceptional class.

The proof is done by induction. We first compute the cycle invariant $(\lam,r)$ of $\s$ and $\s'$, if there are different we stop, if not let $r$ be this rank.\\

\noindent \textbullet$\ $ If $r>2$, by theorem 5.5 of \cite{DS17} there exists a sequence $S_1$ and $S'_1$ of alternation length at most 6 such that $\s_1=S_1\s$ and $\s'_1=S'_1\s'$ have the following property: 
$\s_1=T(\tau_1)$ and $\s'_1=T(\tau'_1)$, $\tau_1$ and $\tau'_1$ have size $n-2$ and have $Z(\tau)\leq5$ and  $Z(\tau')\leq5$.  

Unfortunatly $\tau_1$ and $\tau'_1$ can be in the exceptional class $id$ or $id'$ (depending on the invariant $(\lam,r)$). If both are then this is fine since we know the cayley graph of those classes we can directly output a path from $\tau_1$ to $\tau_1'$. Thus the only problem is if only one is, let us say $\tau_1$. By lemma \ref{lem_exception_class} below we can find a sequence $S_{excp}$ of alternation length at most 19 such that $S_{excp}((\s_1))=T(\tau_{ex})$ with $Z(\tau_{ex})\leq 4$ and $\tau_{ex}$ is not in an exceptional class.

Let us rename $\tau_{ex}$ with $\tau_1$. 

Thus there exists two sequences $S_2, S_2'$  of alternation length at most 5 (or 2 or 4 in the exceptional cases) such that $\tau_2=S_2\tau_1$ and $\tau_2'=S'_2\tau'_1$ are standard. 

By induction hypothesis, if $\tau_2$ and $\tau_2'$ are not in the same class (thus they have different invariant) neither are $\s_1=T(\tau_1)$ and $\s'_1=T(\tau'_1)$ since $T$ increase the rank by 2 and leave invariant both cycle and the sign invariant (thus $\s_1$ and $\s'_1$ must also have different invariant).

If $\tau_2$ and $\tau_2'$ are in the same class, there exist $S_3$ of alternation length at most $14(n-2)$  such that $\tau'_2=S_3\tau_2$.\\

 Let $c_1$ (respectively $c'_1$) be the coloring of $\s_1$ (respectively $\s'_1$) such that the reduction is $\tau_1$ (respectively $\tau'_1$). 

by lemma \ref{lem_T_and_boosted} $B(S_2)(\s_1)=T_i(\tau_2)$  and $B(S'_2)(\s'_1)=T_i(\tau'_2)$ for some $i,i'$. Then let $S_{T,1}$ and $S'_{T,1}$ such $S_{T,1}B(S_2)(\s_1)=T_{min}(\tau_2)$ $S'_{T,1}B(S'_2)(\s'_1)=T_{min}(\tau'_2)$,  by lemma \ref{lem_T_seq} the sequence $S_{T,1}$ and $S'_{T,1}$  exists and have alternation length at most 5. Finally define $(\s_2,c_2)=S_{T,1}B(S_2)(\s_1,c_1)$ and $(\s_2',c_2')=S'_{T,1}S_2''B(S'_2)(\s'_1,c'_1)$), then by lemma \ref{lem_T_and_boosted}  we have $B(S_3)(\s_2,c_2)=(\s_2',c_2')$.

To sum-up we have the following diagram:
\begin{center}
$\begin{tikzcd}
\s \arrow{d}{S_1} \arrow[bend right=50,swap]{dd}{S_1}  \arrow{rrr}{S'^{-1}_1B(S'_2)^{-1}S'^{-1}_{T,1}B(S_3)S_{T,1}B(S_2)S_1} &&&\s' \arrow{d}{S'_1} \arrow[bend left=50]{dd}{S_1'}  \\
T(\tau_1) \arrow{r}{S_{T,1}B(S_2)} \arrow[xshift=0.6ex,dash]{d}  \arrow[xshift=-0.6ex,dash]{d}
&T_{min}(\tau_2)\arrow{r}{B(S_3)} \arrow[xshift=0.6ex,dash]{d}  \arrow[xshift=-0.6ex,dash]{d} 
&T_{min}(\tau'_2) \arrow[xshift=0.6ex,dash]{d}  \arrow[xshift=-0.6ex,dash]{d}
& T(\tau'_1) \arrow[swap]{l}{B(S_2')} \arrow[xshift=0.6ex,dash]{d}  \arrow[xshift=-0.6ex,dash]{d} \\
 (\s_1,c_1) \arrow{r}{S_{T,1}B(S_2)} \arrow{d}{red} 
& (\s_2,c_2) \arrow{r}{B(S_3)}  \arrow{d}{red}
 & (\s'_2,c'_2)  \arrow{d}{red} 
& (\s'_1,c'_1) \arrow[swap]{l}{S'_{T,1}B(S_2')} \arrow{d}{red}\\
\tau_1 \arrow{r}{S_2} 
&\tau_2\arrow{rr}{S_3} 
& \tau_2'  
& \tau'_1 \arrow[swap]{l}{S_2'}
\end{tikzcd}$
\end{center}

and thus $\s'=S'^{-1}_1B(S'_2)^{-1}S'^{-1}_{T,1}B(S_3)S_{T,1}B(S_2)S_1(\s)$.

Moreover since the boosted dynamics does not change the alternation length (since we only increase the exponent of each operator $L,R$ of the sequence) we have \[|S'^{-1}_1B(S'_2)^{-1}S'^{-1}_{T,1}B(S_3)S_{T,1}B(S_2)S_1|=6+5+5+26(n-2)+5+5+6+\underbrace{19}_{\text{if $\tau_1$ in exceptional class} }\leq26n\]

\noindent \textbullet$\ $ If $r=1$ Let $\s_1=L^k(\s)$ and $\s'_1=L^{k'}(\s')$ for $k,k'$ such that they are both of type $X(r,i)$ and $\s_1^{-1}(2)<\s_1^{-1}(n)$ and $\s_1'^{-1}(2)<\s_1'^{-1}(n)$ (those exist by proposition \ref{trivialbutimp}) , then by lemma \ref{lem_correct_cycle}, $\tau_1=d(\s_1)$ and $\tau_1'=d(\s'_1)$ have same cycle invariant $(\lam\setminus \{i\},i)$ and their sign is egal to that of $\s_1$ and $\s_1'$ respectively by lemma \ref{lem_sign_control}. Thus they are connected if and only if $\s_1$ and $\s_1'$ are.

 Let $\s_2=R(\s_1)$ and $\s'_2=R(\s_2)$ then $\s_2=q_1(\tau_1)$ and $\s'_2=q_1(\tau'_1)$.
As a reminder (cf definition 5.9), standard permutations are not in the image of $q_i$ (by choice of definition) but are at distance one of such a permutation: If $\s$ is a standard permutation and has rank $1\leq i\leq 2$ let $\tau=d(\s)$ then $R(\s)=q_i(\tau)$.

Let $\tau_2=R^{\s_1^{-1}(2)-2}\tau_1$ and $\tau_2'=R^{\s_1'^{-1}(2)-2}\tau_1'$ then there are standard. 

By induction if they are not connected we are done. If not let $S$ be a sequence such that $\tau_2'=S(\tau_2)$.

Let $c_1$ (respectively $c'_1$) be the coloring of $\s_1$ (respectively $\s'_1$) such that the reduction is $\tau_1$ (respectively $\tau'_1$). Then by lemma \ref{} we have $B(R(^{\s_1'^{-1}(2)-2})(\s'_1,c'_1)=B(S)B(R^{\s_1^{-1}(2)-2})(\s_1,c_1)$.

Moreover since $S$ has alternation length $26(n-1)$, the sequence connecting $\s$ to $\s'$ has length $26(n-1)+4\leq 26n$.

\noindent \textbullet$\ $ If $r=2$ the case is significantly more complicated than $r=1$. Since the rank is even, $\lambda$ must contain at least one even cycle of length $2\ell$ and $s=0$ by the classification theorem \ref{thm.Main_theorem_2}. (Also by the classification theorem we know that they are connected since there is only one class of a given cycle invariant with sign 0)

We distinguish two cases: either there are exactly one even cycle in $\lam$ or there are at least two cycles of length $2\ell$ and $\ell'$.
\paragraph*{If there are at least two cycles} then the argument is similar to the case $r=1$:

Let $\s_1=L^k(\s)$ and $\s'_1=L^{k'}(\s')$ for $k,k'$ such that they are both of type $X(r,\ell')$ and $\s_1^{-1}(2)<\s_1^{-1}(n)$ and $\s_1'^{-1}(2)<\s_1'^{-1}(n)$ (those exist by proposition \ref{trivialbutimp}) , then by lemma \ref{lem_correct_cycle}, $\tau_1=d(\s_1)$ and $\tau_1'=d(\s'_1)$ have cycle invariant $(\lam\setminus \{\ell'\},r+\ell'-1)$ which contains the even cycle $2\ell$ and their sign must then be egal to 0. Thus they are connected by the classification theorem.

The rest is identical to the case $r=1$.

\paragraph*{If there is exactly one even cycle $2\ell$} we have by the dimension formula $2\ell=n-3$ since $2\ell+r=n-1$ and $r=2$.
The problem here is that if we take $d(\s)$ for a standard permutation $\s$ of type $X(r,2\ell)$ then $d(\s)$ will have invariant $(\varnothing,2\ell+1,s)$ with $s=\pm1$. (we struggled with the same problem in the first case of proposition \ref{pro_structure_2}) 

Thus defining $\tau_1=d(\s)$ and $\tau_1'=d(\s')$ as above we may have that $\tau_1$ and $\tau_1'$ have opposite sign invariant and thus are not connected.

We address this problem thanks to lemma 73 of \cite{D18} which say that in this case if two permutations $\s_1$ and $\s_2$ of type $X(r,2\ell)$ a standard familly are consecutive (in a sense defined in lemma 73) then the sign invariant of $d(\s_1)$ and $d(\s_2)$ are opposed. 

By lemma \ref{lem_only_one_red} we know that at most one permutation of type $X(r,2\ell)$ can have its image by $d(\cdot)$ reducible moreover we also know from appendix A of \cite{D18} that at most one permutation of type $X(r,2\ell)$ can have its image by $d(\cdot)$ be in the exceptional class $id$ thus we might have to discard 2 permutations of type $X(r,2\ell)$ but since $2\ell=n-3$ by the dimension formula that does not matter.

Thus we choose $\s_1$ and $\s'_1$ of type $X(r,2\ell')$, and such that $\tau_1=d(\s_1)$ and $\tau_1'=d(\s_1)$ have the same sign invariant (which is possible by lemma 73 of \cite{D18}) then by lemma \ref{lem_correct_cycle} they have the same invariant $(\varnothing,n-2)$.
We define $\s_2=R(\s_1)$ and $\s'_2=R(\s_2)$ then $\s_2=q_2(\tau_1)$ and $\s'_2=q_2(\tau'_1)$ and we let $\tau_2=S_1(\tau_1)$ and $\tau_2'=S'_1(\tau_1)$ be the standardization of $\tau_1$ and $\tau'_1$, both sequences have alternation length at most 6 by lemma \ref{lem_Z_rank_2}.

The rest is identical to the case $r=1$.

There is an important note here: algorithmically computing the sign is not a trivial task (though it can be done in polynomial time and this proposition in particular provides a way to do it in quadratic time) so given $\tau_1$ and $\tau_1'$ we do not know whether that have the same sign or not. However this is not a problem, it suffices to take $\tau_1^s$ and $\tau_1^{-s}$ the reduction of two consecutive permutations $\s_1$ and $\s_2$ that we know will have different sign, and continue the induction with both. Then for one the algorithm will find that $\tau'_1$ and say $\tau_1^{s}$ are in the same class and output a path and for the other $\tau'_1$ and $\tau_1^{-s}$ the algorithm will find that they are not connected. This introduce a fork but since this happens only if the invariant is $(\{2\ell\},2)$ this can only occur once in the whole induction (thus the algorithm still runs in quadratic time since this only double the total time taken at most).\\

For the algorithm runtime analysis, we do not use recursivity but we transform the procedure so it is iterative: instead of working from $\s$ to $\tau$, we color in gray the two edges of the $T_{min}$ or the single edge of the $q_i$ of $\s$. Clearly by proposition \ref{pro_S_O} at each step ($T$, $q_1$ or $q_2$) we get two permutations $(\s_i,c_i)$ and $(\s'_i,c'_i)$ and two reduced permutations $\tau_i$ and $\tau'_i$ such that $\s_i=S_{O,i}(\tau_i)$ and  $\s_i'=S_{O,i}(\tau'_i)$, moreover $S_{O,i}=T_{min}S_{O,i-1}$ or  $S_{O,i}=q_jS_{O,i-1}$  depending on whether the step is a $T$ or a $q_j$.

We do every step until we reach constant size and we can find directly a sequence $S$ connecting the current $\tau$ and $\tau'.$. Then by proposition \ref{pro_S_O}  $\s=B(S)\s'$ since $\s=S_O(\tau)$ and $\s'=S_O(\tau')$ for some $S_O$.

Each step takes linear time (the sequences we use are constant in alternation size) and we have a linear number of steps so the algorithm is quadratic. 

\qed
\end{pf}

\begin{remark}
We can actually show that the case of $\tau_i$ ending in an exceptional class can only happens three time (twice for $id'$ and once for $id$). Indeed for $id'$ this is easy: $\tau_i$ can be in $id'$ only if its cycle invariant is $(\{n-2\},1)$ or $(\{(n-2)/2,(n-2)/2\},1)$ which can only happens once for each in the recursion. 

The case $id$ is much more tedious (and will not be done here). One can show that in the recursion starting from $\s_{ex,n}$ the case where a $\tau$ at some point is in the exceptional class does not happen. This is possible to prove such statement since we already know the structure of  $\s_{ex,n}$ and thus we can analyse the recursion at each step and verify that each $\tau$ is not in the exceptional class. Thus the case where $\tau_i$ is in $id$ can only happen once since afterward we continue the recursion with $\s_{ex,k}$ for some $k$.

All this put together we can replace in the proof of proposition \ref{} the 20 by simply a 5 which lead to the better upper bound for the diameter of $16n +c$. Another small modification of the procedure yields $14n+c$.
\end{remark}

Combining theorem \ref{thm_upper_bound} and the first standardization sequence whose alternation length is at most $n$, we prove theorem \ref{thm_diameter_up}.
\begin{lemma}\label{lem_exception_class}
Let $\s$ and $\tau$ such that: $\s=T(\tau), Z(\tau)\leq 6$ and $\tau$ is in an exceptional class. Then there exists a sequence $S_{excep}$ of alternation length at most 19 such that $\s'=S_*(\s)=T(\tau_{ex}), Z(\tau_{ex})=2$ or 4 (for $id'$) and $\tau_{ex}$ is not in an exceptional class.
\end{lemma}
\begin{pf}
Let us begin by the id case. We have $\tau$ with $Z(\tau)\leq 6$, $\s=T(\tau)$ and $\tau$ is in $id$. Let $S_0$ be the sequence of alternation length 6 such that $id=S_0(\tau)$, then $B(S_0)(\s)=T_i(id)$ for some $i$, then there exists $S_1$ of alternation length at most  3 such that $\s'=S_1B(S_0)(\s)=T(id)=T_{min}(id)$, Then define $S_3$ as follows:
\begin{center}
$\begin{tikzcd}[column sep=125pt]
\raisebox{-35pt}{\includegraphics[scale=.5]{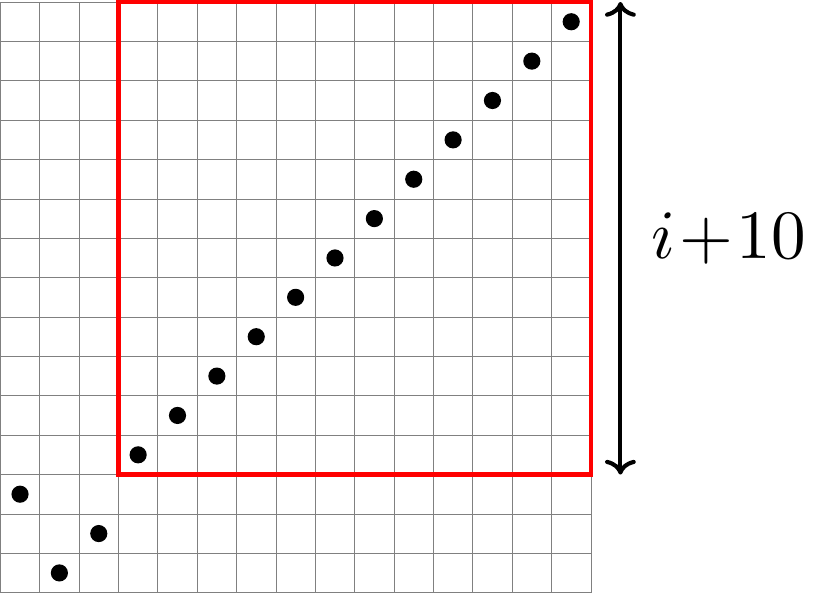}}\arrow{r}{S_3=R^8L^2R^{i\!+\!2}L^{i\!+\!1}R^7L^9R^{i\!+\!3}L^{i\!+\!3}R^2L^8}& \raisebox{-35pt}{\includegraphics[scale=.5]{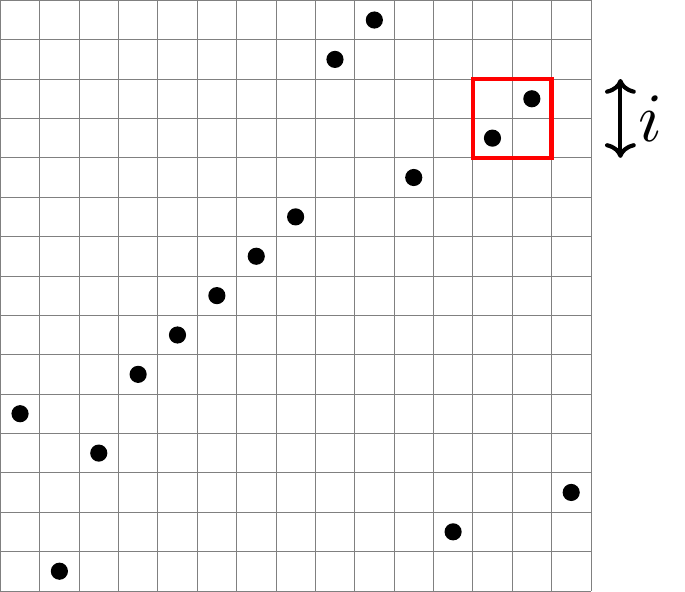}}\\[-6mm]
\!\!\!\!\!\!\!\!\!\!\!\!\!\!\!\!\!\!\!\!T(id) &\!\!\! T(S_{ex,i})
\end{tikzcd}$
\end{center}
Let $\s''=S_3\s'=T(\tau_{ex})$, clearly $\tau_{ex}=S_{ex,i}$ is not in $id$ (it suffices to standardize and check that it is not $id$) and $Z(\tau')=2$. 
The final sequence has alternation length $6+3+10=19$.

The case for $id'$ is similar with 
\begin{center}
$\begin{tikzcd}[column sep=125pt]
\raisebox{-35pt}{\includegraphics[scale=.5]{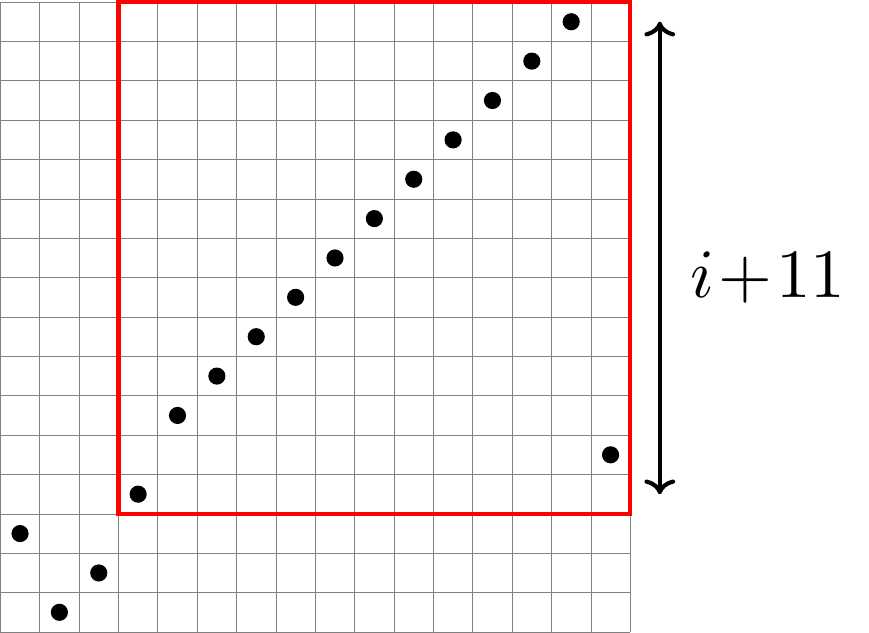}}\arrow{r}{S_3=R^{i\!+\!2}L^3R^9L^3R^{i\!+\!5}L^4RL^{i\!+\!7}R}& \raisebox{-35pt}{\includegraphics[scale=.5]{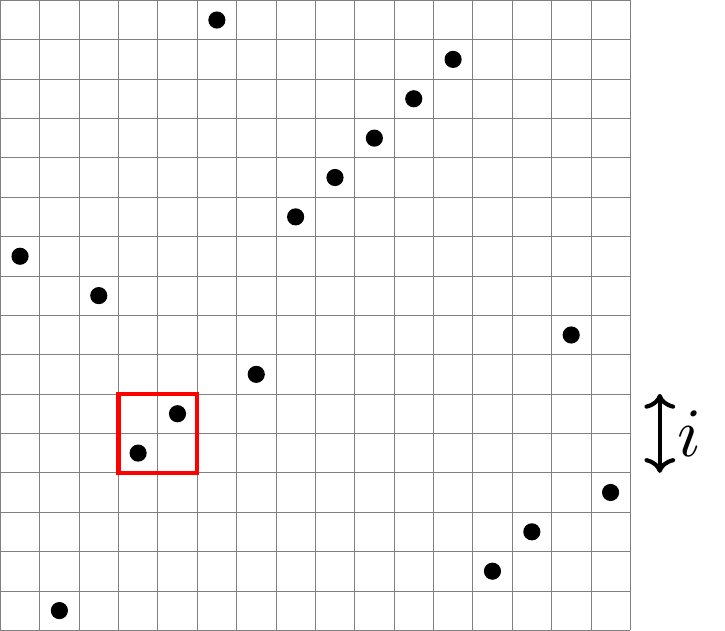}}
\end{tikzcd}$
\end{center}
and $Z(\tau')\leq 4$.
\end{pf}

\end{document}